\documentclass[12pt]{article}
\usepackage{amsmath, amssymb}
\begin{document}
\centerline{\bf Effective Termination of Kohn's Algorithm}
\centerline{\bf for Subelliptic Multipliers}

\bigskip
\centerline{Yum-Tong Siu\ %
\footnote{Partially supported by Grant DMS-0500964 from the
National Science Foundation.} }

\bigbreak\centerline{\it Dedicated to Professor Joseph J. Kohn on his
75th Birthday}

\bigbreak\noindent{\bf Introduction.} In this note we will discuss
the problem of the effective termination of Kohn's algorithm for
subelliptic multipliers for bounded smooth weakly pseudoconvex
domains of finite type [Ko79].  We will give a complete proof for
the case of special domains and will only indicate briefly how this
method is to be extended to the case of general bounded smooth
weakly pseudoconvex domains of finite type. The method is rather
simple and uses some local theory of algebraic geometry. People with
some minimal background in algebraic geometry may find the
algebraic-geometric techniques involved in this note very simple or
even completely trivial. Since this topic is of interest mainly to
the analysts I will use as much as possible the language of function
theory to describe our method.

\bigbreak In a number of conferences in recent years I gave talks on
this topic, but because of time limitation never had the opportunity
to present all the details.  This note is written to make the details
available. This note will appear in the special issue of {\it Pure and Applied Mathematics Quarterly} for Professor Joseph J. Kohn.

\bigbreak The termination of Kohn's algorithm in the real-analytic
case was verified by Diederich-Fornaess [DF78] without
effectiveness. In this note we are also going to formulate Kohn's
algorithm geometrically in terms of the theorem of Frobenius on
integral submanifolds and present a proof from this geometric
viewpoint so that one can see clearly how the procedures of Kohn's
algorithm arise naturally in the geometric context and why the
real-analyticity facilitates the proof of the termination of Kohn's
algorithm. We present this more geometric proof here to provide an
alternative to the proof of Proposition 3 on pp.380-388 of [DF78] 
which is the key step of [DF78] and
which is still quite a bit of a challenge to follow.  Moreover, the
proof of the real-analytic case of the ineffective termination of
Kohn's algorithm from the geometric viewpoint gives a better
understanding of the r\^ole played by the real-analytic assumption
and of what stands as a hurdle between generalizing the ineffective
real-analytic case to the ineffective smooth case.

\medbreak We also hope that this geometric viewpoint will provide an
easier and more transparent setting for further developments of the
application of algebraic-geometric techniques to general partial
differential equations which Kohn has pioneered with his algorithm
for the complex Neumann problem [Ko79].  The key point of Kohn's
theory is the following.  If the distribution of jets, where the
partial-differential-equation estimate fails to hold,
is not integrable even over unreduced points (or
Artinian subschemes) of arbitrarily high order, then there is an
algorithm to apply algebraic-geometric techniques to derive the
partial-differential-equation estimate.  Kohn implemented his theory
for the complex Neumann problem.  His theory should be applicable to
systems of partial differential equations in a more general setting.
Such an expected further development of his theory remains yet to be
carried out.

\bigbreak In a private communication Kohn told me that he has a
direct proof of the result of Diederich-Fornaess [DF78] on the
ineffective termination of the Kohn algorithm for the real-analytic
case by using power series expansion of the real-analytic defining
function and explicitly keeping track of various partial
differentiations in the holomorphic and anti-holomorphic directions.

\bigbreak Andreea Nicoara recently posted an article [Ni07] in which she treats the ineffective termination of Kohn's algorithm for smooth weakly pseudoconvex domains of finite type from the viewpoint of adapting the ineffective argument of Diederich-Fornaess for the real-analytic case to the smooth case by using
Tougeron elements [To72] and Catlin's multitype [Ca84] to examine the possibilities of removing the difficulties of the smooth case such as those arising from the existence of non-identically-zero smooth function germs at a point whose derivatives of all orders vanish at that point.

\bigbreak At the end we include in this note an appendix which
presents some techniques of applying Skoda's theorem on ideal
generation [Sk72, Th.1, pp.555-556] which involve derivatives of
functions and Jacobian determinants.  Though most of these techniques are
not directly used in this note (except the use of (A.2) in (III.7) and the use (A.3) in (III.8)), they may be useful in reducing
vanishing orders of multiplier ideals in Kohn-type algorithms for
more general partial differential equations.

\bigbreak An earlier version of this note was posted in arxiv.org server as
arXiv:0706.4113.  In its Proposition (III.6) on the Multiplicity Estimate from Adjunction Formula, when we differentiate a given multiplier inside the formation of
a Jacobian determinant to construct another multiplier, we should have performed the differentiation as many times as the multiplicity of the given multiplier instead of performing it only once in Proposition (III.6) there.  In this version we put in the correct number of differentiation.   In (III.10) below we explain why
the correct number of differentiation is necessary.  In order to make our argument more transparent
with minimum notational complexity, we first do the special case of complex dimension two.
As a result we change completely the presentation of our argument in this version.  The presentation here is more streamlined than in the earlier version.

\bigbreak Before we go into the main body of this note, we would like to make one remark about the meaning of the effective termination of Kohn's algorithm.  Kohn's algorithm tells us that multipliers can be produced by using Jacobian determinants or by taking roots.  More precisely, taking a root means choosing an element in the radical of the ideal formed by multipliers in the preceding steps.  The challenge in solving the problem of the effective termination of Kohn's algorithm is to come up with a procedure which specifies when a root should be taken and when a Jacobian determinant should be used.  The procedure should specify when and how to choose an element in the radical of the ideal formed by multipliers in the preceding steps.  It also should specify when and how to choose functions (either multipliers from the preceding steps or pre-multipliers in the sense of (III.6)) to form the Jacobian determinant.  If one simply follows the rule of always giving priority to taking roots or simply follows the rule of always giving priority to taking Jacobian determinants, it is very easy to come up with simple counter-examples which show that such simple-minded rules could not in general yield an effective termination of Kohn's algorithm.  Of course, the uselessness of such simple-minded rules does not mean that Kohn's algorithm cannot be effectively terminated.  The problem of the effective termination of Kohn's algorithm is the determination of a selection rule to specify when and how to take an element in the radical of an ideal formed by multipliers from the previous steps and to specify when and how to take Jacobian determinants to form new multipliers.  The purpose of this note is to present such a selection rule and thereby demonstrate the effective termination of Kohn's algorithm.

\eject\noindent{\bf Part I. Global Regularity, Subellipticity,
Kohn's Algorithm, and Special Domains}

\bigbreak\noindent(I.1) {\it The Setting.}  We start out with the
setting of a bounded domain $\Omega$ in ${\mathbb C}^n$ with smooth
boundary defined by $r<0$ for some smooth function $r$ on an open
neighborhood of the topological closure of $\Omega$. We assume that
$\Omega$ is weakly pseudoconvex at all its boundary points in the
sense that at any boundary point $P$ of $\Omega$ the $(1,1)$-form
$\partial\bar\partial r$ is weakly positive definite when restricted
to the complex tangent space $\left\{\partial r=0\right\}$ of the
boundary $\partial\Omega$ of $\Omega$ at $P$. (To be more precise,
we should have said that $\sqrt{-1}\,\partial\bar\partial r$ is
weakly positive definite instead of $\partial\bar\partial r$ being
weakly positive definite, but for notational simplicity we will drop
the factor $\sqrt{-1}$ if there is no risk of any confusion.)

\bigbreak The {\it type} $m$ at a point $P$ of the boundary
$\partial\Omega$ of $\Omega$ is the supremum of the normalized
touching order
$$\frac{{\rm ord}_0\left(r\circ\varphi\right)}{{\rm ord}_0\varphi},$$
to $\partial\Omega$, of all local holomorphic curves
$\varphi:\Delta\to{\mathbb C}^n$ with $\varphi(0)=P$, where $\Delta$
is the open unit $1$-disk and ${\rm ord}_0$ is the vanishing order
at the origin $0$. A point $P$ of the boundary $\partial\Omega$ of
$\Omega$ is said to be of {\it finite type} if the type $m$ at $P$
is finite.  This notion of finite type was introduced by D'Angelo
[DA79].  For notational convenience we also call $m$ the order of
finite type instead of just the ``type'' to indicate that it is in
the sense of D'Angelo's notion of finite type.

\bigbreak Our goal is to obtain {\it global regularity} for any
smooth weakly pseudoconvex domain $\Omega$ of finite type in the
sense that, for any $\bar\partial$-closed smooth $(0,1)$-form $f$ on
$\Omega$ which is smooth up to the boundary of $\Omega$, the
solution $u$ of $\bar\partial u=f$ on $\Omega$ with $u$ orthogonal
to all holomorphic functions on $\Omega$ must also be smooth up to
the boundary of $\Omega$.  In this note by ``smoothness'' we mean
infinite differentiability.  For notational simplicity we formulate
global regularity only for $(0,1)$-forms. It can be similarly
formulated for $(0,q)$-forms for a general $q$.

\bigbreak Global regularity is a consequence of the {\it subelliptic
estimate}, which is defined as follows.  For any
$P\in\partial\Omega$ there exist an open neighborhood $U$ of $P$ in
${\mathbb C}^n$ and positive numbers $\epsilon$ and $C$ satisfying
$$
\||g|\|_\epsilon^2\leq C\left(\|\bar\partial g\|^2+\|\bar\partial^*
g\|^2+\|g\|^2\right)
$$
for any $(0,1)$-form $g$ supported on $U\cap\bar\Omega$ which is in
the domain of $\bar\partial$ and $\bar\partial^*$, where
$\||\cdot|\|_\epsilon$ is the $L^2$ norm on $\Omega$ involving
derivatives up to order $\epsilon$ in the boundary tangential
direction of $\Omega$ and $\|\cdot\|$ is the usual $L^2$ norm on
$\Omega$ without involving any derivatives, and $\bar\partial^*$ is
the actual adjoint of $\bar\partial$ with respect to $\|\cdot\|$.

\bigbreak Kohn [Ko79] introduced the following notion of multipliers
to obtain the subelliptic estimate.  At a point $P$ of
$\partial\Omega$ a smooth function germ $f$ at $P$ is called a {\it
subelliptic multiplier} (or simply called a {\it multiplier}) if
there exist some open neighborhood $U$ of $P$ in ${\mathbb C}^n$ and
some positive numbers $\epsilon$ and $C$ (all three depending on
$f$) such that
$$
\||fg|\|_\epsilon^2\leq C\left(\|\bar\partial
g\|^2+\|\bar\partial^* g\|^2+\|g\|^2\right)\leqno{({\rm I}.1.1)}
$$
for any $(0,1)$-form $g$ supported on $U\cap\bar\Omega$ which is in
the domain of $\bar\partial$ and $\bar\partial^*$.  We call the
positive number $\varepsilon$ an {\it order of subellipticity} for
the multiplier $f$.  We also call a subelliptic multiplier a {\it
scalar multiplier} to emphasize its difference from
vector-multipliers introduced below. The collection of all
multipliers at $P$ forms a ideal which is called the {\it multiplier
ideal} and is denoted by $I_P$.

\medbreak A germ of a smooth $(1,0)$-form $\theta$ at $P$ is called
a {\it subelliptic vector-multiplier} if there exist some open
neighborhood $U$ of $P$ in ${\mathbb C}^n$ and some positive numbers
$\epsilon$ and $C$ (all three depending on $\theta$) such that
$$
\||\bar\theta\cdot g|\|_\epsilon^2\leq C\left(\|\bar\partial
g\|^2+\|\bar\partial^* g\|^2+\|g\|^2\right)\leqno{({\rm I}.1.2)}
$$
for any $(0,1)$-form $g$ supported on $U\cap\bar\Omega$ which is in
the domain of $\bar\partial$ and $\bar\partial^*$, where
$\bar\theta\cdot g$ is the function obtained by taking the pointwise
inner product of the complex-conjugate $\bar\theta$ of $\theta$ with
$g$ with respect to the Euclidean metric of ${\mathbb C}^n$. We call
the positive number $\varepsilon$ an {\it order of subellipticity}
for the vector-multiplier $\theta$.  The collection of all
vector-multipliers at $P$ forms a module which is called the {\it
vector-multiplier module} and is denoted by $A_P$.

\medbreak The subelliptic estimate holds at a boundary point $P$ of
$\Omega$ if a nonzero constant function belongs to $I_P$.  Kohn introduced
the following algorithm to generate elements of $I_P$.

\begin{itemize}\item[(A)] {\it Initial Membership}.

\begin{itemize}\item[(i)] $r\in I_P$.

\item[(ii)] $\partial\bar\partial_j r$ belongs to $A_P$ for every
$1\leq j\leq n-1$ if $\partial r=dz_n$ at $P$ for some local
holomorphic coordinate system $\left(z_1,\cdots,z_n\right)$ centered
at $P$, where $\partial_j$ means $\frac{\partial}{\partial z_j}$.
\end{itemize}
\item[(B)] {\it Generation of New Members}.

\begin{itemize}\item[(i)]  If $f\in I_P$, then $\partial f\in A_P$.

\item[(ii)] If $\theta_1,\cdots,\theta_{n-1}\in A_P$, then the
coefficient of
$$
\theta_1\wedge\cdots\wedge\theta_{n-1}\wedge\partial r$$ is in
$I_P$. \end{itemize}\item[(C)] {\it Real Radical Property}.

\smallbreak\noindent If $g\in I_P$ and
$\left|f\right|^m\leq\left|g\right|$ for some positive integer $m$,
then $f\in I_P$.
\end{itemize}

\bigbreak\noindent The key point of Kohn's algorithm is to allow
certain differential operators to lower the vanishing orders of
multipliers so that eventually one can get a nonzero constant as a
multiplier.  However, there are two limitations on this process of
differentiation to lower vanishing orders.  One is that only
$(1,0)$-differentiation is allowed (as formulated in (B)(i) above).
The other is that only determinants of coefficients of
$(1,0)$-differentials can be used (as formulated in (B)(ii) above).
Besides using differentiation the ``real radical property'' allows
root-taking to reduce vanishing orders (as formulated in (C) above).

\bigbreak For the proof of the effective termination of the Kohn
algorithm, to keep track of the effectiveness we will assign a
positive number to a scalar multiplier (respectively a vector-multiplier)
which is an order of subellipticity for the scalar multiplier (respectively
vector multiplier).  We call such a
positive number an {\it assigned order of subellipticity}.  In
addition to keeping track of the number and the nature of the steps of
the Kohn algorithm used, the effectiveness of the termination of the
Kohn algorithm seeks to keep track of the assigned orders of
subellipticity for the individual scalar multipliers and
vector-multipliers so that the final nowhere zero multiplier has an
effective positive number as its assigned order of subellipticity.
 Note that the assigned order of subellipticity of a scalar
multiplier or vector-multiplier constructed in the Kohn algorithm is
in general not the maximum $\varepsilon$ for the inequality (I.1.1)
or (I.1.2).

\medbreak  We will adopt the following rule of giving to a scalar
multiplier or a vector-multiplier in the Kohn algorithm its assigned
order of multiplicity.  As its assigned order of subellipticity we
give the scalar multiplier $r$ the number $1$. As its assigned
order of subellipticity we give the number $\frac{1}{2}$ to the
vector-multiplier $\partial\bar\partial_j r$ at $P$ for every $1\leq
j\leq n-1$ is $\frac{1}{2}$ if $\partial r=dz_n$ at $P$ for some
local holomorphic coordinate system $\left(z_1,\cdots,z_n\right)$
centered at $P$.  If the assigned order of subellipticity of the
scalar multiplier $f$ is $\varepsilon$, then we give to the
vector-multiplier $\partial f$ as its assigned order of
subellipticity the number $\frac{\varepsilon}{2}$. If the minimum of
the assigned orders of subellipticity of the vector-multipliers
$\theta_1,\cdots,\theta_{n-1}$ is $\varepsilon$, then we give to the
same $\varepsilon$ to the coefficient of
$$
\theta_1\wedge\cdots\wedge\theta_{n-1}\wedge\partial r$$ as a scalar
multiplier as its assigned order of subellipticity. If the assigned
order of subellipticity of the scalar multiplier $g$ is
$\varepsilon$ and if $\left|f\right|^m\leq\left|g\right|$ for some
positive integer $m$, then we give the number
$\frac{\varepsilon}{m}$ to the scalar multiplier $f$ as its assigned
order of subellipticity.

\bigbreak\noindent(I.2) {\it Algebraic-Geometric Description of
Finite Type for Special Domains.} A special domain $\Omega$ in
${\mathbb C}^{n+1}$ (with coordinates $w,z_1,\cdots,z_n$) is a
bounded domain given by
$$
{\rm
Re\,}w+\sum_{j=1}^N\left|F_j\left(z_1,\cdots,z_n\right)\right|^2<0,\leqno{({\rm
I}.2.1)}
$$
where $F_j\left(z_1,\cdots,z_n\right)$ defined on some open
neighborhood of $\bar\Omega$ in ${\mathbb C}^{n+1}$ depends only on
the variables $z_1,\cdots,z_n$ and is holomorphic in
$z_1,\cdots,z_n$ for each $1\leq j\leq N$.  In what follows, when we
consider the subelliptic estimate at a boundary point $P$ of
$\Omega$ and its type of finite order, if there is no confusion we
will assume without loss of generality and without explicit mention
that the point $P$ is the origin of the coordinates
$w,z_1,\cdots,z_n$ and that $F_j$ vanishes at $P$ for each $1\leq
j\leq N$. Such special domains were introduced by Kohn [Ko79, p.115,
\S 7].

\medbreak To quantitatively describe finite type at the origin in a
way which is more convenient to use, we introduce a positive number
$p$ which is the smallest positive integer such that
$$
\left|z\right|^p\leq A_1\sum_{j=1}^N\left|F_j(z)\right|\leqno{({\rm
I}.2.2)}
$$
on some open neighborhood of the origin in ${\mathbb C}^n$ for some
positive constant $A_1$, where $z=\left(z_1,\cdots,z_n\right)$ and
$\left|z\right|^2=\sum_{\ell=1}^n\left|z_\ell\right|^2$.  We will
verify below in (I.3) that the order of finite type at the origin
$P$ is equal to $2p$.

\medbreak We are going to introduce also two other effectively
comparable ways to describe finite type which are both
algebraic-geometrical. The first one is the following.  Let
${\mathfrak m}$ be the maximum ideal ${\mathfrak m}_{{\mathbb
C}^n,0}$ of ${\mathbb C}^n$ at the origin. Let ${\mathcal I}$ be the
ideal on ${\mathbb C}^n$ generated by holomorphic function germs
$F_1,\cdots,F_N$ on ${\mathbb C}^n$ at the origin. Let $q$ be the
smallest positive integer such that
$${\mathfrak m}^q\subset{\mathcal I}.\leqno{({\rm I}.2.3)}.$$
We will verify
below in (I.4) that the number $p$ is related to the number $q$ by
the inequality $p\leq q\leq(n+2)p$. This inequality is far from
being sharp.

\medbreak The second is the following.  Let $s$ be the dimension
over ${\mathbb C}$ of the quotient of ${\mathcal O}_{{\mathbb
C}^n,0}$ by the ideal generated in it by the holomorphic function
germs $F_1,\cdots,F_N$ on ${\mathbb C}^n$ at the origin. We will
verify below in (I.5) that the number $q$ is related to the number
$s$ by the two inequalities $q\leq s$ and $s\leq{n+q-1\choose q-1}$.
Again this pair of inequalities is far from being sharp.

\medbreak When we prove the effective termination of Kohn's
algorithm for special domains, we will in different contexts choose
to use one of the three effectively comparable descriptions $p$,
$q$, or $s$ of the order of finite type.

\medbreak Let $U$ be an open neighborhood of $0$ in ${\mathbb C}$
and $\psi: U\to {\mathbb C}^{n+1}$ be a holomorphic map with
$\psi(0)=0$.  Write $\psi=\left(\psi_0,\psi_1\right)$ such that
$\psi_0:U\to{\mathbb C}$ and $\psi_1:U\to{\mathbb C}^n$. For $j=0,1$
the vanishing order ${\rm ord}_0\psi_j$ at $0$ of $\psi_j$ is the
positive integer $s$ such that
$$\psi_j(0)=\left(d\psi_j\right)(0)=\cdots=\left(d^{s-1}\psi_j\right)(0)=0$$
and $\left(d^s\psi_j\right)(0)\not=0$.  This positive integer $s$
can also be described as the largest integer such that
$$
\left|\psi_j(\zeta)\right|\leq A_{2,j}\left|\zeta\right|^s
$$
for some positive constant $A_{2,j}$ when the coordinate $\zeta$ of
${\mathbb C}$ is small. The vanishing order ${\rm
ord}_0\left(\psi\right)$ at the origin of $\psi$ is equal to the
minimum of the vanishing orders ${\rm ord}_0\left(\psi_0^*w\right)$
and ${\rm ord}_0\left(\psi_1^*z_j\right)$ of the holomorphic
functions $\psi_0^*w$ and $\psi_1^*z_j$ of $\zeta$ at $\zeta=0$ for
$1\leq j\leq n$.

\medbreak In the expansion of the function
$$
\psi^*r={\rm
Re}\left(\left(\psi_0^*w\right)\left(\zeta\right)\right)+\sum_{j=1}^N
\left|F_j\left(\psi_1\left(\zeta\right)\right)\right|^2
$$
in $\zeta,\,\bar\zeta$, only terms of the form $a_\nu\zeta^\nu$ and
$b_\nu\bar\zeta^\nu$ with $\nu\geq 1$  ({\it i.e.,} purely
holomorphic or purely antiholomorphic terms in $\zeta$ ) can occur
in ${\rm Re}\left(\left(\psi_0^*w\right)\left(\zeta\right)\right)$
and only terms of the form $c_{\mu,\nu}\zeta^\mu\bar\zeta^\nu$ with
$\mu\geq 1$ and $\nu\geq 1$ ({\it i.e.,} never purely holomorphic or
purely antiholomorphic terms in $\zeta$) occur in
$\sum_{j=1}^N\left|F_j\left(\psi_1\left(\zeta\right)\right)\right|^2$
(where $a_\nu, b_\nu, c_{\mu,\nu}$ are complex constants).  Since
there is no possibility at all of any term from ${\rm
Re}\left(\left(\psi_0^*w\right)\left(\zeta\right)\right)$ canceling
a term from
$\sum_{j=1}^N\left|F_j\left(\psi_1\left(\zeta\right)\right)\right|^2$,
it follows that the vanishing order at $0$ of $ \psi^*r$ must be
equal always to the minimum of the order at $0$ of ${\rm
Re}\left(\left(\psi_0^*w\right)\left(\zeta\right)\right)$ and the
order at $0$ of
$\sum_{j=1}^N\left|F_j\left(\psi_1\left(\zeta\right)\right)\right|^2$.
Thus
$$
\frac{{\rm ord}_0\psi^*r}{{\rm ord}_0\psi}=\frac{\min\left({\rm
ord}_0\psi_0^*w,\,{\rm
ord}_0\psi_1^*\sum_{j=1}^N\left|F_j\right|^2\right)}{\min\left({\rm
ord}_0\,\psi_0^*w,\frac{1}{2}\,{\rm
ord}_0\psi_1^*\sum_{j=1}^n\left|z_j\right|^2\right)}.
$$

\bigbreak\noindent(I.3) {\it Lemma.}  Let $p$ be the smallest
positive integer which satisfies (I.2.2) for some positive constant
$A_1$. Then the order $t$ of finite type at the origin for the
special domain $\Omega$ given by (I.2.1) is equal to $2p$.

\medbreak\noindent{\it Proof.}  We are going to prove $t=2p$ by
proving the two inequalities $t\leq 2p$ and $t\geq 2p$.  We first
prove the inequality $t\leq 2p$. From the definition of the order
$t$ of finite type we know that there exist some open neighborhood
$U$ of $0$ in ${\mathbb C}$ and some holomorphic map
$\psi=\left(\psi_0,\psi_1\right): U\to {\mathbb C}^{n+1}={\mathbb
C}\times{\mathbb C}^n$ with $\psi(0)=0$ such that
$$
t=\frac{{\rm ord}_0\psi^*r}{{\rm ord}_0\psi}=\frac{\min\left({\rm
ord}_0\psi_0^*w,\,{\rm
ord}_0\psi_1^*\sum_{j=1}^N\left|F_j\right|^2\right)}{\min\left({\rm
ord}_0\,\psi_0^*w,\frac{1}{2}\,{\rm
ord}_0\psi_1^*\sum_{j=1}^n\left|z_j\right|^2\right)}.
$$
We let $$\alpha={\rm ord}_0\psi_0^*w,\quad \beta=\frac{1}{2}\,{\rm
ord}_0\psi_1^*\sum_{j=1}^n\left|z_j\right|^2,\quad\gamma={\rm
ord}_0\psi_1^*\sum_{j=1}^N\left|F_j\right|^2.$$ From the definition
of $p$ we have $\gamma\leq 2p\beta$.  We differentiate among the
following three cases.
\begin{itemize}\item[]Case 1. $\alpha\leq
\beta$.\item[]Case 2. $\beta<\alpha\leq 2p\beta$\item[]Case 3.
$\alpha>2p\beta$.\end{itemize} For Case 1, we have
$$t=\frac{\min(\alpha,\gamma)}{\min(\alpha,\beta)}\leq\frac{\min(\alpha,2p\beta)}{\min(\alpha,\beta)}=\frac{\alpha}{\alpha}\leq 2p.$$
For Case 2, we have
$$t=\frac{\min(\alpha,\gamma)}{\min(\alpha,\beta)}\leq\frac{\min(\alpha,2p\beta)}{\min(\alpha,\beta)}=\frac{\alpha}{\beta}\leq\frac{2p\beta}{\beta}=2p.$$
For Case 3, we have
$$t=\frac{\min(\alpha,\gamma)}{\min(\alpha,\beta)}\leq\frac{\min(\alpha,2p\beta)}{\min(\alpha,\beta)}=\frac{2p\beta}{\beta}=2p.$$
Thus in all three cases $t\leq 2p$.

\medbreak We are now going to prove the other inequality $2p\leq t$.
We use a simultaneous resolution of embedded singularities
$\pi:\tilde W\to W$ for some open neighborhood $W$ of the origin in
${\mathbb C}^n$ with exceptional hypersurfaces
$\left\{Y_j\right\}_{j=1}^J$ in $\tilde W$ in normal crossing so
that the pullback $\pi^*{\mathfrak m}$ of the maximum ideal on $W$
at the origin is equal to the ideal sheaf of the divisor
$\sum_{j=1}^J\sigma_jY_j$ for some nonnegative integers
$\sigma_1,\cdots,\sigma_J$ and the pullback $\pi^*{\mathcal I}$ of
the ideal sheaf ${\mathcal I}$ on ${\mathbb C}^n$ generated by
$F_1,\cdots,F_N$ is equal to the ideal sheaf of the divisor
$\sum_{j=1}^J\tau_jY_j$ for some nonnegative integers
$\tau_1,\cdots,\tau_J$.

\medbreak Since $p$ is the smallest integer which satisfies
condition (I.2.2) for some positive constant $A_1$, it follows that
$p\sigma_j\leq\tau_j$ for all $1\leq j\leq J$. Take any $1\leq
j_0\leq J$ with $\sigma_{j_0}>0$ and $0\in\pi\left(Y_{j_0}\right)$
such that there is a regular point $Q$ of some $Y_{j_0}$ with the
property that $\pi\left(Q\right)=0$ and $Q$ does not belong to any
other $Y_j$ with $j\not=j_0$. Take a local regular complex curve
$\tilde C$ in $\tilde W$ represented by a holomorphic map
$\tilde\varphi:U\to\tilde W$ from some open neighborhood $U$ of $0$
in ${\mathbb C}$ to $\tilde W$ such that $\tilde\varphi(0)=Q$ and
the local complex curve $\tilde C$ is transversal to $Y_{j_0}$ and
disjoint from any other $Y_j$ with $j\not=j_0$.  Now we define a
holomorphic map $\psi=\left(\psi_0,\psi_1\right): U\to {\mathbb
C}^{n+1}={\mathbb C}\times{\mathbb C}^n$ by $\psi_0\equiv 0$ and
$\psi_1=\pi\circ\varphi$.  Then
$$
\displaylines{\frac{{\rm ord}_0\psi^*r}{{\rm
ord}_0\psi}=\frac{\min\left({\rm ord}_0\psi_0^*w,\,{\rm
ord}_0\psi_1^*\sum_{j=1}^N\left|F_j\right|^2\right)}{\min\left({\rm
ord}_0\,\psi_0^*w,\frac{1}{2}\,{\rm
ord}_0\psi_1^*\sum_{j=1}^n\left|z_j\right|^2\right)}\cr =\frac{{\rm
ord}_0\psi_1^*\sum_{j=1}^N\left|F_j\right|^2}{\frac{1}{2}\,{\rm
ord}_0\psi_1^*\sum_{j=1}^n\left|z_j\right|^2}
=\frac{2\tau_{j_0}}{\sigma_{j_0}}\geq 2p.\cr}
$$
By the definition of $t$ we have $t\geq 2p$.  Putting the two
inequalities $t\leq 2p$ and $t\geq 2p$ together, we get $t=2p$.
Q.E.D.

\bigbreak\noindent(I.4) {\it Lemma.}  Let $p$ be the smallest
positive integer which satisfies (I.2.2) for some positive constant
$A_1$. Let $q$ be the smallest positive integer such that the $q$-th
power ${\mathfrak m}^q$ of the maximum ideal ${\mathfrak m}$ of
${\mathbb C}^n$ at the origin is contained in the ideal ${\mathcal
I}$ generated by the holomorphic function germs $F_1,\cdots,F_N$ on
${\mathbb C}^n$ at the origin. Then $p\leq q\leq(n+2)p$.

\medbreak\noindent{\it Proof.} From the definition of $q$ it follows
that
$$z_\ell^q\in{\mathfrak m}^q\subset{\mathcal I}=\sum_{j=1}^N{\mathcal O}_{{\mathbb C}^n,0}F_j$$
for every $1\leq\ell\leq n$ it follows that
$$
\left|z_\ell^q\right|\leq \tilde
A_\ell\sum_{j=1}^N\left|F_j\left(z\right)\right| $$ for some
positive constant $\tilde A_\ell$ for $1\leq\ell\leq n$ on some open
neighborhood of the origin in ${\mathbb C}^n$.  Hence
$$
\left|z\right|^q=\left(\sum_{\ell=1}^n\left|z_\ell\right|^2\right)^{\frac{q}{2}}
\leq n^{\frac{q}{2}}\max_{1\leq\ell\leq n}\left|z_\ell\right|^q\leq
n^{\frac{q}{2}}\left(\max_{1\leq\ell\leq n}\tilde
A_\ell\right)\sum_{j=1}^N\left|F_j\left(z\right)\right|$$ and $p\leq
q$ from the definition of $p$.

\medbreak For the proof of the inequality $p\leq(n+2)q$, we need the
following theorem of Skoda [Sk72, Th.1, pp.555-556].

\medbreak\noindent Let $D$ be a pseudoconvex domain in ${\bf C}^n$
and $\chi$ be a plurisubharmonic function on $D$. Let
$g_1,\ldots,g_m$ be holomorphic functions on $D$. Let $\alpha>1$ and
$\ell=\inf(n,m-1)$. Then for every holomorphic function $F$ on $D$
such that
$$
\int_D|F|^2|g|^{-2\alpha \ell-2}e^{-\chi}<\infty,
$$
there exist holomorphic functions $f_1,\ldots,f_m$ on $\Omega$ such
that
$$
F=\sum_{j=1}^mg_jf_j
$$
and
$$
\int_D|f|^2|g|^{-2\alpha\ell}e^{-\chi}
\leq{\alpha\over{\alpha-1}}\int_D|F|^2|g|^{-2\alpha
\ell-2}e^{-\chi},
$$
where
$$
|g|=\left(\sum_{j=1}^m|g_j|^2\right)^{1\over 2},\qquad
|f|=\left(\sum_{j=1}^m|f_j|^2\right)^{1\over 2}.
$$

\medbreak For nonnegative integers $\gamma_1,\cdots,\gamma_n$ with
$\gamma_1+\cdots+\gamma_n=(n+2)p$ we apply Skoda's theorem to the
case of
$$\displaylines{F=\prod_{j=1}^n z_j^{\gamma_j},\quad m=N+n,\cr
\chi\equiv
0,\quad\left\{g_1,\cdots,g_m\right\}=\left\{F_1,\cdots,F_N,0,\cdots,0\right\},\cr
\ell=n,\quad\alpha=\frac{n+1}{n},\cr}$$ with $D$ being some small
open ball neighborhood of the origin in ${\mathbb C}^n$, to conclude
from (I.2.2) that
$$
\prod_{j=1}^n z_j^{\gamma_j}\in{\mathcal I}=\sum_{j=1}^N{\mathcal
O}_{{\mathbb C}^n,0}F_j.
$$
Hence $q\leq(n+2)p$.  Q.E.D.

\bigbreak\noindent(I.5) {\it Lemma.}  Let $q$ be the smallest
positive integer such that the $q$-th power ${\mathfrak m}^q$ of the
maximum ideal ${\mathfrak m}$ of ${\mathbb C}^n$ at the origin is
contained in the ideal ${\mathcal I}$ generated by the holomorphic
function germs $F_1,\cdots,F_N$ on ${\mathbb C}^n$ at the origin.
Let $s$ be the dimension over ${\mathbb C}$ of the quotient of
${\mathcal O}_{{\mathbb C}^n,0}$ by the ideal ${\mathcal I}$.  Then
$q\leq s$ and $s\leq{n+q-1\choose q-1}$.

 \medbreak\noindent{\it Proof.} Let $R={\mathcal O}_{{\mathbb
C}^n,0}$.  Since ${\mathfrak m}^q\subset{\mathcal I}$, it follows
that
$$
s=\dim_{\mathbb C} R\left/{\mathcal I}\right.\leq\dim_{\mathbb C}
R\left/{\mathfrak m}^q\right.={n+q-1\choose q-1}.
$$
On the other hand, we consider the following sequence of nested
${\mathbb C}$-linear subspaces of the finite-dimensional ${\mathbb
C}$-vector spaces $R\left/{\mathcal I}\right.$ of complex dimension
$s$.
$$
R\left/{\mathcal I}\right.\supset{\mathfrak m}\left(R\left/{\mathcal
I}\right.\right)\supset{\mathfrak m}^2\left(R\left/{\mathcal
I}\right.\right)\supset\cdots\supset{\mathfrak
m}^\ell\left(R\left/{\mathcal I}\right.\right)\supset{\mathfrak
m}^{\ell+1}\left(R\left/{\mathcal I}\right.\right)\supset\cdots
$$
There exists some $1\leq\ell\leq s$ such that
$$
{\mathfrak m}^\ell\left(R\left/{\mathcal I}\right.\right)={\mathfrak
m}^{\ell+1}\left(R\left/{\mathcal I}\right.\right).
$$
By Nakayama's lemma we have
$$
0={\mathfrak m}^\ell\left(R\left/{\mathcal
I}\right.\right)=\left({\mathfrak m}^\ell+{\mathcal
I}\right)\left/{\mathcal I}\right.
$$
and we conclude that ${\mathfrak m}^\ell+{\mathcal I}={\mathcal I}$
and ${\mathfrak m}^\ell\subset{\mathcal I}$, which implies that
$q\leq s$.  Q.E.D.

\bigbreak Later we will need the following corollary of Lemma (I.5)
which is a version of the effective Nullstellensatz in terms of
multiplicity.

\bigbreak\noindent(I.6) {\it Lemma.}  Let ${\mathcal I}$ be an ideal
in ${\mathcal O}_{{\mathbb C}^n,0}$ such that its multiplicity is no
more than some positive integer $m$ in the sense that the complex dimension
of the quotient of ${\mathcal O}_{{\mathbb C}^n,0}$ by ${\mathcal I}$ is no
more than $m$.
Let $f$ be a holomorphic
function germ on ${\mathbb C}^n$ at the origin which vanishes at the origin.  Then $f^{m^2}$
belongs to ${\mathcal I}$.

\medbreak\noindent{\it Proof.} Let ${\mathcal
I}=\bigcap_{j=1}^J{\mathcal Q}_j$ be the primary decomposition of
the ideal ${\mathcal I}$ in the Noetherian ring ${\mathcal
O}_{{\mathbb C}^n,0}$ and let ${\mathcal P}_j$ be the radical of the
primary ideal ${\mathcal Q}_j$.  Since $m$ is the multiplicity of
${\mathcal I}$, we have $J\leq m$.  Let $Z_j$ be the subvariety germ
of ${\mathbb C}^n$ at the origin whose ideal at the origin is
${\mathcal P}_j$ and let $Z_j^0$ be the set of all regular points of
$Z_j$, which without loss of generality we can assume to be
connected.  Let $n_j$ be the complex codimension of $Z_j$ at the
origin.  Take a generic point $Q_j$ in $Z^0_j$ and let $\Pi_j$ be a
complex ${\mathbb C}$-linear subspace of complex dimension $n_j$
${\mathbb C}^n$ which is transversal to $Z^0_j$ at the point $Q_j$.
The multiplicity of the ideal at $Q_j$ induced by ${\mathcal Q}_j$
is no more than $m$ (if we assume without loss of generality that
$Q_j$ is sufficiently close to the origin).   Let ${\mathcal J}_j$
be the ideal on $\Pi_j\approx{\mathbb C}^{n_j}$ at $Q_j$ induced by
${\mathcal Q}_j$.  Then the multiplicity of ${\mathcal
J}_j$ at the point $Q_j$ is no more than $m$.  By Lemma (I.5)
applied to $\Pi_j\approx{\mathbb C}^{n_j}$ and the ideal ${\mathcal
J}_j$ on $\Pi_j\approx{\mathbb C}^{n_j}$, we conclude that the
holomorphic function germ $\left(f|_{\Pi_j}\right)^m$ on
$\Pi_j\approx{\mathbb C}^{n_j}$ at $Q_j$ belongs to the ideal
${\mathcal J}_j$.  Since $Q_j$ is a generic point in $Z^0_j$ and
since ${\mathcal Q}_j$ is a primary ideal, it follows that the
holomorphic function germ $f^m$ on ${\mathbb C}^n$ at the origin
belongs to ${\mathcal Q}_j$. From $J\leq m$ it follows that the
holomorphic function germ $f^{m^2}$ on ${\mathbb C}^n$ at the origin
belongs to the product of ${\mathcal Q}_j$ for $1\leq j\leq J$. In
particular, holomorphic function germ $f^{m^2}$ on ${\mathbb C}^n$
at the origin belongs to ${\mathcal I}$.  Q.E.D.

\eject\noindent{\bf Part II. Algebraic Formulation and Sketches
of Techniques}

\bigbreak\noindent(II.1) {\it Algebraic Formulation of Kohn's
Algorithm for Special Domains}.  The effective termination of Kohn's
algorithm for a special domain is reduced to the purely algebraic-geometric
description of items (i) through (vii) listed below.  For the case of a special domain the setting is as follows.  We have a special domain $\Omega$ in ${\mathbb C}^{n+1}$ (with
coordinates $w,z_1,\cdots,z_n$) defined by
$$
r:={\rm
Re\,}w+\sum_{j=1}^N\left|F_j\left(z_1,\cdots,z_n\right)\right|^2<0,
$$
where for each $1\leq j\leq
N$, $F_j\left(z_1,\cdots,z_n\right)$ is a holomorphic function vanishing at the origin which is independent of $w$ and is defined on some open
neighborhood of $\bar\Omega$ in ${\mathbb C}^{n+1}$.  The boundary point of $\Omega$ under consideration is the origin of ${\mathbb C}^{n+1}$.

\medbreak In this setting, first of all, from $dr=dw$ at the origin and
$$
\partial\bar\partial r=\sum_{j=1}^N dF_j\wedge\overline{dF_j}
$$
we conclude from (I.1)(A)(ii) and standard techniques of estimates
in Kohn's theory of multipliers [Ko79] that $dF_j$ is a vector
multiplier which can be given $\frac{1}{4}$ as its assigned order of
subellipticity, because the vector-multiplier
$$
\partial\left(\frac{\partial r}{\partial{\overline{z_j}}}\right)=\sum_{\ell=1}^N \overline{
\left(\frac{\partial h_\ell}{\partial z_j}\right)}\,dh_\ell
$$
at the origin can be given $\frac{1}{2}$ as its assigned order of
subellipticity for $1\leq j\leq n$.

\medbreak\noindent(i) We start out with the $N$ given holomorphic
function germs $F_1,\cdots,F_N$ on ${\mathbb C}^n$ at the origin
with the origin as their only common zero-point.  The multiplicity $q$
of the ideal generated by
$F_1,\cdots,F_N$ at the origin is what we use for effectiveness statements.  That is, a number is considered effective if it can be estimated by an explicit expression in $q$.

\medbreak\noindent(ii) Select $n$ ${\mathbb C}$-linear combinations
$g_1,\cdots,g_n$ of $F_1,\cdots,F_N$.

\medbreak\noindent(iii) Form the Jacobian determinant of
$g_1,\cdots,g_n$ with respect to $z_1,\cdots,z_n$.

\medbreak\noindent(iv) Take the ideal $I$ generated by all such
Jacobian determinants.

\medbreak\noindent(v) Choose a finite subset
$\varphi_1,\cdots,\varphi_\ell$ of the radical $J$ of $I$ and let
$\sigma$ be a positive number such that $\left(\varphi_j\right)^\sigma\in I$
for $1\leq j\leq\ell$.

\medbreak\noindent(vi) Replace the set $F_1,\cdots,F_N$ by
$F_1,\cdots,F_N,\varphi_1,\cdots,\varphi_\ell$ and repeat the above
procedure.

\medbreak\noindent(vi) Repeat until we get to the point that
$\varphi_1$ can be chosen to be nonzero at the origin.

\medbreak\noindent(vii) Effectiveness means that we have an
effective number of steps and also an effective bound on $\sigma$ at each
step.

\bigbreak\noindent(II.2) {\it Sketch of Proof of Effectiveness for
Special Domains}.  We now give an overview of the logical framework
for the proof of the effective termination of the Kohn algorithm for
special domains.  Details for the derivation of the bounds of the
multiplicities of functions constructed from generic ${\mathbb
C}$-linear combinations and Jacobian determinants which occur in
this overview will not be explained here but will be presented
later in (III.3), (III.4), and (III.5).

\medbreak We start out with an ideal generated by holomorphic function
germs $F_1,\cdots,F_N$ on ${\mathbb C}^n$ at the origin whose common zero-set is the origin.  The multiplicity $q$ of the ideal generated by $F_1,\cdots,F_N$ at the origin is what we use for effectiveness statements. For
$n$ generic ${\mathbb C}$-linear combinations $g_1,\cdots,g_n$ of
$F_1,\cdots,F_N$ the multiplicity of the function $f$ defined by
$dg_1\wedge\cdots\wedge dg_n=f\left(dz_1\wedge\cdots\wedge
dz_n\right)$ is no more than $m_q$ at the origin, where $m_q$ is
some positive integer depending effectively on $q$ (see (III.5)). The
main idea is to use the procedure of replacing ${\mathbb C}^n$ by the subspace $V$ defined by the multiplier $f$ to cut down successively
on the dimension of the zero-set of multipliers while maintaining effectiveness.

\medbreak There are two difficulties here.  One difficulty is that the subspace defined by $f$ is
in general not regular.  The other difficulty is that we are allowed only to form
Jacobian determinants of ${\mathbb C}$-linear combinations $g_1,\cdots,g_n$
of $F_1,\cdots,F_N$
and not allowed to form the
Jacobian determinants of the restrictions of such ${\mathbb C}$-linear combinations
$g_1,\cdots,g_{n-1}$ to $V$.  The two difficulties are related.  If $V$ is nonsingular,
we could compute the Jacobian determinant of $g_1|_V,\cdots,g_{n-1}|_V$ by computing
the coefficient of $dz_1\wedge\cdots\wedge dz_n$ in $dg_1\wedge\cdots\wedge dg_{n-1}\wedge df$.

\medbreak When $V$ is singular at the origin, we have to differentiate $f$ not just once to form $df$ but as many times as the multiplicity of $V$.  To enable us to do it by using Jacobian determinants, we construct a Weierstrass polynomial $\tilde f$ in $z_n$ whose coefficients are functions of $g_1,\cdots,g_{n-1}$ so that $\tilde f$ vanishes on the subspace $V$ and therefore contains $f$ as a factor.  We then differentiate $\tilde f$ as many times as its multiplicity at the origin by applying the operator $dg_1\wedge\cdots
\wedge dg_{n-1}\wedge d\left(\cdot\right)$ to $\tilde f$ and making use of the fact that $\tilde f$ is a Weierstrass polynomial of the type described above.  To continue applying the operator $dg_1\wedge\cdots
\wedge dg_{n-1}\wedge d\left(\cdot\right)$ to $\tilde f$, we need to modify first the result from the previous differentiation by comparing on $V$ the Jacobian determinant $$\frac{\partial\left(g_1,\cdots,g_{n-1}\right)}{\partial\left(z_1,\cdots,z_{n-1}\right)}$$ with an appropriate polynomial $p\left(g_1,\cdots,g_{n-1}\right)$ of $g_1,\cdots,g_{n-1}$ and using the Real Radical Property of Kohn's algorithm in (I.1)(B) to replace $$\frac{\partial\left(g_1,\cdots,g_{n-1}\right)}{\partial\left(z_1,\cdots,z_{n-1}\right)}$$ by $p\left(g_1,\cdots,g_{n-1}\right)$.  The final result of differentiating $\tilde f$ this way as many times as the multiplicity of $\tilde f$ at the origin produces a new multiplier which defines on $V$ a subspace with effective multiplicity at the origin.  This way of cutting down on the dimension of the subspace defined by such effectively constructed multipliers gives the effective termination of Kohn's algorithm for special domains.

\medbreak In the details of the proof for special domains given below in (III.6), (III.7),
(III.8), and (III.9), we actually do not carry out completely the induction of cutting down on the dimension of the zero-set of effectively constructed multipliers.  A short-cut is used to simplify the process to reach the same goal (see (III.9)).

\bigbreak\noindent(II.3) {\it Modification for Effectiveness for
Real-Analytic Case.} Before we give the rigorous details of the
proof of the effective termination of the Kohn algorithm for special
domains, we would like to discuss how the techniques in the above
sketch for special domains in (II.2) can be modified for the general
real-analytic case.  We consider the following real-analytic case
where the weakly pseudoconvex domain of finite type is defined by
$r<0$ with
$r\left(z_1,\cdots,z_n,\overline{z_1},\cdots,\overline{z_n}\right)$
being real-analytic and vanishing at the origin (which is the
boundary point we consider).  The main idea is to let
$w_j=\overline{z_j}$ for $1\leq j\leq n$ and let $R$ be the ring of
convergent power series in $w_1,\cdots,w_n$ and consider the $n+1$
holomorphic function germs $H_0, H_1, \cdots, H_n$ on ${\mathbb
C}^n$ (with coordinates $z_1,\cdots,z_n$) at the origin with
coefficients in the ring $R$ defined as follows.
$$
\displaylines{
H_0\left(z_1,\cdots,z_n\right)=r\left(z_1,\cdots,z_n,w_1,\cdots,w_n\right),\cr
H_j\left(z_1,\cdots,z_n\right)=\frac{\partial r}{\partial
w_j}\left(z_1,\cdots,z_n,w_1,\cdots,w_n\right)\quad{\rm for\ \
}1\leq j\leq n,\cr}
$$
where the coefficients of the power series expansion of $H_j$ in
$z_1,\cdots,z_n$ are all elements of $R$ for $0\leq j\leq n$.  For
the complex Euclidean space ${\mathbb C}^n$ with coordinates
$z_1,\cdots,z_n$ we denote by ${\mathfrak m}_{{\mathbb C}^n,0}$ the
maximum ideal of ${\mathbb C}^n$ at the origin.  Finite type
condition for the domain $\left\{\,r<0\,\right\}$ at the origin
implies the statement that

\medbreak\noindent(II.3.1) there exists some effective positive
integer $q$ such that $R\left({\mathfrak m}_{{\mathbb
C}^n,0}\right)^q$ is contained in the ideal generated by
$H_0,H_1,\cdots,H_n$ in the ring $R\left\{z_1,\cdots,z_n\right\}$ of
convergent power series in $z_1,\cdots,z_n$ with coefficients in
$R$.

\medbreak The statement (II.3.1) simply follows directly from the
definition of finite type.  It can be regarded as the real-analytic
analog of condition (I.2.3) for a special domain.  Note that finite
type is actually much stronger than the statement (II.3.1).

\medbreak  In a way analogous to applying condition (I.2.3) to do an
inductive multiplier-construction process to obtain a nonzero
constant as a multiplier from the Kohn algorithm for a special
domain as described in (II.2), we now apply statement (II.3.1) to do the same inductive
multiplier-construction process with the difference that now the
coefficients of the power series of the function germs involved are
elements of $R=\left\{w_1,\cdots,w_n\right\}$ instead of just
${\mathbb C}$. One modification is needed for the inductive
multiplier-construction process. When we are in the case of a
special domain, we use $n$ generic ${\mathbb C}$-linear combinations
$g_1,\cdots,g_n$ of $F_1,\cdots,F_N$, but here in the real-analytic
case when we choose $n$ $R$-linear combinations $g_1,\cdots,g_n$ of
$H_0,\cdots,H_n$, one of $g_1,\cdots,g_n$ must be chosen to be
$H_0$. The reason for this modification is that we are not in the
special case where the domain is of the form
$$
{\rm
Re\,}w+r\left(z_1,\cdots,z_n,\overline{z_1},\cdots,\overline{z_n}\right)<0
$$
and when our domain is not in this special form we have to use
$\partial r$ to define the tangent space of type $(1,0)$ for the
boundary of the domain.

\medbreak Note that when we take $\partial G_1\wedge\cdots\wedge\partial G_{n-1}
\wedge\partial H_0$ for generic ${\mathbb C}$-linear combinations $G_1,\cdots,G_{n-1}$ of
$H_1,\cdots,H_n$, we are simply using (I.1)(A)(ii) and (I.1)(B)(ii) in Kohn's algorithm.

\medbreak The inductive multiplier-construction process in the
real-analytic case now gives us a nonzero element $f$ of $R$ instead
of a nonzero element of ${\mathbb C}$ in the case of a special
domain.  The main point is that, because of the finite type
condition the multiplicity of this element
$f\left(w_1,\cdots,w_n\right)$ of $R$ at $0$ is bounded effectively
by a constant depending on $n$ and the order of the finite type. Now
we consider the anti-holomorphic function germ $\tilde f$ on
${\mathbb C}^n$ at the origin defined by $\tilde
f=f\left(\overline{z_1},\cdots,\overline{z_n}\right)$ and consider
the complex conjugate $g$ of $\tilde f$.

\medbreak We let $V_1$ be the subspace germ defined
by the holomorphic function germ $g$ on ${\mathbb C}^n$ at the
origin. We then consider $V_1\times\overline{V_1}$ in ${\mathbb
C}^n\times{\mathbb C}^n$ instead of the full $2n$-dimensional
complex Euclidean space ${\mathbb C}^n\times{\mathbb C}^n$ itself
(with $z_1,\cdots,z_n,\overline{z_1},\cdots,\overline{z_n}$ being
the variables of ${\mathbb C}^n\times{\mathbb C}^n$). Let $R_1$ be
the ring of holomorphic function germs on $V_1$ at $0$ when $V_1$ is
considered as a subspace germ of $0$ in ${\mathbb C}^n$ at the origin
with coordinates $w_1,\cdots,w_n$. We now apply the inductive
process to obtain a holomorphic function germ $f_1$ on $V_1$ at $0$
(which is a subspace germ at $0$ of ${\mathbb C}^n$ with
variables $w_1,\cdots,w_n$).

\medbreak Now we consider the function germ $\tilde f_1$ obtained
from $f_1$ by replacing $w_1,\cdots,w_n$ by
$\overline{z_1},\cdots,\overline{z_n}$.  Let $g_1$ be the
complex-conjugate of $\tilde f_1$.  Let $V_2$ be a complete
intersection of codimension two in ${\mathbb C}^n$ at the origin
defined by two holomorphic functions which belong to the radical of
the ideal generated by $g_2$ and the ideal of $V_1$. We then
consider $V_2\times\overline{V_2}$ in ${\mathbb C}^n\times{\mathbb
C}^n$ instead of ${\mathbb C}^n\times{\mathbb C}^n$ itself (with
$z_1,\cdots,z_n,\overline{z_1},\cdots,\overline{z_n}$ being the
variables of ${\mathbb C}^n\times{\mathbb C}^n$). Let $R_2$ be the
ring of holomorphic function germs on $V_2$ at $0$ when $V_2$ is
considered as a subspace germ of $0$ in ${\mathbb C}^n$
with coordinates $w_1,\cdots,w_n$.  We now can continue with this
inductive subspace-construction process which so far yields for us
the subspace of complete intersection $V_1$ and $V_2$.
We continue with this inductive subspace-construction process to
get $V_{\ell+1}$ from $V_\ell$ for $1\leq\ell\leq n-1$ until we get
to the subspace $V_n$ of ${\mathbb C}^n$ which consists
only of the origin.  This then immediately gives us the effective
termination of Kohn's algorithm for the real-analytic case.  Again, instead of carrying
out completely the inductive argument of cutting down the dimension of the
subspace described above, it is also possible to use the analog of the short-cut technique given in (III.9).

\medbreak Another way of describing this modification is to redo the
algebraic-geometric argument for the case of a special domain but to
do it over a parameter space defined by the ring $R$.  The
coordinates for $R$ are the complex-conjugates of the coordinates for
the ambient space ${\mathbb C}^n$.   We can describe the
modification as redoing the algebraic-geometric argument for the
case of a special domain over ${\rm Spec}(R)$ instead of over the
single point ${\rm Spec}\left({\mathbb C}\right)$.  While the case
of a special domain yields effectively a nonzero element of
${\mathbb C}$ as a multiplier, the real-analytic case would yield
effectively a nonzero element of $R$.  Then we replace ${\mathbb
C}^n$ or by the subspace defined by this nonzero element of $R$ and
repeat the argument to get down to lower and lower dimensional
subspaces until we get to a single point, or we use the analog of the short-cut technique
given in (III.9).

\bigbreak\noindent(II.4) {\it Modification for Effectiveness for
Smooth Case.} We are going to have yet another discussion, this
time about modifying further the techniques in the above sketch
for special domains in (II.2) in order to handle the general
smooth case, before going into the rigorous details of the proof
of the effective termination of the Kohn algorithm for special
domains. Now suppose that we have a smooth bounded weakly
pseudoconvex domain $\Omega$ of finite type given by $r<0$ for
some smooth function $r$ defined on some neighborhood of the
topological closure $\bar\Omega$ of $\Omega$ in ${\mathbb C}^n$
and that the origin $0$ of ${\mathbb C}^n$ is a boundary point of
$\Omega$.

\medbreak Let $q$ be the positive integer which is the order of the
finite type of the
origin as a boundary point of $\Omega$. Let $r_N$ be
the $N$-th partial sum of the formal power series expansion of $r$
at the origin with respect to the coordinates $z_1,\cdots,z_n$ of
${\mathbb C}^n$.  We choose $N$ effectively large enough so that the
type of $r_N=0$ at the origin is also $q$.

\medbreak We apply Kohn's algorithm for the real-analytic case
to $r_N$. From the effectiveness for the real-analytic case
(II.3), we can find some positive integer $N_q$ which depends only
on $q$ and $n$ such that the assigned order of subellipticity $\varepsilon$
for the final nonzero constant multiplier from the effective Kohn
algorithm for $r_N$ satisfies $\varepsilon>\frac{1}{N_q}$.

\medbreak When we choose $N$ effectively large enough, for example,
$N>2N_p$, the effective termination of Kohn's algorithm for $r_N$
also gives the effective termination of Kohn's algorithm for $r$
with precisely the same steps and the same assigned order of subellipticity
for each step.  Note that this process of approximating $r$ by $r_N$
is very different from the approximation of a bounded smooth weakly
pseudoconvex domain of finite type by a real-analytic smooth weakly
pseudoconvex domain of finite type, which is in general not
possible.  The $N$-th partial sum $r_N$ is simply used as an
algebraic-geometric comparison guide to guarantee the effective
termination of Kohn's algorithm for the original smooth defining
function $r$.

\medbreak Note that when we do the approximation of $r$ by $r_N$, we are doing this approximation only at the boundary point under consideration and not using the approximation
along the normal directions of the zero-sets of multipliers from Kohn's algorithm for $r$.
The reason is that the purpose of the approximation is to use the effective
termination of Kohn's algorithm for the real-analytic function $r_N$ to conclude for a sufficiently large $N$ that the
corresponding steps result in the effective termination of Kohn's algorithm for the smooth function $r$.  The motivation for choosing this procedure of approximation is twofold.  One is that the notion of finite type at a boundary point of the weakly pseudoconvex domain depending only on the formal power series expansion of the defining function $r$ at that point.  The other is that the zero-sets of multipliers from Kohn's algorithm for $r$ are defined by the vanishing of smooth functions and it is not clear how one can do a real-analytic approximation along the normal directions of such zero-sets.  In our use of the approximation of $r$ by $r_N$, the zero-sets of multipliers from Kohn's algorithm for $r$ are different from the zero-sets of multipliers from Kohn's algorithm for $r_N$.  When we use the ``real radical property'' to produce multipliers from Kohn's algorithm for $r_N$, we simply perform the same operation for the corresponding but different zero-set in Kohn's algorithm for $r$.

\eject\noindent{\bf Part III. Details of Proof of Effective
Termination of Kohn's Algorithm for Special Domains}

\bigbreak\noindent(III.1) {\it Precise Formulation.} Let
$F_1,\cdots,F_N$ be holomorphic function germs on ${\mathbb C}^n$ at
the origin $0$. Assume that
$$s:=\dim_{\mathbb C}\left({\mathcal O}_{{\mathbb
C}^n,0}\left/\sum_{j=1}^N{\mathcal O}_{{\mathbb
C}^n,0}F_j\right.\right)<\infty
$$
so that the subscheme of ${\mathbb C}^n$ defined by
$F_1,\cdots,F_N$ is an Artinian subscheme.  We will call $s$ the {\it multiplicity} of
the ideal generated by $F_1,\cdots,F_N$.  This definition agrees with that given in (III.3) below for ideals generated by $k$ holomorphic function germs whose common zero-set is of complex codimension $k$.  Let
$$
{\mathcal A}_1=\sum_{j=1}^N{\mathcal O}_{{\mathbb
C}^n,0}\left(dF_j\right)
$$
be the ${\mathcal O}_{{\mathbb C}^n,0}$-submodule of the ${\mathcal
O}_{{\mathbb C}^n,0}$-module ${\mathcal O}_{{\mathbb
C}^n,0}\left(T_{{\mathbb C}^n,0}\right)^*$ of all germs of
holomorphic $(1,0)$-forms on ${\mathbb C}^n$ at $0$.  Take a
sequence of positive integers $q_\nu$ for any positive integer
$\nu$.  By induction on the positive integer $\nu$ we define as
follows the ideals ${\mathcal I}_\nu$ and ${\mathcal J}_\nu$ of
${\mathcal O}_{{\mathbb C}^n,0}$ and the ${\mathcal O}_{{\mathbb
C}^n,0}$-submodule ${\mathcal A}_{\nu+1}$ of the ${\mathcal
O}_{{\mathbb C}^n,0}$-module ${\mathcal O}_{{\mathbb
C}^n,0}\left(T_{{\mathbb C}^n,0}\right)^*$\,.

\medbreak For $\nu\geq 1$ the ideal ${\mathcal J}_\nu$ of ${\mathcal
O}_{{\mathbb C}^n,0}$ is generated over ${\mathcal O}_{{\mathbb
C}^n,0}$ by all holomorphic function-germs $f$ on ${\mathbb C}^n$ at
$0$ satisfying
$$
g_1\wedge\cdots\wedge g_n=f\left(dz_1\wedge\cdots\wedge
dz_n\right)
$$
with $g_1,\cdots,g_n\in{\mathcal A}_\nu$.  The ideal ${\mathcal
I}_\nu$ is defined by the set of all holomorphic function germs
$f$ on ${\mathbb C}^n$ at $0$ so that $f^q\in{\mathcal J}_\nu$ for
some $1\leq q\leq q_\nu$.

\medbreak For $\nu\geq 2$ the ${\mathcal O}_{{\mathbb
C}^n,0}$-submodule ${\mathcal A}_\nu$ of the ${\mathcal O}_{{\mathbb
C}^n,0}$-module ${\mathcal O}_{{\mathbb C}^n,0}\left(T_{{\mathbb
C}^n,0}\right)^*$ is generated by all $df$ for $f\in{\mathcal
I}_{\nu-1}$ and all elements of ${\mathcal A}_{\nu-1}$.

\bigbreak\noindent(III.2) {\it Main Theorem.} There exists an
explicit sequence $\left\{q_\nu\right\}_{\nu\in{\mathbb N}}$ and an
explicit number $m$ depending only on $n$ and $s$ such that
${\mathcal I}_m={\mathcal O}_{{\mathbb C}^n,0}$.

\bigbreak To prepare for the proof of the Main Theorem, we put together some
lemmas about selecting ${\mathbb C}$-linear combinations of $F_1,\cdots,F_N$
to generate ideals with effective multiplicity
and about estimating the multiplicity of Jacobian determinants.

\bigbreak\noindent(III.3) {\it Lemma (on Selection of Linear
Combinations of Holomorphic Functions for Effective Multiplicity).}
Let $0\leq q\leq n$.  Let $f_1,\cdots,f_q$ be holomorphic function
germs on ${\mathbb C}^n$ at the origin whose common zero-set $W_q$
is of pure codimension $q$ in ${\mathbb C}^n$ as a subvariety germ,
with the convention that $W_0={\mathbb C}^n$ and
$\sum_{j=1}^0{\mathcal O}_{{\mathbb C}^n,0}f_j=0$ for the case
$q=0$. Let $m$ be the multiplicity of the ideal
$\sum_{j=1}^q{\mathcal O}_{{\mathbb C}^n,0}f_j$ at the origin in the
sense that
$$\dim_{\mathbb
C}\left({\mathcal O}_{{\mathbb
C}^n,0}\left/\left(\sum_{j=1}^q{\mathcal O}_{{\mathbb
C}^n,0}f_j+\sum_{j=1}^{n-q}{\mathcal O}_{{\mathbb
C}^n,0}L_j\right)\right.\right)=m
$$
for any $n-q$ generic ${\mathbb C}$-linear functions
$L_1,\cdots,L_{n-q}$ on ${\mathbb C}^n$.  Let
$F_j\left(z_1,\cdots,z_n\right)$\ ($1\leq j\leq N$) be holomorphic
function germs on ${\mathbb C}^n$ at the origin which vanish at the
origin. Let $p$ be a positive integer and $A$ be a positive number
such that
$$
\left|z\right|^p\leq A\sum_{j=1}^N\left|F_j(z)\right|\leqno{({\rm
III}.3.1)}
$$
for all $z$ in the domain of definition of
$F_j\left(z_1,\cdots,z_n\right)$ ($1\leq j\leq N$).   Then for
generic choices of complex numbers
$$\left\{c_{j,k}\right\}_{1\leq j\leq n-q,1\leq k\leq N}$$ the
${\mathbb C}$-linear combinations
$$
\tilde F_j=\sum_{k=1}^N c_{j,k}F_k\quad(1\leq j\leq n-q)
$$
of $F_1,\cdots,F_N$ satisfy the property that
$$\dim_{\mathbb
C}\left({\mathcal O}_{{\mathbb
C}^n,0}\left/\left(\sum_{j=1}^q{\mathcal O}_{{\mathbb C}^n,0}f_j+
\sum_{j=1}^{n-q}{\mathcal O}_{{\mathbb C}^n,0}\tilde
F_j\right)\right.\right)\leq mp^{n-q}.
$$
That is, the multiplicity of the ideal generated by
$f_1,\cdots,f_q,\tilde F_1\cdots,\tilde F_{n-q-1}$ is $\leq
mp^{n-q}$ at the origin.

\medbreak\noindent{\sc Proof.} We use induction on $1\leq\nu\leq
n-q$ to show that for generic complex numbers
$$\left\{c_{j,k}\right\}_{1\leq j\leq\nu,1\leq k\leq N}$$
the dimension at the origin of the common zero-set $V_\nu$ of the
$f_1,\cdots,f_q$ and the ${\mathbb C}$-linear combinations
$$
\tilde F_j=\sum_{k=1}^N c_{j,k}F_k\quad(1\leq j\leq\nu)
$$
of $F_1,\cdots,F_N$ is precisely $n-q-\nu$ and the multiplicity of
the ideal
$$
\sum_{j=1}^q{\mathcal O}_{{\mathbb
C}^n,0}f_j+\sum_{j=1}^\nu{\mathcal O}_{{\mathbb C}^n,0}\tilde F_j
$$
is no more than $mp^\nu$.

\medbreak We introduce the case of $\nu=0$ and the convention that
$\sum_{j=1}^0{\mathcal O}_{{\mathbb C}^n,0}\tilde F_j=0$ for the
case $\nu=0$.  With this convention, we start out our induction
assumption with the case $\nu=0$ which is trivially true.

\medbreak Suppose the induction process has been carried out for
some $0\leq\nu<n-q$ and we would like to verify it for the next
step when $\nu$ is replaced by $\nu+1$.  We now already have
$\tilde F_1,\cdots,\tilde F_\nu$. Let
$$
{\mathcal I}_\nu=\sum_{j=1}^q{\mathcal O}_{{\mathbb
C}^n}f_j+\sum_{j=1}^\nu{\mathcal O}_{{\mathbb C}^n}\tilde F_j.
$$
The zero-set of ${\mathcal I}_\nu$ is the subvariety $V_\nu$ of pure
dimension $n-q-\nu$.  Let $E_\nu$ be a generic linear subspace of
${\mathbb C}^n$ of codimension $n-q-\nu-1$ defined by $n-q-\nu-1$
generic linear functions $G_1,\cdots,G_{n-q-\nu+1}$ so that the
subvariety $V_\nu\cap E_\nu$ is of pure dimension $1$. Let
$$
{\mathcal J}_\nu=\sum_{j=1}^{n-q-\nu-1}{\mathcal O}_{{\mathbb
C}^n}G_j+\sum_{j=1}^q{\mathcal O}_{{\mathbb
C}^n}f_j+\sum_{j=1}^\nu{\mathcal O}_{{\mathbb C}^n}\tilde F_j.
$$
Let
$$
{\mathcal J}_\nu=\bigcap_{\lambda=1}^\Lambda{\mathcal L}_\lambda
$$
be the primary decomposition of the ideal sheaf ${\mathcal J}_\nu$.
Note that, since the zero-set of ${\mathcal J}_\nu$ is of pure
complex dimension $1$ and ${\mathcal J}_\nu$ is generated by $n-1$
holomorphic functions
$$G_1,\cdots,G_{n-q-\nu+1},f_1,\cdots,f_q,\tilde F_1,\cdots,\tilde
F_\nu,$$ it follows that all the associated prime ideals of
${\mathcal J}_\nu$ are isolated and none are embedded [ZS60, p.397,
Theorem 2].

\medbreak Let $C_\lambda$ be the complex curve-germ which is the
zero-set of the ideal sheaf ${\mathcal L}_\lambda$.  Let
$\mu_\lambda$ be the multiplicity of the curve $C_\lambda$ at the
origin.  Let $\hat\mu_\lambda$ be the multiplicity of the ideal
sheaf ${\mathcal L}_\lambda$ at a generic point $Q\in C_\lambda$,
which can be characterized as the dimension over ${\mathbb C}$ of
$$
{\mathcal O}_{{\mathbb C}^n,Q}\left/\left(\left({\mathcal
L}_\lambda\right)_Q+{\mathcal O}_{{\mathbb C}^n,Q}L\right)\right.,
$$
where $L$ is a generic polynomial of degree $1$ on ${\mathbb C}^n$
vanishing at $Q$ and $\left({\mathcal L}_\lambda\right)_Q$ is the
stalk of the ideal sheaf ${\mathcal L}_\lambda$ at the point $Q$.

\medbreak Without loss of generality we can assume that the
coordinates $\left(z_1,\cdots,z_n\right)$ of ${\mathbb C}^n$ are
chosen so that $C_\lambda$ is defined by
$$
\left\{\begin{matrix}z_1=\zeta^{\mu_\lambda},\hfill\smallskip\cr
z_j=g_{\lambda,j}\left(\zeta\right){\rm\ for\ }2\leq j\leq
n\cr\end{matrix}\right.
$$
for $\zeta$ in some open neighborhood of the origin in ${\mathbb
C}$, where the initial term of $g_{\lambda,j}\left(\zeta\right)$ is
a nonzero complex number times $\zeta^{N_{\lambda,j}}$ for some
$N_{\lambda,j}\geq\mu_\lambda$ for $2\leq j\leq n$. Let
$\pi_\lambda:\tilde C_\lambda\to C_\lambda$ be the normalization of
$C_\lambda$ defined by
$$\pi_\lambda:\zeta\mapsto z=\left(\zeta^{\mu_\lambda},
g_{\lambda,2}\left(\zeta\right),\cdots,
g_{\lambda,n}\left(\zeta\right)\right),$$ where $\tilde C_\lambda$
is an open neighborhood of $0$ in $\mathbb C$ with $\zeta$ as
coordinate. The pullback $\pi_\lambda^*{\mathfrak m}_{C_\lambda,0}$
to $\tilde C_\lambda$ of the maximum ideal ${\mathfrak
m}_{C_\lambda,0}$ of $C_\lambda$ at the origin is generated by
$\zeta^{\mu_\lambda},
g_{\lambda,2}\left(\zeta\right),\cdots,g_{\lambda,n}\left(\zeta\right)$.
Since $\pi_\lambda^*{\mathfrak m}_{C_\lambda,0}$ is a principal
ideal, it must be generated by $\zeta^{\mu_\lambda}$.

\medbreak The inequality (III.3.1), when pulled back by
$\pi_\lambda$, becomes
$$
\left|\zeta\right|^{p\mu_\lambda}\leq
A_\lambda\sum_{j=1}^N\left|\left(F_j\circ\pi_\lambda\left(\zeta\right)\right)
\right|\leqno{({\rm III}.3.2)_\nu}
$$
for $1\leq\lambda\leq\Lambda$, where $A_\lambda$ is a positive
number.  Take a generic point
$$\left(c_{\nu+1,1},\cdots,c_{\nu+1,N}\right)\in{\mathbb
C}^N
$$
and let
$$
\tilde F_{\nu+1}=\sum_{k=1}^N c_{\nu+1,k}F_k.
$$
By $({\rm III}.3.2)_\nu$, for each $1\leq\lambda\leq\Lambda$ the
vanishing order of $\left(\tilde
F_{\nu+1}\circ\pi_\lambda\right)\left(\zeta\right)$ at $\zeta=0$ is
some number $\tilde\mu_\lambda$ which is no more than
$p\mu_\lambda$. For a small generic nonzero $\eta\in{\mathbb C}$ the
number of zeros of $\eta+\left(\tilde
F_{\nu+1}\circ\pi_\lambda\right)\left(\zeta\right)$ on $\tilde
C_\lambda$ is precisely $\tilde\mu_\lambda$ with multiplicity $1$
for each $1\leq\lambda\leq\Lambda$.  Since the map
$\pi_\lambda:\tilde C_\lambda\to C_\lambda$ is one-to-one, it
follows that for any small generic nonzero $\eta\in{\mathbb C}$ the
number of zeroes of $\eta+\tilde F_{\nu+1}$ on $C_\lambda$ is
precisely $\tilde\mu_\lambda$ with multiplicity $1$.

\medbreak Since the multiplicity of the ideal sheaf ${\mathcal
L}_\lambda$ at a generic point $Q\in C_\lambda$ is
$\hat\mu_\lambda$, it follows from $\tilde\mu_\lambda\leq
p\mu_\lambda$ that the dimension over ${\mathbb C}$ of the vector
space
$${\mathcal O}_{{\mathbb C}^n,0}\left/\left(
\sum_{j=1}^{n-q-\nu-1}{\mathcal O}_{{\mathbb
C}^n,0}G_j+\sum_{j=1}^q{\mathcal O}_{{\mathbb C}^n,0}f_j
+\sum_{j=1}^{\nu+1}{\mathcal O}_{{\mathbb C}^n,0}\tilde
F_j\right)\right.
$$
is no more than $
p\sum_{\lambda=1}^\Lambda\mu_\lambda\hat\mu_\lambda$.  By
induction hypothesis the multiplicity of
$$
{\mathcal I}_\nu=\sum_{j=1}^q{\mathcal O}_{{\mathbb C}^n}f_j
+\sum_{j=1}^\nu{\mathcal O}_{{\mathbb C}^n}\tilde F_j
$$
is no more than $mp^\nu$ at the origin.  The multiplicity of
${\mathcal J}_\nu$ at the origin, which can be computed from
${\mathcal I}_\nu$ by adding generic ${\mathbb C}$-linear functions
$G_1,\cdots,G_{n-q-\nu-1}$ on ${\mathbb C}^n$, is also no more than
$mp^\nu$.  We can compute the multiplicity of ${\mathcal J}_\nu$ at
the origin by adding to ${\mathcal J}_\nu$ a generic ${\mathbb
C}$-linear function $L$ on ${\mathbb C}^n$ and considering the sum
of the multiplicities at points of intersection of the zero-set with
$L+\eta$ for some small generic $\eta\in{\mathbb C}$.  From the
decomposition ${\mathcal J}_\nu=\bigcap_{\lambda=1}^\Lambda{\mathcal
L}_\lambda$ and the multiplicity $\hat\mu_\lambda$ of ${\mathcal
L}_\lambda$ at the origin we conclude that
$$
\sum_{\lambda=1}^\Lambda\mu_\lambda\hat\mu_\lambda\leq mp^\nu.$$
Thus
$$\dim_{\mathbb
C}\left({\mathcal O}_{{\mathbb
C}^n,0}\left/\left(\sum_{j=1}^{n-q-\nu+1}{\mathcal O}_{{\mathbb
C}^n,0}G_j+\sum_{j=1}^q{\mathcal O}_{{\mathbb C}^n,0}f_j
+\sum_{j=1}^{\nu+1}{\mathcal O}_{{\mathbb C}^n,0}\tilde
F_j\right)\right.\right)\leq m p^{\nu+1}.
$$
Since $G_1,\cdots,G_{n-q-\nu+1}$ are generic linear functions on
${\mathbb C}^n$, it follows that the multiplicity of the ideal
$$
\sum_{j=1}^q{\mathcal O}_{{\mathbb C}^n,0}f_j
+\sum_{j=1}^{\nu+1}{\mathcal O}_{{\mathbb C}^n,0}\tilde F_j
$$
at the origin is no more than $mp^{\nu+1}$.  This finishes the
induction process. Q.E.D.

\bigbreak\noindent(III.4) {\it Corollary.} Let
$F_j\left(z_1,\cdots,z_n\right)$\ ($1\leq j\leq N$) be holomorphic
function germs on ${\mathbb C}^n$ at the origin which vanish at the
origin. Let $p$ be a positive integer and $A$ be a positive number
such that
$$
\left|z\right|^p\leq A\sum_{j=1}^N\left|F_j(z)\right|
$$
for all $z$ in the domain of definition of
$F_j\left(z_1,\cdots,z_n\right)$ ($1\leq j\leq N$).   Then for
generic choices of complex numbers
$$\left\{c_{j,k}\right\}_{1\leq j\leq n,1\leq k\leq N}$$ the
${\mathbb C}$-linear combinations
$$
\tilde F_j=\sum_{k=1}^N c_{j,k}F_k\quad(1\leq j\leq n)
$$
of $F_1,\cdots,F_N$ satisfy the property that
$$\dim_{\mathbb
C}\left({\mathcal O}_{{\mathbb C}^n,0}\left/\sum_{j=1}^n{\mathcal
O}_{{\mathbb C}^n,0}\tilde F_j\right.\right)\leq p^n.
$$

\medbreak\noindent{\it Proof.}  Introduce one more complex
variable $w$ and consider $F_j$ as a holomorphic function germ on
${\mathbb C}^{n+1}$ at $0$ in the variables $z_1,\cdots,z_n,w$
though it is independent of the variable $w$.  Add the function
$w$ to the functions $F_1,\cdots,F_N$. Let $f_1=w$ and apply Lemma (III.3) on
Selection of Linear Combinations of Holomorphic Functions for
Effective Multiplicity to the case $m=1$ with ${\mathbb C}^n$
replaced by ${\mathbb C}^{n+1}$.  Q.E.D.

\bigbreak\noindent(III.5) {\it Lemma (Multiplicity Estimate for Jacobian Determinant).}
Let $g_1,\cdots,g_n$ be holomorphic function germs on ${\mathbb C}^n$ at the origin such that
$$\dim_{\mathbb
C}\left({\mathcal O}_{{\mathbb C}^n,0}\left/\sum_{j=1}^n{\mathcal
O}_{{\mathbb C}^n,0}g_j\right.\right)\leq m.$$  Let $dg_1\wedge\cdots\wedge dg_n=f\left(dz_1,\cdots,z_n\right)$.  Then the multiplicity of $f$ at the origin is $\leq m$.

\medbreak\noindent{\it Proof.} We can find a connected open neighborhood $U$ of $0$ in ${\mathbb C}^n$ and an open ball neighborhood $W$ of $0$ in ${\mathbb C}^n$ such that the map $\pi: U\to W$ defined by $g_1,\cdots,g_n$ is a proper holomorphic map.  This is possible, because the common zero-set of $g_1,\cdots,g_n$ consists only of the origin in a sufficiently small neighborhood of the origin in ${\mathbb C}^n$.  The number of sheets in the analytic cover map $\pi: U\to V$ is $\leq m$.  Let $Y$ be the divisor of $f$ in $U$ and $Z$ be the image of $Y$ in $W$.  Let $Z_0$ be the set of regular points of (the reduction of $Z$).  Let $L$ be a generic complex line in the target space ${\mathbb C}^n$ such that $L\cap Z$ is a single point $P$ in $Z_0$ and $L$ intersects $Z_0$ transversally at $P$.  For a sufficiently small neighborhood $D$ of $P$ in $W$ the map $U\cap\pi^{-1}(D)\to D$ induced by $\pi$ is just a cyclic branched cover on each topological component of $U\cap\pi^{-1}(D)$. Thus the multiplicity of the intersection of the regular curve $\pi^{-1}(L)$ and the divisor $Y$ is no more than the number of sheets of $\pi:U\to W$.  Since the line $L$ is generic, it follows that the multiplicity of the divisor $Y$ is more than $m$. Q.E.D.

\bigbreak\noindent(III.6) {\it Preparatory Remarks on Proof of Main Theorem.}
We now start the setting for the proof of the Main Theorem (III.2). Let $F_1,\cdots,F_N$ be
holomorphic function germs on ${\mathbb C}^n$ at the origin whose common zero-set is the origin. Let
$q$ be a positive integer. Assume that, for some positive number
$A$,
$$
\left|z\right|^q\leq A\sum_{j=1}^N\left|F_j(z)\right|\leqno{({\rm
III}.6.1)}
$$
for all $z$ in the domain of definition of
$F_j\left(z_1,\cdots,z_n\right)$ ($1\leq j\leq N$).   Because of the discussion in (II.1), for the case of special domains we need only consider multipliers which are holomorphic and we need only consider vector-multipliers which are holomorphic $(1,0)$-forms.  Though $F_1,\cdots, F_N$ are not multipliers, their differentials $dF_1,\cdots,dF_N$ are vector-multipliers and, in order to form Jacobian determinants to generate multipliers, we can also use $\ell$ ${\mathbb C}$-linear combinations of $F_1,\cdots,F_N$ and $n-\ell$ multipliers for $0\leq\ell\leq n$ to form a Jacobian determinant which will then be a multiplier.  We will refer to any ${\mathbb C}$-linear combination of $F_1,\cdots,F_N$ and multipliers as {\it pre-multipliers} so that the $(1,0)$-differential of a pre-multiplier is a vector-multiplier.  Note that the pre-multipliers form a ${\mathbb C}$-vector space but do not form an ideal.  The product of a multiplier and a holomorphic function germ is again a multiplier, but the product of a pre-multiplier and a holomorphic function germ in general is not a pre-multiplier.  In our proof of the Main Theorem (III.2) we will not use vector-multipliers, because we will directly form the Jacobian determinants of the holomorphic pre-multipliers to generate new multipliers to bypass the process of forming vector-multipliers by differentiation and then using Cramer's rule.

\medbreak In order not to be encumbered by complicated expressions of constants,
we will not explicitly keep track of the various effective bounds occurring in the proof.  We introduce the
following terminology.  A multiplier is called {\it effectively constructed} if there is an effective upper bound for its multiplicity and there is an effective positive lower bound for its assigned order of subellipticity.  Effective means some explicit function of the multiplicity of the ideal generated by the pre-multipliers $F_1,\cdots,F_N$, which means an explicit function of $q$ given in (III.6.1).  The goal is to show that the function-germ with constant value $1$ can be effectively constructed.

\medbreak To make the argument more transparent and to minimize notational clutters, we start out with the proof of the simple case where $n=2$.

\bigbreak\noindent(III.7) {\it Proof of Main Theorem for Dimension Two.}  We now assume that $n=2$ and we have holomorphic function germs $F_1,\cdots,F_N$ on ${\mathbb C}^2$ at the origin whose zero-set is the origin of ${\mathbb C}^n$.  The multiplicity of the ideal generated by $F_1,\cdots,F_N$ is the number used to express effectiveness.  By applying Corollary (III.4) and (III.5) to get two ${\mathbb C}$-linear combinations of $F_1,\cdots,F_N$ and form their Jacobian determinant, we get an effectively constructed multiplier $\tilde h_2\left(z_1,z_2\right)$ at the origin, which vanishes at the origin.  Because the multiplicity of $\tilde h_2\left(z_1,z_2\right)$ is effectively bounded at the origin, by replacing $\tilde h_2\left(z_1,z_2\right)$ by the product of holomorphic function germs defining the branch germs of the reduction of the subspace defined by $\tilde h_2\left(z_1,z_2\right)$, we can assume without loss of generality that the subspace germ $C_2$ defined by $\tilde h_2$ is a reduced curve germ in ${\mathbb C}^2$ at the origin with effectively bounded multiplicity.  Note that in general the curve germ $C_2$ is not irreducible, though $C_2$ is a reduced curve.  A reduced curve means that its structure sheaf does not contain any nonzero nilpotent elements.  For example, it means that $h_2$ does not vanish to order higher than one at any regular point of $C_2$.

\medbreak Now the ideal generated by the functions $\frac{\partial F_i}{\partial z_j}$ for $1\leq i\leq N$ and $j=1,2$ has effectively bounded multiplicity at the origin, because by Proposition(A.2) in Appendix A, for each fixed $1\leq i\leq N$, the function germ
  $\left(F_i\right)^3$ at the origin belongs to the ideal generated by $\frac{\partial F_i}{\partial z_j}$ for $j=1,2$.  We consider the pre-multiplier $h_1=\sum_{j=1}^N c_jF_j$ with generic $c_j\in{\mathbb C}$ for $1\leq j\leq N$ and consider a new generic linear coordinate system $\left(w_1, w_2\right)$ which is related to $\left(z_1,z_2\right)$ by $w_i=\sum_{j=1}^2 b_{ij}z_j$ with generic $b_{ij}\in{\mathbb C}$ for $1\leq i,j\leq 2$.  By (III.3) we can find generic $c_j\in{\mathbb C}$ for $1\leq j\leq N$ and generic $b_{ij}\in{\mathbb C}$ for $1\leq i,j\leq 2$ such that
  \begin{itemize}\item[(i)] the ideal generated by $h_1$ and $\tilde h_2$ has effectively bounded multiplicity at the origin,
  \item[(ii)] the ideal generated by $\frac{\partial h_1}{\partial w_1}$ and $\tilde h_2$ has effectively bounded multiplicity at the origin, where the partial derivative $\frac{\partial h_1}{\partial w_1}$ is computed with $w_2$ being kept constant,
\item[(iii)] the projection $P\mapsto g(P)$ makes $C_2$ an analytic cover over ${\mathbb C}$ locally at the origin as germs, and

\item[(iv)] the projection $\left(w_1,w_2\right)\mapsto\left(h_1,w_2\right)$ makes ${\mathbb C}^2$ an analytic cover over ${\mathbb C}^2$ locally at the origin as germs.
\end{itemize}
Without loss of generality we can assume that the coordinate system $\left(w_1,w_2\right)$ is just the coordinate system $\left(z_1, z_2\right)$.  Note that $h_1$ is only a pre-multiplier and in general may not be a multiplier.  The function germ $\frac{\partial h_1}{\partial z_1}$ is in general not a multiplier and not even a pre-multiplier.

\medbreak Consider the image $\hat C_2$ of $C_2$ under the projection $\left(z_1,z_2\right)\mapsto\left(h_1,z_2\right)$ and let $$h_2=z_2^\lambda+\sum_{j=0}^{\lambda-1}a_j\left(h_1\right)z_2^j$$ be the Weierstrass polynomial in ${\mathbb C}^2$ with coordinates $\left(h_1,z_2\right)$ whose vanishing defines the curve-germ $\hat C_2$ at the origin in ${\mathbb C}^2$.  This is possible, because the projection $P\mapsto h_1(P)$ makes $C_2$ an analytic cover over ${\mathbb C}$ locally at the origin as germs.  When regarded as a function-germ in the variables $\left(z_1,z_2\right)$ the function-germ $h_2$ contains $\tilde h_2$ as a factor, because the inverse image of $\hat C_2$ under the projection $\left(z_1,z_2\right)\mapsto\left(h_1,z_2\right)$ contains $C_2$ and $C_2$ is a reduced curve.  Since $\tilde h_2$ is a multiplier, it follows that $h_2$ is also a multiplier and is, in fact, an effectively constructed multiplier.  The multipliers in the effective procedure presented here and also in (III.8) and (III.9) are all effectively constructed multipliers (unless explicitly pointed out otherwise) and we will drop the description ``effectively constructed'' when we mention these multipliers here and in (III.8) and (III.9).  Sometimes, to highlight certain aspects of effectiveness, we may mention ``the assigned order of subellipticity having an effective positive lower bound'' or ``the multiplicity having an effective upper bound'' in conjunction with such multipliers, though according to the convention given here such multipliers are all effective constructed unless explicitly pointed out otherwise.

\medbreak Since the ideal generated by $\frac{\partial h_1}{\partial z_1}$ and $\tilde h_2$ has effectively bounded multiplicity at the origin and since $h_1$ vanishes at the origin, it follows that, for some effectively bounded positive integer $s$, the function germ $\left(h_1\right)^s$ belongs to the ideal generated by $\frac{\partial h_1}{\partial z_1}$ and $\tilde h_2$.  In particular, $$\left|\left(h_1\right)^s\right|\stackrel{<}{\sim}\left|\frac{\partial h_1}{\partial z_1}\right|+\left|\tilde h_2\right|.\leqno{({\rm III}.7.1)}$$
  Here and for the rest of this note the symbol $\stackrel{<}{\sim}$ means ``less than some constant times'' and is being used to avoid introducing new symbols for constants.  We now form $dh_1\wedge dh_2$ and get
$$
\displaylines{dh_1\wedge dh_2=dh_1\wedge\left(\lambda z_2^{\lambda-1}dz_2+\sum_{j=1}^{\lambda-1}j a_j\left(h_1\right)z_2^{j-1}dz_2+\sum_{j=0}^{\lambda-1}a_j^\prime\left(h_1\right)z_2^jdh_1\right)\cr
=dh_1\wedge\left(\lambda z_2^{\lambda-1}dz_2+\sum_{j=1}^{\lambda-1}j a_j\left(h_1\right)z_2^{j-1}dz_2\right)\cr
=\left(\frac{\partial h_1}{\partial z_1}dz_1+\frac{\partial h_1}{\partial z_2}dz_2\right)\wedge\left(\lambda z_2^{\lambda-1}dz_2+\sum_{j=1}^{\lambda-1}j a_j\left(h_1\right)z_2^{j-1}dz_2\right)\cr
=\frac{\partial h_1}{\partial z_1}\left(\lambda z_2^{\lambda-1}+\sum_{j=1}^{\lambda-1}j a_j\left(h_1\right)z_2^{j-1}\right)dz_1\wedge dz_2,\cr}
$$
where $a^\prime_j\left(h_1\right)$ is the derivative of $a_j\left(h_1\right)$ as a function of $h_1$.  Since $h_1$ is a pre-multiplier, the coefficient of $dz_1\wedge dz_2$ in $dh_1\wedge dh_2$ is a multiplier.  Thus
$$\frac{\partial h_1}{\partial z_1}\left(\lambda z_2^{\lambda-1}+\sum_{j=1}^{\lambda-1}j a_j\left(h_1\right)z_2^{j-1}\right)$$
is a multiplier.  Since $\tilde h_2$ is a multiplier, it follows that
$$\tilde h_2\left(\lambda z_2^{\lambda-1}+\sum_{j=1}^{\lambda-1}j a_j\left(h_1\right)z_2^{j-1}\right)$$
is a multiplier.  From (III.7.1) it follows that
$$\displaylines{\left|\left(h_1\right)^s\left(\lambda z_2^{\lambda-1}+\sum_{j=1}^{\lambda-1}j a_j\left(h_1\right)z_2^{j-1}\right)\right|\cr\stackrel{<}{\sim}\left|\frac{\partial h_1}{\partial z_1}\left(\lambda z_2^{\lambda-1}+\sum_{j=1}^{\lambda-1}j a_j\left(h_1\right)z_2^{j-1}\right)\right|+\left|
\tilde h_2\left(\lambda z_2^{\lambda-1}+\sum_{j=1}^{\lambda-1}j a_j\left(h_1\right)z_2^{j-1}\right)\right|.\cr}$$
Hence
$$\left(h_1\right)^s\left(\lambda z_2^{\lambda-1}+\sum_{j=1}^{\lambda-1}j a_j\left(h_1\right)z_2^{j-1}\right)$$
is a multiplier.  Let $h_2^{(0)}=h_2$ and
$$
h^{(1)}_2=\left(h_1\right)^s\left(\lambda z_2^{\lambda-1}+\sum_{j=1}^{\lambda-1}j a_j\left(h_1\right)z_2^{j-1}\right)
$$
and for $1<\nu\leq\lambda$ define
$$
h^{(\nu)}_2=\left(h_1\right)^{s\nu}\left(\frac{\lambda!}{\left(\lambda-\nu\right)!} z_2^{\lambda-\nu}+\sum_{j=\nu}^{\lambda-1}\frac{j!}{\left(j-\nu\right)!} a_j\left(h_1\right)z_2^{j-\nu}\right).
$$
We are going to verify by induction on $\nu$ that $h^{(\nu)}_2$ is a multiplier.  We know that both $h^{(0)}_2$ and $h^{(1)}_2$ are multipliers.  Assume that we have already verified that
$h^{(0)}_2,\cdots,h^{(\nu-1)}_2$ are multipliers.  Then $dh_1\wedge dh^{(\nu-1)}_2$ is equal to
$$
\displaylines{dh_1\wedge\; d\left(
\left(h_1\right)^{s\left(\nu-1\right)}\left(\frac{\lambda!}{\left(\lambda-\nu+1\right)!} z_2^{\lambda-\nu+1}+\sum_{j=\nu-1}^{\lambda-1}\frac{j!}{\left(j-\nu+1\right)!} a_j\left(h_1\right)z_2^{j-\nu+1}\right)\right)\cr
=dh_1\wedge \left(
\left(h_1\right)^{s\left(\nu-1\right)}\left(\frac{\lambda!}{\left(\lambda-\nu\right)!} z_2^{\lambda-\nu}+\sum_{j=\nu}^{\lambda-1}\frac{j!}{\left(j-\nu\right)!} a_j\left(h_1\right)z_2^{j-\nu}\right)dz_2\right)\cr
=\frac{\partial h_1}{\partial z_1}\left(
\left(h_1\right)^{s\left(\nu-1\right)}\left(\frac{\lambda!}{\left(\lambda-\nu\right)!} z_2^{\lambda-\nu}+\sum_{j=\nu}^{\lambda-1}\frac{j!}{\left(j-\nu\right)!} a_j\left(h_1\right)z_2^{j-\nu}\right)\right)dz_1\wedge dz_2\cr
}
$$
Since the coefficient of $dz_1\wedge dz_2$ in $dh_1\wedge dh^{(\nu-1)}_2$ is a multiplier, it follows that
$$\frac{\partial h_1}{\partial z_1}\left(
\left(h_1\right)^{s\left(\nu-1\right)}\left(\frac{\lambda!}{\left(\lambda-\nu\right)!} z_2^{\lambda-\nu}+\sum_{j=\nu}^{\lambda-1}\frac{j!}{\left(j-\nu\right)!} a_j\left(h_1\right)z_2^{j-\nu}\right)\right)$$
is a multiplier. Since $\tilde h_2$ is a multiplier, it follows that
$$\tilde h_2\left(
\left(h_1\right)^{s\left(\nu-1\right)}\left(\frac{\lambda!}{\left(\lambda-\nu\right)!} z_2^{\lambda-\nu}+\sum_{j=\nu}^{\lambda-1}\frac{j!}{\left(j-\nu\right)!} a_j\left(h_1\right)z_2^{j-\nu}\right)\right)$$
is a multiplier.  From (III.7.1) it follows that
$$\displaylines{\left|\left(h_1\right)^s\left(
\left(h_1\right)^{s\left(\nu-1\right)}\left(\frac{\lambda!}{\left(\lambda-\nu\right)!} z_2^{\lambda-\nu}+\sum_{j=\nu}^{\lambda-1}\frac{j!}{\left(j-\nu\right)!} a_j\left(h_1\right)z_2^{j-\nu}\right)\right)\right|\cr\stackrel{<}{\sim}\left|\frac{\partial h_1}{\partial z_1}\left(
\left(h_1\right)^{s\left(\nu-1\right)}\left(\frac{\lambda!}{\left(\lambda-\nu\right)!} z_2^{\lambda-\nu}+\sum_{j=\nu}^{\lambda-1}\frac{j!}{\left(j-\nu\right)!} a_j\left(h_1\right)z_2^{j-\nu}\right)\right)\right|\cr+\left|
\tilde h_2\left(
\left(h_1\right)^{s\left(\nu-1\right)}\left(\frac{\lambda!}{\left(\lambda-\nu\right)!} z_2^{\lambda-\nu}+\sum_{j=\nu}^{\lambda-1}\frac{j!}{\left(j-\nu\right)!} a_j\left(h_1\right)z_2^{j-\nu}\right)\right)\right|.\cr}$$
Hence
$$\left(h_1\right)^s\left(
\left(h_1\right)^{s\left(\nu-1\right)}\left(\frac{\lambda!}{\left(\lambda-\nu\right)!} z_2^{\lambda-\nu}+\sum_{j=\nu}^{\lambda-1}\frac{j!}{\left(j-\nu\right)!} a_j\left(h_1\right)z_2^{j-\nu}\right)\right)$$
is a multiplier and $h^{(\nu)}_2$ is a multiplier.  When $\nu=\lambda$, we end up with
$$
h_2^{(\lambda)}=\left(h_1\right)^{s\lambda}\lambda!
$$
being a multiplier.

\medbreak Note that this step of forming Jacobian determinants $\lambda$ times to construct $h_2^{(\lambda)}$ from $\tilde h_2$ is the step of differentiating a multiplier as many times as its multiplicity to form a new multiplier, which is referred to at the end of the Introduction of this note.   Also note that though this step only requires $h_2$ to be a pre-multiplier, yet $\tilde h_2$ has to be a multiplier instead of just a pre-multiplier, otherwise we cannot conclude that $h_2$ is a pre-multiplier, because in general the set of all pre-multipliers do not form an ideal.

\medbreak Since the multiplicity of the ideal generated by $h_1$ and $\tilde h_2$ is effectively bounded, there exists some positive integer $\sigma$ which is effectively bounded such that $z_1^\sigma$ and $z_2^\sigma$ both belong to the ideal generated by $h_2^{(\lambda)}$ and $\tilde h_2$.  Hence both $z_1^\sigma$ and $z_2^\sigma$ are multipliers.  We take the $\sigma$-th roots of both $z_1^\sigma$ and $z_2^\sigma$ to produce multipliers $z_1$ and $z_2$.  We finally form the Jacobian determinant of the two holomorphic function germs $z_1$ and $z_2$ to conclude that Kohn's algorithm effectively produces the function $F\equiv 1$ as a multiplier.  This finishes the proof of the effective termination of Kohn's algorithm in complex dimension $2$.

\bigbreak\noindent(III.8) {\it Construction of a New Multiplier in Higher Dimensional Case by Fiberwise Differentiating a Given Multiplier as Many Times as its Multiplicity.} We now look at the higher dimensional case.  As a preparation for the proof of the Main Theorem for the higher dimensional case, we do the argument here for the construction a new multiplier in higher dimensional case by fiberwise differentiating a given multiplier as many times as its multiplicity.  The argument is the same as the $2$-dimensional case with corresponding modifications in notations.

\medbreak We have holomorphic function germs $F_1,\cdots,F_N$ on ${\mathbb C}^n$ at the origin which generate an ideal of multiplicity $q$ whose zero-set is the origin of ${\mathbb C}^n$.  By applying Corollary (III.4) and (III.5) to get $n$ ${\mathbb C}$-linear combinations of $F_1,\cdots,F_N$ and form their Jacobian determinant, we get an effectively constructed multiplier $\tilde h_n\left(z_1,\cdots,z_n\right)$ at the origin, which vanishes at the origin.  The divisor of $\tilde h_n$ is a subspace germ $V_n$ of codimension $1$ in ${\mathbb C}^n$ at the origin with effectively bounded multiplicity.
Because $V_n$ has effectively bounded multiplicity, by replacing $\tilde h_n$ by the product of holomorphic function germs defining the branch germs of the reduction of $V_n$, we can assume without loss of generality that $V_n$ is a reduced hypersurface germ in ${\mathbb C}^n$ at the origin with effectively bounded multiplicity.   Again this does not mean that $V_n$ is irreducible.  It only means that the divisor of $\tilde h_n$ has coefficient $1$ for every one of its irreducible components.

\medbreak \medbreak By Proposition (A.3) in Appendix A, the ideal generated by
$$
\frac{\partial\left(F_{i_1},\cdots,F_{i_{n-1}}\right)}{\partial\left(z_{j_1},\cdots,z_{j_{n-1}}\right)}
$$
for $1\leq i_1<\cdots<i_{n-1}\leq N$ and $1\leq j_1<\cdots<j_{n-1}\leq n$ contains an effective power of the maximum ideal ${\mathfrak m}_{{\mathbb C}^n,0}$ of ${\mathbb C}^n$ at the origin. Just like the argument given in the $2$-dimensional case in (III.7), after a generic ${\mathbb C}$-linear coordinate change and after taking $n-1$ generic ${\mathbb C}$-linear combinations $h_1,\cdots,h_{n-1}$ of $F_1,\cdots,F_N$ we have the following situation. \begin{itemize}\item[(i)] The ideal generated by $h_1,\cdots,h_{n-1}$ and $\tilde h_n$ has effectively bounded multiplicity at the origin,
  \item[(ii)] The ideal generated by $\tilde h_n$ and  $$\frac{\partial\left(h_1,\cdots,h_{n-1}\right)}{\partial\left(z_1,\cdots,z_{n-1}\right)}$$ has effectively bounded multiplicity at the origin.
\item[(iii)] The projection $P\mapsto\left(h_1(P),\cdots,h_{n-1}(P)\right)$ makes $V_n$ an analytic cover over ${\mathbb C}^{n-1}$ locally at the origin as germs.
\item[(iv)] The projection $\left(z_1,\cdots,z_n\right)\mapsto\left(h_1,\cdots,h_{n-1},z_n\right)$ makes ${\mathbb C}^n$ an analytic cover over ${\mathbb C}^n$ locally at the origin as germs.
\end{itemize}
Consider the image $\hat V_n$ of $V_n$ under the projection $\left(z_1,\cdots,z_n\right)\mapsto\left(h_1,\cdots,h_{n-1},z_n\right)$ and let $$h_n=z_n^\lambda+\sum_{j=0}^{\lambda-1}a_j\left(h_1,\cdots,h_{n-1}\right)z_n^j$$ be the Weierstrass polynomial in the target space ${\mathbb C}^n$ with coordinates $\left(h_1,\cdots,h_{n-1},z_n\right)$ whose vanishing defines the subspace germ $\hat V_n$ at the origin in ${\mathbb C}^n$.  This is possible, because the projection $P\mapsto\left(h_1(P),\cdots,h_{n-1}(P)\right)$ makes $V_n$ an analytic cover over ${\mathbb C}^{n-1}$ locally at the origin as germs.  When regarded as a function-germ in the variables $\left(z_1,\cdots,z_n\right)$ the function-germ $h_n$ contains $\tilde h_n$ as a factor, because the inverse image of $\hat V_n$ under the projection $\left(z_1,\cdots,z_n\right)\mapsto\left(h_1,\cdots,h_{n-1},z_n\right)$ contains $V_n$ and because $V_n$ which is defined by $\tilde h_n$ is reduced.  Since $\tilde h_n$ is a multiplier, it follows that $h_n$ is also a multiplier.  Since the ideal generated by $\tilde h_n$ and  $$\frac{\partial\left(h_1,\cdots,h_{n-1}\right)}{\partial\left(z_1,\cdots,z_{n-1}\right)}$$ has effectively bounded multiplicity at the origin and since $h_1\cdots,h_{n-1}$ all vanish at the origin, it follows that there exists some polynomial $p\left(h_1,\cdots,h_{n-1}\right)$ such that
\begin{itemize}\item[(i)] the ideal generated by $p\left(h_1,\cdots,h_{n-1}\right)$ and $\tilde h_n$ has effectively bounded multiplicity at the origin, and\item[(ii)] $p\left(h_1,\cdots,h_{n-1}\right)$ belongs to the ideal generated by $\tilde h_n$ and $$\frac{\partial\left(h_1,\cdots,h_{n-1}\right)}{\partial\left(z_1,\cdots,z_{n-1}\right)}.$$\end{itemize}  In particular, $$\left|p\left(h_1,\cdots,h_{n-1}\right)\right|\stackrel{<}{\sim}
\left|\frac{\partial\left(h_1,\cdots,h_{n-1}\right)}{\partial\left(z_1,\cdots,z_{n-1}\right)}\right|+\left|\tilde h_n\right|.\leqno{({\rm III}.8.1)}$$
One way to obtain the polynomial $p\left(h_1,\cdots,h_{n-1}\right)$ is to use the direct image of the ideal generated by generated by $\tilde h_n$ and $$\frac{\partial\left(h_1,\cdots,h_{n-1}\right)}{\partial\left(z_1,\cdots,z_{n-1}\right)}$$
under the local projection $P\mapsto\left(h_1(P),\cdots,h_{n-1}(P)\right)$ from ${\mathbb C}^n$ to ${\mathbb C}^{n-1}$ and obtain $p\left(h_1,\cdots,h_{n-1}\right)$ from the zero-set of this direct image by taking an effective power.

\medbreak We now form $dh_1\wedge\cdots\wedge dh_n$ and get
$$
\displaylines{dh_1\wedge\cdots\wedge dh_n=dh_1\wedge\cdots\wedge dh_{n-1}\wedge\left(\lambda z_n^{\lambda-1}dz_n+\sum_{j=1}^{\lambda-1}j a_j\left(h_1,\cdots,h_{n-1}\right)z_n^{j-1}dz_n\right)\cr
=\frac{\partial\left(h_1,\cdots,h_{n-1}\right)}{\partial\left(z_1,\cdots,z_{n-1}\right)}\left(\lambda z_n^{\lambda-1}+\sum_{j=1}^{\lambda-1}j a_j\left(h_1,\cdots,h_{n-1}\right)z_n^{j-1}\right)dz_1\wedge\cdots\wedge dz_n.\cr}
$$
Since the coefficient of $dz_1\wedge\cdots\wedge dz_n$ in $dh_1\wedge\cdots\wedge dh_n$ is a multiplier, it follows that
$$\frac{\partial\left(h_1,\cdots,h_{n-1}\right)}{\partial\left(z_1,\cdots,z_{n-1}\right)}\left(\lambda z_n^{\lambda-1}+\sum_{j=1}^{\lambda-1}j a_j\left(h_1,\cdots,h_{n-1}\right)z_n^{j-1}\right)$$
is a multiplier.  Since $\tilde h_n$ is a multiplier, it follows that
$$\tilde h_n\left(\lambda z_n^{\lambda-1}+\sum_{j=1}^{\lambda-1}j a_j\left(h_1,\cdots,h_{n-1}\right)z_n^{j-1}\right)$$
is a multiplier.  From (III.8.1) it follows that
$$\displaylines{\left|p\left(h_1,\cdots,h_{n-1}\right)
\left(\lambda z_n^{\lambda-1}+\sum_{j=1}^{\lambda-1}j a_j\left(h_1,\cdots,h_{n-1}\right)z_n^{j-1}\right)\right|\cr
\stackrel{<}{\sim}\left|
\frac{\partial\left(h_1,\cdots,h_{n-1}\right)}{\partial\left(z_1,\cdots,z_{n-1}\right)}\left(\lambda z_n^{\lambda-1}+\sum_{j=1}^{\lambda-1}j a_j\left(h_1,\cdots,h_{n-1}\right)z_n^{j-1}\right)\right|\cr+\left|
\tilde h_n\left(\lambda z_n^{\lambda-1}+\sum_{j=1}^{\lambda-1}j a_j\left(h_1,\cdots,h_{n-1}\right)z_n^{j-1}\right)\right|.\cr}$$
Hence by the Real Radical Property of Kohn's algorithm in (I.1)(C),
$$p\left(h_1,\cdots,h_{n-1}\right)
\left(\lambda z_n^{\lambda-1}+\sum_{j=1}^{\lambda-1}j a_j\left(h_1,\cdots,h_{n-1}\right)z_n^{j-1}\right)$$
is a multiplier.  Let $h_n^{(0)}=h_n$ and
for $1\leq\nu\leq\lambda$ define
$$
h^{(\nu)}_n=p\left(h_1,\cdots,h_{n-1}\right)^\nu
\left(\frac{\lambda!}{\left(\lambda-\nu\right)!} z_n^{\lambda-\nu}+\sum_{j=\nu}^{\lambda-1}\frac{j!}{\left(j-\nu\right)!} a_j\left(h_1,\cdots,h_{n-1}\right)z_n^{j-\nu}\right).
$$
We are going to verify by induction on $\nu$ that $h^{(\nu)}_n$ is a multiplier.  We know that both $h^{(0)}_n$ and $h^{(1)}_n$ are multipliers.  Assume that we have already verified that
$h^{(0)}_n,\cdots,h^{(\nu-1)}_n$ are multipliers.  Then $dh_1\wedge\cdots\wedge dh_{n-1}\wedge dh^{(\nu-1)}_n$ is equal to
$$
\displaylines{dh_1\wedge\cdots\wedge dh_{n-1}\wedge d\left(
p\left(h_1,\cdots,h_{n-1}\right)^{\nu-1}\left(\frac{\lambda!}{\left(\lambda-\nu+1\right)!} z_n^{\lambda-\nu+1}\right.\right.\cr\hfill\left.\left.+\sum_{j=\nu-1}^{\lambda-1}\frac{j!}{\left(j-\nu+1\right)!} a_j\left(h_1,\cdots,h_{n-1}\right)z_n^{j-\nu+1}\right)\right)\cr
=dh_1\wedge\cdots\wedge dh_{n-1}\wedge \left(
p\left(h_1,\cdots,h_{n-1}\right)^{\nu-1}\left(\frac{\lambda!}{\left(\lambda-\nu\right)!} z_n^{\lambda-\nu}\right.\right.\cr\hfill\left.\left.+\sum_{j=\nu}^{\lambda-1}\frac{j!}{\left(j-\nu\right)!} a_j\left(h_1,\cdots,h_{n-1}\right)z_n^{j-\nu}\right)dz_n\right)\cr
=\frac{\partial\left(h_1,\cdots,h_{n-1}\right)}{\partial\left(z_1,\cdots,z_{n-1}\right)}\left(
p\left(h_1,\cdots,h_{n-1}\right)^{\nu-1}\left(\frac{\lambda!}{\left(\lambda-\nu\right)!} z_n^{\lambda-\nu}\right.\right.\cr\hfill+\left.\left.\sum_{j=\nu}^{\lambda-1}\frac{j!}{\left(j-\nu\right)!} a_j\left(h_1,\cdots,h_{n-1}\right)z_n^{j-\nu}\right)\right)dz_1\wedge\cdots\wedge dz_n.\cr
}
$$
Since the coefficient of $dz_1\wedge\cdots\wedge dz_n$ in $dh_1\wedge\cdots\wedge dh_{n-1}\wedge dh^{(\nu-1)}_n$ is a multiplier, it follows that
$$\displaylines{\frac{\partial\left(h_1,\cdots,h_{n-1}\right)}{\partial\left(z_1,\cdots,z_{n-1}\right)}\left(
p\left(h_1,\cdots,h_{n-1}\right)^{\nu-1}\left(\frac{\lambda!}{\left(\lambda-\nu\right)!} z_n^{\lambda-\nu}\right.\right.\cr\hfill+\left.\left.\sum_{j=\nu}^{\lambda-1}\frac{j!}{\left(j-\nu\right)!} a_j\left(h_1,\cdots,h_{n-1}\right)z_n^{j-\nu}\right)\right)\cr}$$
is a multiplier. Since $\tilde h_n$ is a multiplier, it follows that
$$\tilde h_n\left(
p\left(h_1,\cdots,h_{n-1}\right)^{\nu-1}\left(\frac{\lambda!}{\left(\lambda-\nu\right)!} z_n^{\lambda-\nu}+\sum_{j=\nu}^{\lambda-1}\frac{j!}{\left(j-\nu\right)!} a_j\left(h_1,\cdots,h_{n-1}\right)z_n^{j-\nu}\right)\right)$$
is a multiplier.  From (III.8.1) it follows that
$$\displaylines{\left|p\left(h_1,\cdots,h_{n-1}\right)\left(
p\left(h_1,\cdots,h_{n-1}\right)^{\nu-1}\left(\frac{\lambda!}{\left(\lambda-\nu\right)!} z_n^{\lambda-\nu}\right.\right.\right.\cr\hfill\left.\left.\left.+\sum_{j=\nu}^{\lambda-1}\frac{j!}{\left(j-\nu\right)!} a_j\left(h_1,\cdots,h_{n-1}\right)z_n^{j-\nu}\right)\right)\right|\cr\stackrel{<}{\sim}\left|\frac{\partial\left(h_1,\cdots,h_{n-1}\right)}{\partial\left(z_1,\cdots,z_{n-1}\right)}\left(
p\left(h_1,\cdots,h_{n-1}\right)^{\nu-1}\left(\frac{\lambda!}{\left(\lambda-\nu\right)!} z_n^{\lambda-\nu}\right.\right.\right.\cr\hfill\left.\left.\left.+\sum_{j=\nu}^{\lambda-1}\frac{j!}{\left(j-\nu\right)!} a_j\left(h_1,\cdots,h_{n-1}\right)z_n^{j-\nu}\right)\right)\right|\cr+\left|
\tilde h_n\left(
p\left(h_1,\cdots,h_{n-1}\right)^{\nu-1}\left(\frac{\lambda!}{\left(\lambda-\nu\right)!} z_n^{\lambda-\nu}\right.\right.\right.\cr\hfill\left.\left.\left.+\sum_{j=\nu}^{\lambda-1}\frac{j!}{\left(j-\nu\right)!} a_j\left(h_1,\cdots,h_{n-1}\right)z_n^{j-\nu}\right)\right)\right|.\cr}$$
Hence by the Real Radical Property of Kohn's algorithm in (I.1)(C),
$$\displaylines{p\left(h_1,\cdots,h_{n-1}\right)\left(
p\left(h_1,\cdots,h_{n-1}\right)^{\nu-1}\left(\frac{\lambda!}{\left(\lambda-\nu\right)!} z_n^{\lambda-\nu}\right.\right.\cr\hfill\left.\left.+\sum_{j=\nu}^{\lambda-1}\frac{j!}{\left(j-\nu\right)!} a_j\left(h_1,\cdots,h_{n-1}\right)z_n^{j-\nu}\right)\right)\cr}$$
is a multiplier and $h^{(\nu)}_n$ is a multiplier.  When $\nu=\lambda$, we end up with
$$
h^{(\lambda)}_n=p\left(h_1,\cdots,h_{n-1}\right)^{\lambda}\lambda!
$$
being a multiplier.  Since the multiplicity of the ideal generated by $p\left(h_1,\cdots,h_{n-1}\right)$ and $\tilde h_n$ is effectively bounded at the origin, it follows that the multiplicity of the ideal generated by $p\left(h_1,\cdots,h_{n-1}\right)^{\lambda}$ and $\tilde h_n$ is effectively bounded
at the origin.  We can conclude that $p\left(h_1,\cdots,h_{n-1}\right)$ is a multiplier admitting an order of subellipticity with an effective positive lower bound.

\bigbreak\noindent(III.9) {\it Effective Termination of Kohn's Algorithm in the Higher Dimensional Case.}  Recall that in (III.8) we have the multiplier $h_n$ constructed from $F_1,\cdots,F_N$ and $\tilde h_n$ by choosing $n-1$ good ${\mathbb C}$-linear combinations of $F_1,\cdots,F_N$.  Now we enhance the construction of $h_n$ by choosing $n$ good ${\mathbb C}$-linear combinations of $F_1,\cdots,F_N$ so that any subset of $n-1$ of them are good ${\mathbb C}$-linear combinations for our purpose.  More precisely, as in (III.8) we choose $n$ generic ${\mathbb C}$-linear combinations $H_1,\cdots,H_n$ of $F_1,\cdots,F_N$ such that\begin{itemize}\item[(i)] the map $\pi:{\mathbb C}^n\to{\mathbb C}^n$ defined by $H_1,\cdots,H_n$ is an analytic cover map locally at the origin whose number of sheets is effectively bounded, and\item[(ii)] for any $1\leq j\leq n$ we can use $H_1,\cdots,H_{j-1},H_{j+1},\cdots,H_n$ as $h_1,\cdots,h_{n-1}$ for the argument in (III.8) to produce a polynomial $p_j\left(H_1,\cdots,H_{j-1},H_{j+1},\cdots,H_n\right)$ of $H_1,\cdots,H_{j-1},H_{j+1},\cdots,H_n$ which is a multiplier and whose multiplicity at the origin is effectively bounded and whose assigned order of subellipticity has an effective positive lower bound.
\end{itemize}
The argument in (III.8) shows that each $p_j\left(H_1,\cdots,H_{j-1},H_{j+1},\cdots,H_n\right)$ is an effectively constructed multiplier for $1\leq j\leq n$.

\medbreak We introduce coordinates $z_1,\cdots,z_n$ in the domain space ${\mathbb C}^n$ of the map $\pi:{\mathbb C}^n\to{\mathbb C}^n$.  We use coordinates $w_1,\cdots,w_n$ in the target space ${\mathbb C}^n$ of the map $\pi:{\mathbb C}^n\to{\mathbb C}^n$. Since the polynomial $p_j\left(w_1,\cdots,w_{j-1},w_{j+1},\cdots,w_n\right)$ as a function of $w_1,\cdots,w_{j-1},w_{j+1},\cdots,w_n$ has effectively bounded multiplicity at the origin. it follows that an effectively bounded positive power of the maximum ideal of the target space ${\mathbb C}^n$ of $\pi$ at the origin is contained in the ideal of the target space ${\mathbb C}^n$ of $\pi$ at the origin generated by the $n$ polynomials $p_j\left(w_1,\cdots,w_{j-1},w_{j+1},\cdots,w_n\right)$ for $1\leq j\leq n$.

\medbreak Since
the map $\pi:{\mathbb C}^n\to{\mathbb C}^n$ defined by $H_1,\cdots,H_n$ is an analytic cover map locally at the origin whose number of sheets is effectively bounded, it follows that an
an effectively bounded positive power of the maximum ideal of the domain space ${\mathbb C}^n$ of $\pi$ at the origin is contained in the ideal of the domain space ${\mathbb C}^n$ of $\pi$ at the origin generated by the $n$ holomorphic function germs $p_j\left(H_1,\cdots,H_{j-1},H_{j+1},\cdots,H_n\right)$ for $1\leq j\leq n$.  Since each $p_j\left(H_1,\cdots,H_{j-1},H_{j+1},\cdots,H_n\right)$ is an effectively constructed multiplier for $1\leq j\leq n$, it follows that each of the coordinates $z_1,\cdots,z_n$ of the domain space ${\mathbb C}^n$ of $\pi$ is a multiplier with effective assigned order of subellipticity.  By forming the Jacobian determinant of the multipliers $z_1,\cdots,z_n$, we conclude that the function $F\equiv 1$ is a multiplier whose assigned order of subellipticity has an effective positive lower bound.  This finishes the proof of the effective termination of Kohn's algorithm and concludes the proof of Main Theorem (III.2).

\bigbreak\noindent(III.10) {\it Remark on the Need to Fiberwise Differentiate as Many Times as the Multiplicity of the Given Multiplier.}  An earlier version of this paper puts in the proof only one fiberwise differentiation for the given multiplier instead of the number of fiber differentiations equal to the multiplicity of the multiplier.  This version adds the required number of differentiations.  Let us explain the need to fiberwise differentiate as many times as the multiplicity of the multiplier by considering the following simple situation in complex dimension $2$.

\medbreak Let $f(z,w)$ be a Weierstrass polynomial of degree $q$ in $w$, which is a monic polynomial in $w$ whose coefficients, except the leading one, are holomorphic function germs in $z$ vanishing at the origin.  Denote by $f_w(z,w)$ the derivative of $f(z,w)$ with respect to $w$.   Let $D(z)$ be the discriminant of $f(z,w)$ as a polynomial in $w$.  Then Euclid's algorithm gives $D(z)=a(z,w)f(z,w)+b(z,w)f_w(z,w)$, where $a(z,w)$ and $b(z,w)$ are holomorphic function germs on ${\mathbb C}^2$ at the origin.

\medbreak Note that if $q$ is small, we can only conclude that the multiplicity $q$ of $f(z,w)$ at the origin is small and we cannot conclude that the coefficients of powers of $w$, other than the leading one, have low vanishing order in $z$ at $z=0$.

\medbreak Suppose $f(z,w)$ is a multiplier and $z$ is a pre-multiplier.  When we apply the operator $dz\wedge d\left(\cdot\right)$ to $f(z,w)$ to get $f_w\left(z,w\right)dz\wedge dw$, we conclude that $f_w(z,w)$ is a multiplier.  From $D(z)=a(z,w)f(z,w)+b(z,w)f_w(z,w)$ it follows that the discriminant $D(z)$ is a also multiplier which in general is not effectively constructed.  The vanishing order of $D(z)$ in $z$ at $z=0$ in general does not have anything to do with $q$ and certainly in general cannot be bounded by an effective function of $q$.  Thus the ideal generated by the multipliers $f(z,w)$ and $D(z)$ may have high multiplicity at the origin if $D(z)$ has high vanishing order in $z$ at $z=0$.  This function germ $D(z)$ is obtained by one single fiberwise differentiation of the multiplier $f(z,w)$.

\medbreak The discriminant $D(z)$ is given by $\prod_{i\not=j}\left(w_i(z)-w_j(z)\right)^2$, where $\left\{w_1(z),\cdots,w_q(z)\right\}$ (without any well-defined ordering) is the collection of the $q$ roots of $f(z,w)$ in $w$ with the multiplicities of the roots counted.  If the minimum distance of two points in $\left\{w_1(z),\cdots,w_q(z)\right\}$ as a function of $z$ vanishes to high order in $z$ at $z=0$, we would have high vanishing order for $D(z)$.  The process of getting $D(z)$ by differentiating once does not help in our goal of achieving an effective termination of Kohn's algorithm when two of the roots from the set $\left\{w_1(z),\cdots,w_q(z)\right\}$ are becoming close very fast as $z$ approaches $0$.  Since we have no control over how fast some of the roots $\left\{w_1(z),\cdots,w_q(z)\right\}$ are getting close as $z\to 0$, we need to differentiate $q$ times in order to achieve our goal of an effective termination of Kohn's algorithm.  This explains why we need to fiberwise differentiate as many times as the multiplicity of the multiplier.

\bigbreak\noindent(III.11) {\it Motivation of the Proof of Termination of Kohn's Algorithm from the Fundamental Theorem in Multivariate Calculus for Fubini's Iterated Integration.} We would like to remark that the motivation for the above proof of the termination of Kohn's algorithm for special domains comes from the fundamental theorem in multivariate calculus for the following theorem of Fubini on iterated integration.  The reason for this motivation is that Jacobian determinants occur in the change-of-variables formula for integrals of several variables and that an induction process can be used when we convert an integral of several variables to an iterated integral by Fubini's theorem.

\medbreak\noindent(III.11.1) {\it Fubini's Iterated Integration.}  Let $y_1,\cdots,y_{n-1}$ be functions defining a projection from an $n$-space $G$ with coordinates $x_1,\cdots,x_n$ to an $(n-1)$-space $D$ with coordinates $y=\left(y_1,\cdots,y_{n-1}\right)$ so that $x_n$ can be used to be a local coordinate for the fiber $L_y$ of the projection over the point $y\in D$.  Then for a function $f$ on $G$ the formula
    $$
    \int_G f=\int_{y\in D}\left(\int_{L_y}f\right)
    $$
    holds with the use of appropriate measures.

\medbreak Like the fundamental theorem of calculus of a single real variable, the fundamental theorem in multivariate calculus for the above theorem of Fubini on iterated integration changes integration to differentiation.  If we write the function $f$ in the form $f\left(y_1,\cdots,y_{n-1},x_n\right)$, then
    $$
    dy_1\wedge\cdots\wedge dy_{n-1}\wedge df=
    dy_1\wedge\cdots\wedge dy_{n-1}\wedge\left(\frac{\partial f}{\partial x_n}\right)dx_n
    $$
    so that fiberwise integration over $L_y$ with respect to $x_n$ in (III.11.1) changes over to fiberwise differentiation on $L_y$ with respect to $x_n$.

\medbreak When we use a multiplier as $f$ and pre-multipliers as $y_1,\cdots,y_{n-1}$ to form the Jacobian determinant with respect to $x_1,\cdots,x_n$, we get
       $$\displaylines{
    dy_1\wedge\cdots\wedge dy_{n-1}\wedge df\cr=
    \frac{\partial\left(y_1,\cdots,y_{n-1}\right)}{\partial\left(x_1,\cdots,x_{n-1}\right)}
    \left(\frac{\partial f}{\partial x_n}\right)\left(dx_1\wedge\cdots\wedge dx_n\right).}
    $$
The occurrence of the factor $\frac{\partial f}{\partial x_n}$ enables us to reduce the vanishing order of $f$ by differentiation and the occurrence, as a factor, of the Jacobian determinant $$\frac{\partial\left(y_1,\cdots,y_{n-1}\right)}{\partial\left(x_1,\cdots,x_{n-1}\right)}$$
involving one fewer variable makes it possible to use an induction process.

\eject\noindent{\bf Part IV. Geometric Formulation of Kohn's
Algorithm in Terms of Frobenius Theorem on Integral Submanifolds and
the R\^ole of Real-Analyticity}

\bigbreak Kohn's conjecture for the real-analytic case without
effectiveness was proved by Diederich-Fornaess [DF78].  We are going
to formulate Kohn's algorithm geometrically in terms of the theorem
of Frobenius on integral submanifolds and present a proof of the
real-analytic case of the ineffective termination of Kohn's
algorithm from the geometric viewpoint. This geometric formulation
of Kohn's algorithm in terms of the theorem of Frobenius enables one
to see clearly how the procedures of Kohn's algorithm come about
naturally in the geometric context.  Moreover, the proof of the
real-analytic case of the ineffective termination of Kohn's
algorithm from the geometric viewpoint gives a better understanding
of the r\^ole played by the real-analytic assumption and of the hurdles
standing in the way of generalizing the ineffective
real-analytic case to the ineffective smooth case.

\bigbreak\noindent(IV.1) {\it Usual Theorem of Frobenius on Integral
Submanifolds for ${\mathbb R}^m$.}  The setting of the usual
Frobenius theorem on integral submanifolds of real dimension $k$
starts out with a domain $U$ in ${\mathbb R}^m$ and a distribution
$$x\mapsto W_x\subset T_{{\mathbb R}^m,x}={\mathbb R}^m\ \ {\rm for\ \ }x\in
U$$ which is smooth, where $W_x$ is a $k$-dimensional ${\mathbb
R}$-linear subspace of the tangent space $T_{{\mathbb R}^n,x}$ of
${\mathbb R}^n$ at $x$.

\medbreak The vector-field version of Frobenius's theorem states
that the distribution $x\mapsto W_x$ is locally integrable (in the
sense that locally $U$ is foliated by smooth real submanifolds of
real dimension $k$ whose tangent space at the point $x$ is precisely
$W_x$) if and only if $\left[W_x,W_x\right]\subset W_x$ for all
$x\in U$ (in the sense that for all $x\in U$ the value at $x$ of the
Lie bracket of two local vector fields whose values at $y$ in their
domains of definition are in $W_y$ for each $y$ belongs to $W_x$).

\medbreak The equivalent differential-form version of Frobenius's
theorem states that the distribution $x\mapsto W_x$ is locally
integrable if and only if for any local smooth differential $1$-form
$\omega_1,\cdots,\omega_{m-k}$ whose common kernel is $W_x$ there
exist local smooth differential $1$-forms $\eta_1,\cdots,\eta_{m-k}$
such that $d\omega_j=\sum_{\ell=1}^{m-k}\omega_\ell\wedge\eta_\ell$
for $1\leq j\leq m-k$.

\medbreak The vector-field version of Frobenius's theorem is related
to its differential-form version by Cartan's formula relating Lie
brackets of vector fields and the exterior differentiation of
differential forms (see, for example, [He62, p.21, Formula (9)]).

\bigbreak\noindent(IV.2) {\it Setting of CR Hypersurface for the
Theorem of Frobenius Theorem.}  In the formulation of Kohn's
algorithm in terms of Frobenius's theorem the setting is the
boundary $S$ of a bounded smooth domain $\Omega$ in ${\mathbb C}^n$
and the distribution on $S$ is $P\mapsto T_{S,P}^{\mathbb R}\cap
J\left(T_{S,P}^{\mathbb R}\right)$, where $T_{S,P}^{\mathbb R}$ is
the space of all real tangent vectors in $S$ at $P$ and $J$ is the
almost-complex operator of ${\mathbb C}^n$.

\medbreak In this setting the condition of the theorem of Frobenius
is equivalent to $S$ being Levi-flat, in which case $S$ is locally
foliated by local complex submanifolds of complex dimension $n-1$.

\medbreak The finite type condition of d'Angelo is the opposite of
being Levi-flat.  The finite type condition of d'Angelo can be
interpreted in this context as the impossibility of finding even
Artinian subschemes of arbitrarily high order in $S$ which are
tangential to the distribution $x\mapsto W_x$ of $S$.  The
underlying point set of an Artinian subscheme is just a single
point, but its structure sheaf may be more than the complex number
field ${\mathbb C}$ and can be an Artinian ring ({\it i.e.} a ring
of finite dimension) which is the quotient of the structure sheaf of
$S$.

\medbreak We will not go into the precise definition an Artinian
subscheme here.  Its definition depends on the structure sheaf of
$S$ which in the real-analytic case is the sheaf of germs of all
real-analytic functions and in the smooth case is the sheaf of germs
of all smooth functions.  As an illustration we give here the
following two examples of Artinian subschemes ${\mathcal A}$ of the
ringed space $\left({\mathbb C}^n,{\mathcal O}_{{\mathbb
C}^n}\right)$ supported at the origin of ${\mathbb C}^n$, where
${\mathcal O}_{{\mathbb C}^n}$ is the sheaf of all holomorphic
function germs on ${\mathbb C}^n$.

\medbreak Let ${\mathfrak m}_{{\mathbb C}^n,0}$ be the maximum ideal
at the origin of ${\mathbb C}^n$. Then the ringed space
$\left(\left\{0\right\},\,{\mathcal O}_{{\mathbb
C}^n}\left/\left({\mathfrak m}_{{\mathbb
C}^n,0}\right)^q\right.\right)$ for any positive integer $q$ is an
Artinian subscheme of the ringed space $\left({\mathbb
C}^n,{\mathcal O}_{{\mathbb C}^n}\right)$ supported at the origin of
${\mathbb C}^n$.  For any ideal ${\mathcal I}$ of ${\mathcal
O}_{{\mathbb C}^n,0}$ with $\left({\mathfrak m}_{{\mathbb
C}^n,0}\right)^N\subset{\mathcal I}$ for some positive integer $N$,
the ringed space $\left(\left\{0\right\},\,{\mathcal O}_{{\mathbb
C}^n}\left/{\mathcal I}\right.\right)$ is also an Artinian subscheme
of the ringed space $\left({\mathbb C}^n,{\mathcal O}_{{\mathbb
C}^n}\right)$ supported at the origin of ${\mathbb C}^n$.

\bigbreak\noindent(IV.3) {\it Steps of Kohn's Algorithm from
Constructing Integral Complex Curves.}  We are going to see how the
steps of Kohn's algorithm naturally arise in the study of conditions
necessary for the construction of integral complex curves in the
boundary of a smooth bounded domain.  Again the setting is a weakly
pseudoconvex smooth bounded domain $\Omega$ with boundary $S$ and
again on $S$ we consider the distribution $P\mapsto T_{S,P}^{\mathbb
R}\cap J\left(T_{S,P}^{\mathbb R}\right)$ for $P\in S$.  We will
later specialize to the case where the boundary $S$ of the bounded
domain $\Omega$ is real-analytic and will investigate precisely the
r\^ole played by the assumption of real-analyticity of $S$.  To
anticipate the later specialization into the case of the boundary
$S$ being real-analytic, we would like to explore
conditions which give as a consequence the existence of some local
complex curve in $S$. What we would like to do is to assume that
Kohn's algorithm does not terminate and seek to produce geometrically
a local complex curve in $S$ in the real-analytic case.  For this
purpose, in our discussion, from time to time we will restrict
ourselves to some appropriate open subsets of $S$ in order to
exclude the singularity of real-analytic subsets which arise in our
discussion.

\medbreak Let $N^{(1,0)}_S$ be the set of all $(1,0)$-vectors of $S$
which is in the null space of the Levi form of $S$.  Let ${\mathcal
N}$ be the real part of $N^{(1,0)}_S$ in the sense that at a point
$P$ of $S$ the space ${\mathcal N}$ consists of all ${\rm Re\,}\xi$
with $\xi\in N^{(1,0)}_S$ at the point $P$. Let $T^{\mathbb R}_S$ be
the vector bundle of all real tangent vectors of $S$.  One key
property of ${\mathcal N}$ is the following.

\bigbreak\noindent (IV.3.1) Let $P$ be a point of $S$ and $U$ be
an open neighborhood of $P$ in $S$. Let $\xi$ and $\eta$ be smooth
sections of $T^{\mathbb R}_S$ over $U$. That is, $\xi$ and $\eta$
are real tangent vector fields of $S$ defined on $U$. Assume that
both $\xi$ and $\eta$ belongs to ${\mathcal N}$ at $P$.  Then the
value of the Lie bracket $\left[\xi,\eta\right]$ at $P$ belongs to
${\mathcal N}$.

\bigbreak\noindent Another way to state (IV.3.1) is the following.

\bigbreak\noindent(IV.3.2) The real part of the null space of
$(1,0)$-vectors for the Levi form of a weakly pseudoconvex
boundary is closed under the Lie bracket after extension of the
pointwise vectors to vector fields.

\bigbreak\noindent The statement (IV.3.1) holds mainly because on
$S$ the first derivative of the Levi form for a $(1,0)$-vector
field vanishes at any of its zero points.  The details for its
proof are given in [DF78, Proposition 1].

\medbreak Let $G$ be a nonempty open subset of $S$ where the real
dimension of ${\mathcal N}$ is constant, say $\ell$. For the case
$\ell\geq 1$, it follows from Frobenius theorem and condition
(IV.3.1) that over $G$ we can integrate ${\mathcal N}$ to get local
integral submanifolds $M$ of $G$ of real dimension $\ell$ so that
the tangent space of $M$ at any point $P$ of $M$ is equal to the
real vector space ${\mathcal N}$ at $P$.

\medbreak Since at every point $P$ of $S$ the space $N^{(1,0)}_S$ is
a vector space over the complex number field ${\mathbb C}$, we know
that its real part ${\mathcal N}$ must be invariant under the
almost-complex-structure operator $J$. Thus we can conclude that the
tangent space ${\mathcal N}$ of each local integral submanifold $M$
is invariant under the almost-complex-structure operator $J$.  This
implies that each $M$ is complex-analytic.  As a consequence, one
has the following trivial remark.

\bigbreak\noindent(IV.4) {\it Remark.} Suppose $S$ is a local smooth
weakly pseudoconvex hypersurface in ${\mathbb C}^n$.  If at each
point of $S$ the null space of the Levi form of $S$ is nontrivial,
then some nonempty open subset $G$ of $S$ is foliated by local
complex submanifolds of positive dimension.

\bigbreak\noindent(IV.5) {\it Natural Occurrence of the Steps of
Kohn's Algorithm.} The algorithm of Kohn comes into the picture only
when we do not have a nonempty open subset $G$ of $S$ where the real
dimension of ${\mathcal N}$ is some positive constant $\ell$. We
consider the set $E$ of points of $S$ where the real dimension of
${\mathcal N}$ is some positive constant $\ell$.  The case of
interest is when $E$ does not contain an open subset of $S$.  This
step of introducing $E$ corresponds to introducing the coefficients
of the $(n,n-1)$-form $
\partial r\wedge\left(\partial\bar\partial r\right)^{n-1}
$
as multipliers in Kohn's algorithm. We are going to assume that $E$
is a smooth submanifold of real dimension $m\geq 1$ and that
${\mathcal N}|_E$ is a smooth vector bundle over $E$. In the
real-analytic case because of the stratification of real-analytic
subvarieties we can always get to a real-analytic submanifold and a
real-analytic bundle by replacing the point under consideration by
another point nearby. In the general smooth case there is no such
stratification and the situation becomes complicated and calls for
other techniques than those discussed here.

\medbreak We want to apply Frobenius's theorem to $E$ with the
distribution of vector spaces ${\mathcal N}$ over it.  The trouble
is that the vector space ${\mathcal N}$ at a point $P$ of $E$ may
not be inside the real tangent space $T^{\mathbb
R}_{E,P}=\left(T^{\mathbb R}_E\right)_P$ of $E$ at $P$. To apply
Frobenius's theorem to $E$ we must work with a distribution of
subspaces of the tangent spaces of $E$.  We are forced to replace
${\mathcal N}$ by ${\mathcal N}\cap T^{\mathbb R}_E$ at each point
$P$ of $E$.  We also want to keep the smaller new vector space
${\mathcal N}\cap T^{\mathbb R}_E$ invariant under the
almost-complex-structure operator $J$, because we are interested in
producing local complex curves inside $S$.   We use the even smaller
vector space ${\mathcal N}\cap T^{\mathbb R}_E\cap JT^{\mathbb
R}_E$.  (Note that ${\mathcal N}$ is invariant under $J$.)  Then we
consider the new subset $E_1$ of $E$ where the real dimension of
${\mathcal N}\cap T^{\mathbb R}_E\cap JT^{\mathbb R}_E$ is positive.

\medbreak For the real-analytic case this step corresponds to
introducing real-valued real-analytic function germs $f$ vanishing on $E$ as
multipliers and also $\partial f$ as vector-multipliers.  The reason
is that taking intersection with $T^{\mathbb R}_E$ is the same as
considering the kernel of the differential $df$ of local real-valued
functions $f$ on $S$ which vanish identically in $E$. Taking the
further intersection with $JT^{\mathbb R}_E$ is to consider also the
kernel of the $J$-image $Jdf$ of the differential $df$ of local
real-valued functions $f$ on $S$ which vanish on $E$. Taking both
intersections together is the same as considering the kernel of
$\partial f$ for local real-valued functions $f$ on $S$ which vanish on $E$. Th
use of all local real-analytic function germs vanishing on $E$ tells
us how the step (I.1)(C) of Kohn's algorithm naturally arises from
the geometric viewpoint.  The use of $\partial f$ tells us how the
step (I.1)(B)(i) of Kohn's algorithm naturally arises from the
geometric viewpoint.

\medbreak  As the initial steps of an inductive process we set
$E_0=E$ and ${\mathcal N}^{(0)}={\mathcal N}$ and ${\mathcal
N}^{(1)}={\mathcal N}^{(0)}\cap T^{\mathbb R}_{E_0}\cap JT^{\mathbb
R}_{E_0}$.  Then we inductively define $${\mathcal
N}^{(\nu+1)}={\mathcal N}^{(\nu)}\cap T^{\mathbb R}_{E_\nu}\cap
JT^{\mathbb R}_{E_\nu}={\mathcal N}\cap T^{\mathbb R}_{E_\nu}\cap
JT^{\mathbb R}_{E_\nu}$$ and define $E_{\nu+1}$ to be the subset of
$E$ where the real dimension of ${\mathcal N}^{(\nu)}$ is positive.
We obtain the limiting common intersection $E_\infty$ defined by
$E_\infty=\cap_\nu E_\nu$.  By replacing $E_\infty$ by a nonempty
open subset in the real-analytic case, we can assume that $E_\infty$
is regular and ${\mathcal N}^{(\infty)}:={\mathcal N}\cap T^{\mathbb
R}_{E_\infty}\cap JT^{\mathbb R}_{E_\infty}$ is a real-analytic
vector bundle over $E_\infty$.  Note that, though we go to the
regular part of $E_\nu$ in order to describe more easily the tangent
bundle $T_{E_\nu}$ of $E_\nu$, when we take the intersection
$E_\infty=\cap_\nu E_\nu$ we have to make sure that the intersection
$E_\infty=\cap_\nu E_\nu$ is defined in the real-analytic case as a
real-analytic subvariety, which forces us to consider $T_{E_\nu}$
also at singular points of $E_\nu$ where it is defined as the common
kernel of differentials of all local real-analytic functions
vanishing on $E_\nu$.

\medbreak Note that the definition of $E_{\nu+1}$ as the subset of
$E$ where the real dimension of ${\mathcal N}^{(\nu)}$ is positive
involves the existence of a nontrivial solution in a system of
homogeneous linear equations or equivalently the vanishing of the
determinant of the coefficient matrix or equivalently the vanishing
of the corresponding exterior product of co-vectors.  This tells us
how the step of Kohn's algorithm described in (I.1)(B)(ii) naturally
arises from the geometric viewpoint.

\medbreak Now the distribution of vector spaces ${\mathcal
N}^{(\infty)}={\mathcal N}\cap T^{\mathbb R}_{E_\infty}\cap
JT^{\mathbb R}_{E_\infty}$ is contained in the tangent space of
$E_\infty$ and each ${\mathcal N}^{(\infty)}={\mathcal N}\cap
T^{\mathbb R}_{E_\infty}\cap JT^{\mathbb R}_{E_\infty}$ is
$J$-invariant.  For the purpose of understanding how the procedures
of Kohn's algorithm come about naturally in the geometric context,
we assume that each fiber of ${\mathcal N}^{(\infty)}={\mathcal
N}\cap T^{\mathbb R}_{E_\infty}\cap JT^{\mathbb R}_{E_\infty}$ is of
positive dimension and we also assume that we are in the
real-analytic case so that we have the benefit of stratification.
Under such assumptions and after restriction to a dense open subset
if necessary, $E_\infty$ is a CR manifold and has holomorphic
dimension at least $1$.  However, for ${\mathcal
N}^{(\infty)}={\mathcal N}\cap T^{\mathbb R}_{E_\infty}\cap
JT^{\mathbb R}_{E_\infty}$ we may not have the involutive condition
of the theorem of Frobenius (which means closure under Lie bracket).
In order to apply the theorem of Frobenius we generate a larger
linear subspace of the tangent space of $E_\infty$ by taking
iterated Lie brackets of local sections of the vector bundle
${\mathcal N}^{(\infty)}={\mathcal N}\cap T^{\mathbb
R}_{E_\infty}\cap JT^{\mathbb R}_{E_\infty}$ to generate a new
distribution $\tilde{\mathcal N}$.  This new distribution
$\tilde{\mathcal N}$ now satisfies the following three conditions.
\begin{itemize}\item[(i)] $\tilde{\mathcal N}$ is contained in the
tangent space of $E_\infty$. \item[(ii)] $\tilde{\mathcal N}$ is
involutive in the sense that it is closed under Lie bracket.
\item[(iii)] $\tilde{\mathcal N}$ belongs to the real part of the
null space of the Levi form of $S$.
\end{itemize}
Note that Condition (iii) is a consequence of (IV.3.1).  However, in
general $\tilde{\mathcal N}$ is no longer $J$-invariant.  An
integral submanifold $M$ of $\tilde{\mathcal N}$ has the following
property.  At each point of $M$ the holomorphic dimension of $S$ is
at least $1$.  An open dense subset of $M$ is a CR manifold, but in
general $M$ is not complex-analytic.  We are going to show, with our
present assumption of real-analyticity, that when the Kohn algorithm
does not terminate, we are able to produce some local complex curve
inside $M$. One key point here is that the tangent space of $M$ is
contained in the null space ${\mathcal N}$ of the Levi-form which is
$J$-invariant.

\medbreak Since we have assumed that we are in the real-analytic
case, at a generic point of $M$ we can consider the smallest complex
submanifold germ $V$ in ${\mathbb C}^n$ which contains the germ of
$M$ at that point.  We then have the following situation. At a
generic point $P_0$ of $M$ there exist
\begin{itemize}\item[(i)] an open neighborhood $U_0$ of $P_0$ in ${\mathbb C}^n$,
\item[(ii)]
a complex submanifold $V$ in $U_0$, and \item[(iii)] real-valued
real-analytic functions $\rho_1,\cdots,\rho_\ell$ on $V$
\end{itemize}
such that
\begin{itemize}\item[(a)]
$M\cap U_0$ is the common zero-set of $\rho_1,\cdots,\rho_\ell$,
\item[(b)]
$\partial\rho_1,\cdots,\partial\rho_\ell$ are ${\mathbb C}$-linearly
independent at points of $M\cap U_0$, and
\item[(c)] at any point of $M\cap U_0$ the tangent space of $V$ is contained
in ${\mathcal N}$.\end{itemize} Condition (b) means that, besides
the ${\mathbb R}$-linear independence of $d\rho_1,\cdots,d\rho_\ell$
at points of $M\cap U_0$, we also have the ${\mathbb R}$-linear
independence of
$$\left(Jd\rho_1\right)|_{T^{\mathbb
R}_M},\cdots,\left(Jd\rho_\ell\right)|_{T^{\mathbb R}_M}$$ at points
of $M\cap U_0$.  The complex dimension of $T^{\mathbb R}_M\cap
JT^{\mathbb R}_M$ is equal to $\dim_{\mathbb C}V-\ell$, which is
$\geq 1$.  The reason why the smallest complex submanifold germ $V$
of ${\mathbb C}^n$ at $P_0$ containing the germ of $M$ at $P_0$
satisfies condition (c) is that $T^{\mathbb R}_M$ is contained in
the $J$-invariant vector space ${\mathcal N}$ at any point of $M$
and we can determine $V$ as the zero-set of holomorphic function
germs on ${\mathbb C}^n$ at $P_0$ obtained by extending CR
real-analytic functions on $M$ by using the condition of their
annihilation by $\bar\partial$ to define the infinite jets of their
extensions.

\medbreak Condition (c) means that $V$ is tangential to $S$ at
points of $M$. There are two possibilities. One is that $V$ is
contained in $S$, in which case $S$ contains a local complex curve
and we are done.  The other possibility is that $V$ is not contained
in $S$. We are going to assume the second possibility and derive a
contradiction for the real-analytic case so that we can conclude in
the real-analytic case that $S$ must contain a local complex curve.

\medbreak For clarity in the later discussion we digress at this
point to say something about the well-known alternative description
of the Levi form and also about the process of polarization.

\bigbreak\noindent(IV.6) {\it Alternative Description of Levi Form.}
Recall the following formula of Cartan for exterior differentiation
of differential forms
$$\displaylines{(p+1)\left(d\omega\right)\left(X_1,\cdots,X_{p+1}\right)=\sum_{i=1}^{p+1}(-1)^{i+1}X_i\left(\omega\left(X_1,\cdots,
\hat{X_i},\cdots,X_{p+1}\right)\right)\cr+\sum_{i<j}(-1)^{i+j}\omega
\left(\left[X_i,X_j\right],X_1,\cdots,\hat{X_i},\cdots,\hat{X_j},\cdots,X_{p+1}\right),\cr}
$$
where $\omega$ is a differential $p$-form and $X_1,\cdots,X_{p+1}$
are vector fields and $\hat X_i$ and $\hat X_j$ indicate that $X_i$
and $X_j$ have been removed (see, for example, [He62, p.21, Formula
(9)]). For the special case of $p=1$ we have
$$2d\omega\left(X_1,X_2\right)=X_1\left(\omega\left(X_2\right)\right)
-X_2\left(\omega\left(X_1\right)\right)-\omega
\left(\left[X_1,X_2\right]\right).
$$
For a function $\rho$ and vector fields $\xi$ and $\eta$, by
applying this to the case of $\omega=Jd\rho$ (where $J$ is the
almost-complex-structure operator), we get
$$2dJd\rho\left(\xi,\eta\right)=\xi\left(Jd\rho\left(\eta\right)\right)
-\eta\left(Jd\rho\left(\xi\right)\right)-Jd\rho
\left(\left[\xi,\eta\right]\right).
$$
Note that $\partial\rho=\frac{1}{2}\left(1-\sqrt{-1}\,J\right)d\rho$
and $\bar\partial\rho=\frac{1}{2}\left(1+\sqrt{-1}\,J\right)d\rho$
so that
$$\partial\bar\partial\rho=d\bar\partial\rho=\frac{1}{2}\,dd\rho
+\frac{\sqrt{-1}}{2}\,dJd\rho=\frac{\sqrt{-1}}{2}\,dJd\rho$$
and $dJd\rho=-2\sqrt{-1}\partial\bar\partial\rho$.  When
$\left(d\rho\right)\left(\xi\right)=0$, we have
$$\left(\partial\rho\right)\left(\xi\right)=
\left(\frac{1}{2}\left(1-\sqrt{-1}\,J\right)d\rho\right)\left(\xi\right)
=-\frac{\sqrt{-1}}{2}\left(Jd\rho\right)\left(\xi\right).$$ When both
$\left(d\rho\right)\left(\xi\right)\equiv 0$ and
$\left(Jd\rho\right)\left(\xi\right)\equiv 0$, we have $$Jd\rho
\left(\left[\xi,\eta\right]\right)=4\sqrt{-1}\left(\partial\bar\partial\rho\right)\left(\xi,\eta\right).$$
When we compute the Levi form of $\rho$ we limit ourselves to
vectors of type $(1,0)$ which are tangential to $\rho=0$.  A vector
$\xi$ of type $(1,0)$ means that $J(\xi)=\sqrt{-1}\,\xi$.  Tangency of $\xi$
to $\rho=0$ means that $\left(d\rho\right)(\xi)=0$, which implies
automatically
$\left(Jd\rho\right)(\xi)=\left(d\rho\right)(J\xi)=\sqrt{-1}\left(d\rho\right)(\xi)=0$,
because by definition the operator $J$ acting on $1$-forms is the
adjoint of the operator $J$ acting on tangent vectors.  Likewise,
for a vector $\bar\xi$ of type $(0,1)$ tangential to $\rho=0$ we
have $\left(d\rho\right)(\bar\xi)=0$ and
$\left(Jd\rho\right)(\bar\xi)=\left(d\rho\right)(J\bar\xi)=-\sqrt{-1}\left(d\rho\right)(\bar\xi)=0$.
Thus for vector fields $\xi$ and $\eta$ of type $(1,0)$ or $(1,0)$
tangential to $\rho=0$ we have
$$Jd\rho
\left(\left[\xi,\eta\right]\right)=4\sqrt{-1}\left(\partial\bar\partial\rho\right)\left(\xi,\eta\right).$$

\bigbreak\noindent(IV.7) {\it Polarization.} Let $Y$ be a CR
submanifold of some open subset of ${\mathbb C}^n$.  Let
$\xi_1,\xi_2$ be real-valued vector fields in $T^{\mathbb R}_Y\cap
JT^{\mathbb R}_Y$. The condition that $\xi_j$ is in $T^{\mathbb
R}_Y\cap JT^{\mathbb R}_Y$ is equivalent to the condition that we
can write $\xi_j=\tau_j+\overline{\tau_j}$ for some complex-valued
vector fields $\tau_j$ in $T^{(1,0)}_Y$ for $j=1,2$. We have
$$\left[\xi_1,\xi_2\right]=\left[\tau_1+\overline{\tau_1},\tau_2+\overline{\tau_2}\right]
=\left[\tau_1,\tau_2\right]+\left[\tau_1,\overline{\tau_2}\right]-
\left[\tau_2,\overline{\tau_1}\right]+\left[\overline{\tau_1},\overline{\tau_2}\right].
$$
For the purpose of later computation of the Levi form, we now
introduce the standard linear polarization process to express
$\left[\tau_2,\overline{\tau_1}\right]$ in terms of
$\left[\tau,\overline{\tau}\right]$ for some vector field $\tau$
of type $(1,0)$ tangential to $Y$ so that $\tau$ is expressed
linearly and explicitly in terms of $\xi_1, \xi_2, J\xi_1, J\xi_2$
modulo ${\mathbb C}\otimes_{\mathbb R}\left(T^{\mathbb R}_Y\cap
JT^{\mathbb R}_Y\right)=T^{(1,0)}_Y\oplus T^{(0,1)}_Y$. From
$$
\left[\tau_1+\tau_2,\overline{\tau_1+\tau_2}\right]=
\left[\tau_1,\overline{\tau_1}\right]+
\left[\tau_1,\overline{\tau_2}\right]+
\left[\tau_2,\overline{\tau_1}\right]+
\left[\tau_2,\overline{\tau_2}\right]
$$
we subtract the expression with $\tau_2$ changed to $-\tau_2$ to get
$$\left[\tau_1+\tau_2,\overline{\tau_1+\tau_2}\right]-
\left[\tau_1-\tau_2,\overline{\tau_1-\tau_2}\right]
=2\left[\tau_1,\overline{\tau_2}\right]+
2\left[\tau_2,\overline{\tau_1}\right].
$$
Then we add to it $\sqrt{-1}$ times the expression which is obtained by
replacing $\tau_2$ by $\sqrt{-1}\tau_2$ and we get
$4\left[\tau_1,\overline{\tau_2}\right]$ equal to
$$\displaylines{\qquad\qquad\qquad\left[\tau_1+\tau_2,\overline{\tau_1+\tau_2}\right]-
\left[\tau_1-\tau_2,\overline{\tau_1-\tau_2}\right]\hfill\cr\hfill
+\sqrt{-1}\left(\left[\tau_1+\sqrt{-1}\tau_2,\overline{\tau_1+\sqrt{-1}\tau_2}\right]-
\left[\tau_1-\sqrt{-1}\tau_2,\overline{\tau_1-\sqrt{-1}\tau_2}\right]\right).\cr}
$$
Since $\left[\tau_1,\tau_2\right]$ is in $T^{(1,0)}_Y$ and
$\left[\overline{\tau_1},\overline{\tau_2}\right]$ is in
$T^{(0,1)}_Y$, we conclude that modulo ${\mathbb C}\otimes_{\mathbb
R}\left(T^{\mathbb R}_Y\cap JT^{\mathbb
R}_Y\right)=T^{(1,0)}_Y\oplus T^{(0,1)}_Y$ the Lie bracket
$\left[\xi_1,\xi_2\right]$ is equal to
$\left[\tau_1,\overline{\tau_2}\right]-\left[\tau_2,\overline{\tau_1}\right]
$ which is in turn equal to $\frac{1}{4}$ times
$$\displaylines{\left[\tau_1+\tau_2,\overline{\tau_1+\tau_2}\right]-
\left[\tau_1-\tau_2,\overline{\tau_1-\tau_2}\right]\cr
+\sqrt{-1}\left(\left[\tau_1+\sqrt{-1}\tau_2,\overline{\tau_1+\sqrt{-1}\tau_2}\right]-
\left[\tau_1-\sqrt{-1}\tau_2,\overline{\tau_1-\sqrt{-1}\tau_2}\right]\right)\cr
-\bigg\{\left[\tau_1+\tau_2,\overline{\tau_1+\tau_2}\right]-
\left[\tau_2-\tau_1,\overline{\tau_2-\tau_1}\right]\cr
+\sqrt{-1}\left(\left[\tau_2+\sqrt{-1}\tau_1,\overline{\tau_2+\sqrt{-1}\tau_1}\right]-
\left[\tau_2-\sqrt{-1}\tau_1,\overline{\tau_2-\sqrt{-1}\tau_1}\right]\right)\bigg\}\cr
=2\sqrt{-1}\left(\left[\tau_1+\sqrt{-1}\tau_2,\overline{\tau_1+\sqrt{-1}\tau_2}\right]-
\left[\tau_1-\sqrt{-1}\tau_2,\overline{\tau_1-\sqrt{-1}\tau_2}\right]\right).\cr }$$
Thus modulo ${\mathbb C}\otimes_{\mathbb R}\left(T^{\mathbb R}_Y\cap
JT^{\mathbb R}_Y\right)=T^{(1,0)}_Y\oplus T^{(0,1)}_Y$ the Lie
bracket $\left[\xi_1,\xi_2\right]$ is equal to
$$\frac{\sqrt{-1}}{2}\left(\left[\tau_1+\sqrt{-1}\tau_2,\overline{\tau_1+\sqrt{-1}\tau_2}\right]-
\left[\tau_1-\sqrt{-1}\tau_2,\overline{\tau_1-\sqrt{-1}\tau_2}\right]\right),$$
where
$$
\displaylines{\tau_1+\sqrt{-1}\tau_2=\frac{1}{2}\left(\xi_1+J\xi_2\right)+\frac{\sqrt{-1}}{2}\left(\xi_2-J\xi_1\right),\cr
\tau_1-\sqrt{-1}\tau_2=\frac{1}{2}\left(\xi_1-J\xi_2\right)-\frac{\sqrt{-1}}{2}\left(\xi_2+J\xi_1\right),\cr}
$$
because $\tau_j=\frac{1}{2}\left(\xi_j-\sqrt{-1}\,J\xi_j\right)$ for
$j=1,2$.  Suppose $\rho$ is a real-valued function in some
neighborhood of $Y$.  Then by (IV.6) we have
$$\displaylines{Jd\rho
\left(\left[\tau_1,\tau_2\right]\right)=4\sqrt{-1}\left(\partial\bar\partial\rho\right)\left(\tau_1,\tau_2\right)=0,\cr
Jd\rho \left(\left[\overline{\tau_1},\overline{\tau_2}\right]\right)
=4\sqrt{-1}\left(\partial\bar\partial\rho\right)\left(\overline{\tau_1},\overline{\tau_2}\right)=0\cr}$$
and as a consequence
$$\displaylines{\left(Jd\rho\right)\left(\left[\xi_1,\xi_2\right]\right)=
\frac{\sqrt{-1}}{2}\left(Jd\rho\right)\left(\left[\tau_1+\sqrt{-1}\tau_2,\overline{\tau_1+\sqrt{-1}\tau_2}\right]\right)\cr\hfill-
\frac{\sqrt{-1}}{2}\left(Jd\rho\right)\left(\left[\tau_1-\sqrt{-1}\tau_2,\overline{\tau_1-\sqrt{-1}\tau_2}\right]\right).\cr}$$
When
$\left|\left(Jd\rho\right)\left(\left[\xi_1,\xi_2\right]\right)\right|=C$
for some $C>0$, we have
$\left|\left(Jd\rho\right)\left(\left[\tau,\overline{\tau}\right]\right)\right|\geq
C$ for one of the following two values of $\tau$.
$$
\displaylines{\tau_1+\sqrt{-1}\tau_2=\frac{1}{2}\left(\xi_1+J\xi_2\right)+\frac{\sqrt{-1}}{2}\left(\xi_2-J\xi_1\right),\cr
\tau_1-\sqrt{-1}\tau_2=\frac{1}{2}\left(\xi_1-J\xi_2\right)-\frac{\sqrt{-1}}{2}\left(\xi_2+J\xi_1\right),\cr}
$$

\bigbreak\noindent(IV.8) {\it Locating Holomorphic Direction at
Which Precisely One Levi-Form Is Nonzero.} After the above
digression on the alternative description of the Levi form and the
process of polarization, we now go back to the situation of the CR
submanifold $M$ at the end of (IV.5).  According to the construction
of $M$ as an integral submanifold of $\tilde{\mathcal N}$ the
tangent bundle $T_M^{\mathbb R}$ of $M$ is generated by iterated Lie
brackets of vector fields of ${\mathcal N}\cap T_{E_\infty}\cap
JT_{E_\infty}$ defined on $M$. Moreover, we have
$$
{\mathcal N}\cap T_{E_\infty}\cap JT_{E_\infty}\subset T^{\mathbb
R}_M\cap JT^{\mathbb R}_M=\bigcap_{j=1}^\ell{\rm
Ker}\left(\left(Jd\rho_j\right)|_{T^{\mathbb R}_M}\right)\subset
T^{\mathbb R}_M.
$$
When we take vector fields in ${\mathcal N}\cap T_{E_\infty}\cap
JT_{E_\infty}$ defined on $M$ and form their iterated Lie brackets
in order to generate $T^{\mathbb R}_M$, there is a first time the
vector field fails to be inside $T^{\mathbb R}_M\cap JT^{\mathbb
R}_M=\bigcap_{j=1}^\ell{\rm
Ker}\left(\left(Jd\rho_j\right)|_{T^{\mathbb R}_M}\right)$.  Thus
we can find real-valued vector fields $\xi_1,\xi_2$ in $T^{\mathbb
R}_M\cap JT^{\mathbb R}_M=\bigcap_{j=1}^\ell{\rm
Ker}\left(\left(Jd\rho_j\right)|_{T^{\mathbb R}_M}\right)$ defined
on $M$ such that their Lie bracket $\left[\xi_1,\xi_2\right]$ is
not in $T^{\mathbb R}_M\cap JT^{\mathbb
R}_M=\bigcap_{j=1}^\ell{\rm
Ker}\left(\left(Jd\rho_j\right)|_{T^{\mathbb R}_M}\right)$.  There
exists $1\leq j\leq\ell$ such that
$\left(Jd\rho_j\right)\left(\left[\xi_1,\xi_2\right]\right)$ is
nonzero.  Without loss of generality we assume that $j=1$ so that
$\left(Jd\rho_1\right)\left(\left[\xi_1,\xi_2\right]\right)$ is
nonzero.  Since $\xi_1,\xi_2$ are both in $T^{\mathbb R}_M\cap
JT^{\mathbb R}_M$, we can write $\xi_j=\tau_j+\overline{\tau_j}$
for some complex-valued vector fields $\tau_j$ in $T^{(1,0)}_M$
for $j=1,2$.  As explained above in (IV.7), the polarization
process gives us
$$\displaylines{\left(Jd\rho_1\right)\left(\left[\xi_1,\xi_2\right]\right)=
\frac{\sqrt{-1}}{2}\left(Jd\rho_1\right)\left(\left[\tau_1+\sqrt{-1}\tau_2,\overline{\tau_1+\sqrt{-1}\tau_2}\right]\right)\cr\hfill-
\frac{\sqrt{-1}}{2}\left(Jd\rho_1\right)\left(\left[\tau_1-\sqrt{-1}\tau_2,\overline{\tau_1-\sqrt{-1}\tau_2}\right]\right),\cr}$$
where
$$
\displaylines{\tau_1+\sqrt{-1}\tau_2=\frac{1}{2}\left(\xi_1+J\xi_2\right)+\frac{\sqrt{-1}}{2}\left(\xi_2-J\xi_1\right),\cr
\tau_1-\sqrt{-1}\tau_2=\frac{1}{2}\left(\xi_1-J\xi_2\right)-\frac{\sqrt{-1}}{2}\left(\xi_2+J\xi_1\right),\cr}
$$
One of
$\left(Jd\rho_1\right)\left(\left[\tau_1+\sqrt{-1}\tau_2,\overline{\tau_1+\sqrt{-1}\tau_2}\right]\right)$
and
$\left(Jd\rho_1\right)\left(\left[\tau_1-\sqrt{-1}\tau_2,\overline{\tau_1-\sqrt{-1}\tau_2}\right]\right)$
must be nonzero at $P_0$.  We can choose $\tau$ to be either
$\tau_1+\sqrt{-1}\tau_2$ or $\tau_1-\sqrt{-1}\tau_2$ so that
$\left(Jd\rho_1\right)\left(\left[\tau,\overline{\tau}\right]\right)$
is nonzero at $P_0$.  Since $\xi_1,\xi_2$ belong to $T^{\mathbb
R}_M\cap JT^{\mathbb R}_M=\bigcap_{j=1}^\ell{\rm
Ker}\left(\left(Jd\rho_j\right)|_{T^{\mathbb R}_M}\right)$, it
follows that $\left(\partial\rho_j\right)\left(\tau\right)=0$ at
$P_0$ for $1\leq j\leq\ell$.

\medbreak Now for $2\leq j\leq\ell$ we replace $\rho_j$ by
$$
\rho_j-\frac{\partial\rho_j\left(\left[\tau,\overline{\tau}\right]\right)}
{\partial\rho_1\left(\left[\tau,\overline{\tau}\right]\right)}\,\rho_1
$$
so that we can assume without loss of generality that
$$\displaylines{0\equiv\left(\partial\rho_j\right)\left(\left[\tau,\overline{\tau}\right]\right)
=\left(\frac{1}{2}\left(1-\sqrt{-1}\,J\right)d\rho_j\right)\left(\left[\tau,\overline{\tau}\right]\right)\cr
=-\,\frac{\sqrt{-1}}{2}\left(Jd\rho_j\right)\left(\left[\tau,\overline{\tau}\right]\right)
\quad{\rm for\ \ }2\leq j\leq\ell.\cr}
$$
We can write
$$
r|_V=\sum_{\nu_1+\cdots+\nu_\ell=k}\sigma_{\nu_1,\cdots,\nu_\ell}\left(\rho_1\right)^{\nu_1}\cdots\left(\rho_\ell\right)^{\nu_\ell}
+O\left(\left(\sum_{j=1}^\ell\left(\rho_j\right)^2\right)^{\frac{k+1}{2}}\right)
$$
for some integer $k\geq 2$, where $\sigma_{\nu_1,\cdots,\nu_\ell}$
is a real-analytic function on $U_0$ (after shrinking $U_0$ as an
open neighborhood of $P_0$ in ${\mathbb C}^n$ if necessary) and
$\sigma_{\nu_1^*,\cdots,\nu_\ell^*}$ is nonzero at $P_0$ for some
$\nu_1^*+\cdots+\nu_\ell^*=k$.

\bigbreak\noindent(IV.9) {\it Argument of Different Vanishing Orders
for Complex Hessian on the Complex Tangent Space Along Vector Fields
Tangential or Normal to the Intersection with the Weakly
Pseudoconvex Boundary.} To make the argument more transparent and
more understandable, we will first consider the special case
$\ell=1$ so that $M=V\cap\left\{\rho_1=0\right\}$ and $V$ is a
complex submanifold in some open neighborhood of some point $P_0$ of
$M$. For this special case, for notational simplicity we drop the
subscript $1$ from $\rho_1$ and simply denote $\rho_1$ by $\rho$. By
replacing $\rho$ by its product with a local nowhere zero
real-analytic function we can assume without loss of generality that
$r=\rho^k$ on $V$.

\medbreak Let $m$ be the complex dimension of $V$.  We choose a
local holomorphic coordinate system $\left(z_1,\cdots,z_n\right)$ on
the open neighborhood $U_0$ of $P_0$ in ${\mathbb C}^n$ centered at
$P_0$ (after shrinking $U_0$ if necessary) such that $S\cap
U_0\cap\left\{z_{m+1}=\cdots=z_n=0\right\}$ is regular and
$V=\left\{z_{m+1}=\cdots=z_n=0\right\}\cap U_0$.  Since our argument
will be confined to an open neighborhood of $P_0$ in ${\mathbb
C}^n$, for notational simplicity, by replacing ${\mathbb C}^n$ by
${\mathbb C}^{m+1}$ and $S$ by
$S\cap\left\{z_{m+1}=\cdots=z_n=0\right\}$ we can assume without
loss of generality that $n=m+1$ and we have the following setup.

\begin{itemize}\item[(i)]
$dr=\left(0,0,\cdots,0,1\right)$ at the origin so that the complex
submanifold $V$ of the neighborhood $U_0$ of $P_0$ in ${\mathbb
C}^n$ is an open subset of the complex tangent space of $S$ at the
origin which is defined by $z_n=0$.
\item[(ii)] $r=x_n+O\left(|z|^2\right)$ near the origin, where $x_n$
is the real part of the coordinate $z_n$.
\item[(iii)] The intersection $M=V\cap S$ of $V$ and $S$ is a CR
manifold whose complex tangent space $T^{\mathbb R}_M\cap
J\left(T^{\mathbb R}_M\right)$ has positive complex dimension at
every point of $M$.
\item[(iv)] $M$ is defined by $\rho=0$ in $V$ for some real-valued real-analytic
function $\rho$ on $V$ such that $r|_V=\rho^k$ for some positive
integer $k$ and $d\rho$ is nowhere zero on $M$.
\item[(v)] For some nonzero tangent vector $\tau$ of type $(1,0)$ tangential to $M$ at the origin
the value of the Levi form of $\rho$ at $\tau$ is nonzero.
\end{itemize}
We are going to derive a contradiction.  First we sketch the main
idea of the argument.  On $V$ we will introduce two vector fields of
type $(1,0)$. One is tangential to $M$ at points of $M$ and the
other is normal to $M$ at points of $M$.  when we compute the
complex Hessian of $\rho^k$ at these two vector fields of type
$(1,0)$ on $V$, we get two different orders of vanishing as we
approach $M$ from $V-M$, one of order $k-1$ and the other of order
$k-2$.  Because the touching order between $V$ and $S$ is $k$ along
$M$, when we extend these two vector fields of type $(1,0)$ on $M$
to an open neighborhood of $P_0$ in ${\mathbb C}^n$ so that the two
extensions are tangential to $S$ at points of $S$, the Levi forms of
$r$ with respect to the two extensions give again the two different orders of
vanishing as we approach $M$ from $S-M$.  Since one of the two
orders is odd, the weak pseudoconvexity of $S$ is violated, yielding
a contradiction. Now we give below the details of this argument of
different vanishing orders for the complex Hessian of $r$ on $V$
along vector fields tangential or normal to its intersection $M$
with the weakly pseudoconvex boundary $S$.

\medbreak There is some open neighborhood $U_1$ of the origin $0$ in
$U_0$ on which
$$
r\left(z_1,\cdots,z_n\right)=\phi\left(z_1,\cdots,z_n\right)z_n+\overline{\phi\left(z_1,\cdots,z_n\right)}\,\overline{z_n}+\rho\left(z_1,\cdots,z_{n-1}\right)^k
$$
for some smooth complex-valued function
$\phi\left(z_1,\cdots,z_n\right)$ on $U_1$, because on
$V=\left\{z_n=0\right\}$ the function $r$ is of the form $\rho^k$.
Let $\xi$ be any smooth vector field of type $(1,0)$ on $U_1$ whose
$n$-th component is $\xi_n$. Then
$$
\partial r=\partial\phi\,z_n+\phi
dz_n+\left(\partial\bar\phi\right)\overline{z_n}+k\rho^{k-1}\partial\rho.
$$
$$ \left<\partial
r,\xi\right>=\left<\partial\phi\,,\xi\right>z_n+\phi
\xi_n+\left<\partial\bar\phi,\xi\right>\overline{z_n}+k\rho^{k-1}\left<\partial\rho,\xi
\right>.\leqno{({\rm IV}.9.1)}
$$
$$\displaylines{\bar\partial
r=\bar\partial\phi\,z_n+\left(\bar\partial\bar\phi\right)\overline{z_n}+\bar\phi\,d\overline{z_n}
+k\rho^{k-1}\bar\partial\rho.\cr\partial\bar\partial
r=\partial\bar\partial\phi\,z_n-\bar\partial\phi
dz_n+\partial\bar\partial\bar\phi\,\overline{z_n}+\partial\bar\phi\,d\overline{z_n}+
k(k-1)\rho^{k-2}\partial\rho\bar\partial\rho+k\rho^{k-1}\partial\bar\partial\rho.\cr}
$$
$$\displaylines{
({\rm IV}.9.2)\qquad\left<\partial\bar\partial
r,\xi\wedge\bar\xi\right>=\left<\partial\bar\partial\phi\,,\xi\wedge\bar\xi\right>z_n-\left<\bar\partial\phi,
\bar\xi\right>
\xi_n+\left<\partial\bar\partial\bar\phi\,,\xi\wedge\bar\xi\right>\overline{z_n}+\left<\partial\bar\phi,\,\xi\right>\bar\xi_n\hfill\cr\hfill
+k(k-1)\rho^{k-2}\left<\partial\rho,\xi\right>\,\left<\bar\partial\rho,\bar\xi\right>+k\rho^{k-1}\left<\partial\bar\partial\rho\,,\xi\wedge\bar\xi\right>.\cr}
$$
At a point of $r=0$ in $U_1$ we have
$$
\phi z_n+\bar\phi\overline{z_n}+\rho^k=0.
$$
Let $A$ and $B$ be respectively the real and imaginary parts of
$2\bar\phi$. Then $\phi=\frac{A-Bi}{2}$ and
$\bar\phi=\frac{A+Bi}{2}$ so that
$$
\phi z_n+\bar\phi\overline{z_n}=Ax_n+By_n
$$
(where $y_n$ is the imaginary part of $z_n$) and at a point in $U_1$
we have
$$
Ax_n+By_n+\rho^k=0.
$$
Since $dr=\left(0,0,\cdots,0,1\right)$ at the origin, it follows
that $\phi=\frac{1}{2}$ at the origin and $A=1$ and $B=0$ at the
origin. Let $Y$ be the set defined by $y_n=0$.  At any point of
$S\cap Y\cap U$ where $A$ is nonzero, we have
$$x_n=-\frac{\rho^k}{A},\qquad
z_n=-\frac{\rho^k}{A},\qquad\overline{z_n}=-\frac{\rho^k}{A}.
$$
We can choose an open neighborhood $U$ of the origin $0$ in $U_1$ of
the form $U=W\times G$ with $W\subset{\mathbb C}^{n-1}$ and
$G\subset{\mathbb C}$ such that \begin{itemize}\item[(i)] $A$ is
nowhere zero on $U$ and for $Q\in W$ the set $G$ contains the point
$z_n=-\frac{\rho(Q)^k}{A}$ and \item[(ii)]
$\phi+\left(\partial_n\phi\right)z_n+\left(\partial_n\bar\phi\right)\overline{z_n}$ is
nowhere zero on $U_1$.\end{itemize} On $S\cap Y\cap U$ the two
functions $z_n$ and $\overline{z_n}$ are of the order
$O\left(\rho^k\right)$.

\bigbreak We now derive our contradiction by choosing $\xi$ in two
different ways.  The first way is to choose $\xi$ equal to $\tau$ at
the origin.  Since $\tau$ (from (IV.8)) is a vector of ${\mathbb
C}^n$ of type $(1,0)$ at the origin which is tangential to $E=V\cap
S$ and since $V=\left\{z_n=0\right\}$, it follows that the $n$-th
component of the $n$-vector $\tau$ is zero. Since the differential
$d\rho$ of the real-valued function $\rho$ on $V\cap U$ is nowhere
zero at every point of $E=V\cap S$, we can extend $\tau$ to some
smooth $(1,0)$-vector field
$\xi=\left(\xi_1,\cdots,\xi_{n-1}\right)$ of $W_1$ for some open
neighborhood $W_1$ of $0$ in $W$ such that
$\left<d\rho,\,\xi\right>\equiv 0$ on $W_1$.

\medbreak We regard $\xi_j=\xi_j\left(z_1,\cdots,z_{n-1}\right)$ as
functions of $\left(z_1,\cdots,z_{n-1},z_n\right)\in W_1\times G$
for $1\leq j\leq n-1$ (which means denoting also by $\xi_j$ the
composite of $\xi_j$ and the natural projection $W_1\times G\to W_1$
for $1\leq j\leq n-1$).  Since
$\phi+\left(\partial_n\phi\right)z_n+\left(\partial_n\bar\phi\right)\overline{z_n}$ is
nowhere zero on $U_1$, we can define $\xi_n$ on $W_1\times G$ by
$$
\xi_n=\frac{-1}{\phi+\left(\partial_n\phi\right)z_n+\left(\partial_n\bar\phi\right)\overline{z_n}}
\left(\sum_{j=1}^{n-1}\left(\partial_j\phi\right)
\xi_j\,z_n+\sum_{j=1}^{n-1}\left(\partial_j\bar\phi\right)\xi_j\,\overline{z_n}\right)\leqno{({\rm
IV}.9.3)}
$$
so that the vector field
$\xi=\left(\xi_1,\cdots,\xi_{n-1},\xi_n\right)$ on $W_1\times G$
satisfies $\left<\partial r,\,\xi\right>\equiv 0$ because of
(IV.9.1). Since on $S\cap Y\cap U$ the two functions $z_n$ and
$\overline{z_n}$ are of the order $O\left(\rho^k\right)$, it follows
from (IV.9.3) that $\xi_n$ is of the order $O\left(\rho^k\right)$ on
$S\cap Y\cap\left(W_1\times G\right)$. By (IV.9.2)
$$\left<\partial\bar\partial
r,\xi\wedge\bar\xi\right>=k\rho^{k-1}\left<\partial\bar\partial\rho\,,\xi\wedge\bar\xi\right>+O\left(\rho^k\right)$$
on $S\cap Y\cap\left(W_1\times G\right)$. Since at the origin
$\left<\partial\bar\partial\rho,\,\xi\wedge\bar\xi\right>=\left<\partial\bar\partial\rho,\,\tau\wedge\bar\tau\right>$
is nonzero and since $S$ is weakly pseudoconvex at every point of
$S$, it follows that $k$ must be odd.

\medbreak We now introduce our second way of choosing $\xi$ with the
goal of deriving from it the conclusion that $k$ is even.  We choose
some smooth vector field $\left(\xi_1,\cdots,\xi_{n-1}\right)$ of
type $(1,0)$ on some open neighborhood $W_2$ of $0$ in $W$ such that
$\left<d\rho,\,\left(\xi_1,\cdots,\xi_{n-1}\right)\right>$ is
nowhere zero on $W_2$.  We now define $\xi_n$ on $W_2\times G$ by
$$
\xi_n=\frac{-1}{\phi+\left(\partial_n\phi\right)
z_n+\left(\partial_n\bar\phi\right)\overline{z_n}}
\left(\sum_{j=1}^{n-1}\left(\partial_j\phi\right)
\xi_j\,z_n+\sum_{j=1}^{n-1}\left(\partial_j\bar\phi\right)\xi_j\,\overline{z_n}
+k\rho^{k-1}\sum_{j=1}^{n-1}\left(\partial_j\rho\right)\xi_j\right)\leqno{({\rm
IV}.9.4)}
$$
so that the vector field
$\xi=\left(\xi_1,\cdots,\xi_{n-1},\xi_n\right)$ on $W_2\times G$
satisfies $\left<\partial r,\,\xi\right>\equiv 0$ because of
(IV.9.1). Since on $S\cap Y\cap U$ the two functions $z_n$ and
$\overline{z_n}$ are of the order $O\left(\rho^k\right)$, it follows
from (IV.9.4) that $\xi_n$ is of the order
$O\left(\rho^{k-1}\right)$ on $S\cap Y\cap\left(W_2\times G\right)$.
By (IV.9.2)
$$\left<\partial\bar\partial
r,\xi\wedge\bar\xi\right>=k(k-1)\rho^{k-2}\left<\partial\rho,\xi\right>\,\left<\bar\partial\rho,\bar\xi\right>
+O\left(\rho^{k-1}\right)$$ on $S\cap Y\cap\left(W_2\times
G\right)$. Since at the origin
$\left<d\rho,\,\left(\xi_1,\cdots,\xi_{n-1}\right)\right>$ is
nonzero and since $S$ is weakly pseudoconvex at every point of $S$,
it follows that $k$ must be even.  Thus we have a contradiction,
because earlier we have the conclusion that $k$ must be odd.

\bigbreak\noindent(IV.10) {\it Another Special Case to Illustrate
the Argument of Different Tangential and Normal Vanishing Orders for
Complex Hessian When Approaching CR Submanifold of Higher
Holomorphic Codimension.} We now consider another special case for
the more general situation where locally $M$ is defined by
real-valued real-analytic functions $\rho_1,\cdots,\rho_\ell$ on $V$
with $\ell>1$ and $\partial\rho_1,\cdots,\partial\rho_\ell$ are
${\mathbb C}$-linearly independent at points of $M$.  We use this
special case to further illustrate the argument of different
tangential and normal vanishing orders for the complex Hessian.  We
first explain what this special case is.

\medbreak As discussed above in (IV.8), there exist some $\tau\in
T^{(1,0)}_M$ such that
$\left(\partial\rho_j\right)\left(\tau\right)=0$ at $P_0$ for $1\leq
j\leq\ell$ and
$\left(\partial\bar\partial\rho_1\right)\left(\tau,\overline{\tau}\right)$
is nonzero but
$\left(\partial\bar\partial\rho_j\right)\left(\tau,\overline{\tau}\right)$
is zero for $2\leq j\leq\ell$.  We can write
$$
r|_V=\sum_{\nu_1+\cdots+\nu_\ell=k}\sigma_{\nu_1,\cdots,\nu_\ell}\left(\rho_1\right)^{\nu_1}\cdots\left(\rho_\ell\right)^{\nu_\ell}
+O\left(\left(\sum_{j=1}^\ell\left(\rho_j\right)^2\right)^{\frac{k+1}{2}}\right)
$$
for some integer $k\geq 2$, where $\sigma_{\nu_1,\cdots,\nu_\ell}$
is a real-analytic function on $U_0$ (after shrinking $U_0$ as an
open neighborhood of $P_0$ in ${\mathbb C}^n$ if necessary) and
$\sigma_{\nu_1^*,\cdots,\nu_\ell^*}$ is nonzero at $P_0$ for some
$\nu_1^*+\cdots+\nu_\ell^*=k$.  This special case which we now
consider is when $\sigma_{\nu_1,\cdots,\nu_\ell}$ is nonzero at
$P_0$ for some $\nu_1+\cdots+\nu_\ell=k$ with $\nu_1\not=0$.

\medbreak For this special case, just as for the case of $\ell=1$ we
can find a smooth vector field $\xi$ of type $(1,0)$ in some open
neighborhood of $P_0$ in ${\mathbb C}^n$ which are tangential to
$\partial\Omega$ such that the value of $\xi$ at $P_0$ agrees with
$\tau$. By computing the Levi form of $r$ at the vector field $\xi$
and its vanishing order at $M$ by using
$$
r|_V=\sum_{\nu_1+\cdots+\nu_\ell=k}\sigma_{\nu_1,\cdots,\nu_\ell}\left(\rho_1\right)^{\nu_1}\cdots\left(\rho_\ell\right)^{\nu_\ell}
+O\left(\left(\sum_{j=1}^\ell\left(\rho_j\right)^2\right)^{\frac{k+1}{2}}\right),
$$
as in the case of $\ell=1$ we can conclude that $k$ must be odd.
Thus we have a contradiction.  However, this argument depends on the
additional assumption that $\sigma_{\nu_1,\cdots,\nu_\ell}$ is
nonzero at $P_0$ for some $\nu_1+\cdots+\nu_\ell=k$ with
$\nu_1\not=0$ for a specially chosen set of defining functions
$\rho_1,\cdots,\rho_\ell$.

\medbreak Note that to rule out the case of an odd $k$, we do not
need this additional assumption that
$\sigma_{\nu_1,\cdots,\nu_\ell}$ is nonzero at $P_0$ for some
$\nu_1+\cdots+\nu_\ell=k$ with $\nu_1\not=0$.  There is also another
way to rule out the case of an odd $k$ by using bounded strictly
plurisubharmonic exhaustion functions for weakly pseudoconvex
domains in the following way.

\bigbreak\noindent(IV.11) {\it Handling the Case of Odd Vanishing
Order by Using Bounded Strictly Plurisubharmonic Exhaustion
Functions for Weakly Pseudoconvex Domains.}  First let us introduce
the following trivial statement about the vanishing order of a
negative subharmonic function at a boundary segment, which is
related to Hopf's lemma or the strong maximum principle [GT83, p.34,
Lemma 3.4].

\medbreak\noindent(IV.11.1) Let $D$ be a connected open subset of
${\mathbb C}$ and $C$ is a smooth connected curve in $D$ defined by
$\rho=0$ with $d\rho$ nowhere zero at points of $C$ such that $D-C$
consists of two nonempty components $W_1$ and $W_2$ with $\rho<0$ on
$W_1$. Let $\eta>1$ and $\varphi$ be a smooth negative subharmonic
function on $W_1$. Then it is impossible to write
$-\varphi=\left(-\rho\right)^\eta$ on $W_1$.

\medbreak The reason is as follows.  We compute
$$
\displaylines{\bar\partial\left(-\varphi\right)=-\eta\left(-\rho\right)^{\eta-1}\bar\partial\rho,
\cr
\partial\bar\partial\left(-\varphi\right)=\eta\left(\eta-1\right)\left(-\rho\right)^{\eta-2}\partial\rho\bar\partial\rho
-\eta\left(-\rho\right)^{\eta-1}\partial\bar\partial\rho.\cr}
$$
Since $\partial\bar\partial\varphi\geq 0$ on $W_1$, it follows that
$$
0\geq
\partial\bar\partial\left(-\varphi\right)=\eta\left(\eta-1\right)\left(-\rho\right)^{\eta-2}\partial\rho\bar\partial\rho
-\eta\left(-\rho\right)^{\eta-1}\partial\bar\partial\rho
$$
and
$$
\partial\bar\partial\rho\geq
\partial\bar\partial\left(-\varphi\right)=\frac{\eta-1}{-\rho}\,\partial\rho\bar\partial\rho,
$$
which is a contradiction, because the left-hand side evaluated at a
point of $W_1$ stays bounded as the point approaches some point of
$C$ but the right-hand side evaluated at the same point becomes
$\infty$ as the point approaches some point of $C$.

\medbreak We now recall the following theorem of Diederich-Fornaess
on bounded strictly plurisubharmonic exhaustion functions for weakly
pseudoconvex domains [DF77, p.133, Remark b].

\bigbreak\noindent Let $\Omega$ be a domain in ${\mathbb C}^n$ and
$P_0$ belong to the boundary of $\Omega$ so that for some open
neighborhood $D$ of $P_0$ in ${\mathbb C}^n$ the boundary of
$\Omega\cap D$ in $D$ is smooth and weakly pseudoconvex.  Let
$\delta$ be the distance function from a point of $\Omega$ to
${\mathbb C}^n-\Omega$. Let $\psi$ be a smooth strictly
plurisubharmonic function on ${\mathbb C}^n$ (or just defined on
some open neighborhood of some point of $\partial\Omega$ in
${\mathbb C}^n$).  Then for any choice of $0<\gamma<1$ there is a
suitable choice of a sufficiently small $L>0$ such that the complex
Hessian $\partial\bar\partial\left(-\delta^\gamma e^{-L\psi}\right)$
is strictly positive on $\Omega\cap D^\prime$ for some open
neighborhood $D^\prime$ of $P_0$ in $D$.

\medbreak Suppose we have the case of an odd $k$ in the following
expansion which we would like to rule out.
$$
r|_V=\sum_{\nu_1+\cdots+\nu_\ell=k}\sigma_{\nu_1,\cdots,\nu_\ell}\left(\rho_1\right)^{\nu_1}\cdots\left(\rho_\ell\right)^{\nu_\ell}
+O\left(\left(\sum_{j=1}^\ell\left(\rho_j\right)^2\right)^{\frac{k+1}{2}}\right)
$$
for some integer $k\geq 2$, where $\sigma_{\nu_1,\cdots,\nu_\ell}$
is a real-analytic function on $U_0$ (after shrinking $U_0$ as an
open neighborhood of $P_0$ in ${\mathbb C}^n$ if necessary) and
$\sigma_{\nu_1^*,\cdots,\nu_\ell^*}$ is nonzero at $P_0$ for some
$\nu_1^*+\cdots+\nu_\ell^*=k$.  Assume that $k$ is odd.  We can find
a tangent vector $\eta$ of $V$ at the point $P_0$ of $M$ normal to
$M$ such that $J\eta$ is tangential to $M$ and the $k$-derivative of
$r$ in the direction of $\eta$ is nonzero.  Let $C$ be local complex
curve in ${\mathbb C}^n$ through $P_0$ such that the complex tangent
vector to $C$ of type $(1,0)$ at $P_0$ is equal to $\eta-\sqrt{-1} J\eta$ and
$C\cap\partial\Omega\subset M$ and $C\cap M$ is a regular curve in
$C$. Since $k$ is odd, after replacing $C$ by an open neighborhood
of $P_0$ in $C$ we can assume without loss of generality that $C-M$
consists of two nonempty connected components $C\cap\Omega$ and
$C-\overline\Omega$.

\medbreak Let $\kappa=-\delta^\gamma e^{-L\psi}$ and we restrict
$\kappa$ to $C\cap\Omega$.  Let $\phi$ be a smooth function on $C$
whose zero-set is $C\cap M$ and which is negative on $C\cap\Omega$
with $d\phi$ nowhere zero on $C\cap M$. Since $-\kappa$ is equal to
$\sigma\left(-r\right)^\gamma=\tilde\sigma\left(-\phi\right)^{k\gamma}$
on $C\cap\Omega$ for some positive-valued smooth functions $\sigma$
and $\tilde\sigma$ on $C$ (after replacing $C$ by an open
neighborhood of $P_0$ in $C$ if necessary), from the
plurisubharmonicity of $\kappa$ on $\Omega$ we have a contradiction
to (IV.11.1) when $0<\gamma<1$ is chosen to satisfy $k\gamma>1$,
because $\kappa|_{C\cap\Omega}$ is subharmonic on $C\cap\Omega$ and
$-\kappa|_{C\cap\Omega}$ is equal to
$\left(-\tilde\sigma^{\frac{1}{k\gamma}}\phi\right)^{k\gamma}$ and
$\tilde\sigma^{\frac{1}{k\gamma}}\phi$ is smooth on $C$ and is $0$
at $C\cap M$ and $d\tilde\sigma^{\frac{1}{k\gamma}}\phi$ is nowhere
zero on $M$. This argument avoids the process in (IV.10) of
constructing the analog of the second vector field, at the end of
(IV.9), of type $(1,0)$ in a neighborhood of $P_0$ in ${\mathbb
C}^n$ tangential to $S$ and not tangential to $M$ at $P_0$.

\bigbreak\noindent(IV.12) {\it Handling the Case of Even Vanishing
Order by Stratification According to Iterated Lie Brackets.} We now
deal with the general case by choosing the set of defining functions
$\rho_1,\cdots,\rho_\ell$ by stratification according to iterated
Lie brackets.  Recall that iterated Lie brackets of vector fields on
$E_\infty$ with coefficients in ${\mathcal N}^{(\infty)}={\mathcal
N}\cap T^{\mathbb R}_{E_\infty}\cap JT^{\mathbb R}_{E_\infty}$
generate the distribution $\tilde{\mathcal N}$ and $M$ is an
integral submanifold of $E_\infty$ whose tangent space at every
point is equal to the subspace distribution $\tilde{\mathcal N}$ at
that point.  Because of the Jacobi identity for the Lie brackets of
three vector fields, we can select vector fields
$\tau_0,\tau_1,\cdots,\tau_\ell$ on $E_\infty$ with values in
${\mathcal N}^{(\infty)}={\mathcal N}\cap T^{\mathbb
R}_{E_\infty}\cap JT^{\mathbb R}_{E_\infty}$ such that inductively,
$\tilde\tau_1=\left[\tau_0,\tau_1\right]$ and
$\tilde\tau_j=\left[\tilde\tau_{j-1},\tau_j\right]$ for $2\leq
j\leq\ell$ and $\tilde\tau_j\left(P_0\right)$ is not spanned by
$\left({\mathcal N}^{(\infty)}\right)_{P_0}$, $\tilde
\tau_1\left(P_0\right),\,\cdots,\, \tilde\tau_{j-1}\left(P_0\right)$
for $1\leq j\leq\ell$.  We now choose $\rho_1,\cdots,\rho_\ell$ such
that, modulo $\left({\mathcal N}^{(\infty)}\right)_{P_0}$, the
$1$-forms $\left(Jd\rho_j\right)\left(P_0\right)$ at $P_0$ for
$1\leq j\leq\ell$, when restricted to the tangent space $T^{\mathbb
R}_{M,P_0}$ of $M$ at $P_0$ precisely form a dual basis for
$\tilde\tau_1\left(P_0\right),\cdots,\tilde\tau_\ell\left(P_0\right)$.
In other words, the ${\mathbb R}$-linear functionals defined by
$\left(Jd\rho_j\right)\left(P_0\right)$ at $P_0$ for $1\leq
j\leq\ell$ on the quotient space $T^{\mathbb
R}_{M,P_0}\left/\left({\mathcal N}^{(\infty)}\right)_{P_0}\right.$
form the dual basis of the elements in $T^{\mathbb
R}_{M,P_0}\left/\left({\mathcal N}^{(\infty)}\right)_{P_0}\right.$
induced by
$\tilde\tau_1\left(P_0\right),\cdots,\tilde\tau_\ell\left(P_0\right)$.

\medbreak  Let
$\xi_j=\frac{1}{2}\left(\tilde\tau_{j-1}-\sqrt{-1}J\tilde\tau_{j-1}\right)$
and $\eta_j=\frac{1}{2}\left(\tau_j-\sqrt{-1}J\tau_j\right)$ on $M$ for
$1\leq j\leq\ell$ so that both $\xi_j$ and $\eta_j$ are of type
$(1,0)$ tangential to $V$ with the real part of $\xi_j$ being
$\frac{1}{2}\tilde\tau_j$ and the real part of $\eta_j$ being
$\frac{1}{2}\tau_j$. Take $2\leq j\leq\ell$.  As verified above in
(IV.7), from $\tilde\tau_{j-1}=\xi_j+\overline{\xi_j}$ and
$\tau_j=\eta_j+\overline{\eta_j}$ we get
$$\left[\tilde\tau_{j-1},\tau_j\right]=\frac{\sqrt{-1}}{2}
\left(\left[\xi_j+\sqrt{-1}\eta_j,\overline{\xi_j+\sqrt{-1}\eta_j}\right]-
\left[\xi_j-\sqrt{-1}\eta_j,\overline{\xi_j-\sqrt{-1}\eta_j}\right]\right)$$ modulo
${\mathbb C}\otimes_{\mathbb R}\left(T^{\mathbb R}_M\cap JT^{\mathbb
R}_M\right)=T^{(1,0)}_M\oplus T^{(0,1)}_M$, where
$$
\displaylines{\xi_j+\sqrt{-1}\eta_j=\frac{1}{2}\left(\tilde\tau_{j-1}+J\tau_j\right)
+\frac{\sqrt{-1}}{2}\left(\tau_j-J\tilde\tau_{j-1}\right),\cr
\xi_j-\sqrt{-1}\eta_j=\frac{1}{2}\left(\tilde\tau_{j-1}-J\tau_j\right)
-\frac{\sqrt{-1}}{2}\left(\tau_j+J\tilde\tau_{j-1}\right).\cr}
$$
At the point $P_0$ we have
$$
\displaylines{1=\left(Jd\rho_j\right)\left(\tilde\tau_j\right)\cr=
\frac{\sqrt{-1}}{2}\left(Jd\rho_j\right)\left(\left[\xi_j+\sqrt{-1}\eta_j,\overline{\xi_j+\sqrt{-1}\eta_j}\right]-
\left[\xi_j-\sqrt{-1}\eta_j,\overline{\xi_j-\sqrt{-1}\eta_j}\right]\right)\cr}
$$
and at least one of
$$
\left(Jd\rho_j\right)\left(\left[\xi_j+\sqrt{-1}\eta_j,\overline{\xi_j+\sqrt{-1}\eta_j}\right]\right)
\quad{\rm
and}\quad \left(Jd\rho_j\right)\left(
\left[\xi_j-\sqrt{-1}\eta_j,\overline{\xi_j-\sqrt{-1}\eta_j}\right]\right)
$$
has absolute value at least $1$ and is nonzero at $P_0$.  We set
$\zeta_j$ to be one of the two possibilities
$$
\displaylines{\xi_j+\sqrt{-1}\eta_j=\frac{1}{2}\left(\tilde\tau_{j-1}+J\tau_j\right)
+\frac{\sqrt{-1}}{2}\left(\tau_j-J\tilde\tau_{j-1}\right),\cr
\xi_j-\sqrt{-1}\eta_j=\frac{1}{2}\left(\tilde\tau_{j-1}-J\tau_j\right)
-\frac{\sqrt{-1}}{2}\left(\tau_j+J\tilde\tau_{j-1}\right).\cr}
$$
so that
$\left|\left(Jd\rho_j\right)\left(\left[\zeta_j,\overline{\zeta_j}\right]\right)\right|\geq
1$ at $P_0$.  From the way we define the $1$-jet of $\rho_j$ at
$P_0$ we know that among the following vectors
$$
\displaylines{ \tau_0,\tau_1,\cdots,\tau_\ell,
J\tau_0,J\tau_1,\cdots,J\tau_\ell,\cr
\tilde\tau_1,\cdots,\tilde\tau_\ell,
J\tilde\tau_1,\cdots,J\tilde\tau_\ell\cr}
$$
at $P_0$ the only one at which $d\rho_j$ is nonzero is
$J\tilde\tau_j$ where the value of $d\rho_j$ is $1$, because the
vectors
$$
\tau_0,\tau_1,\cdots,\tau_\ell, J\tau_0,J\tau_1,\cdots,J\tau_\ell,
\tilde\tau_1,\cdots,\tilde\tau_\ell
$$
all belong to the tangent space $T^{\mathbb R}_{M,P_0}$ of $M$ at
$P_0$ which is equal to $\left({\mathcal N}^{(\infty)}\right)_{P_0}$
and $\rho_j$ vanishes on $M$ and because the ${\mathbb R}$-linear
functionals defined by
$$\left(Jd\rho_1\right)\left(P_0\right),\cdots,\left(Jd\rho_\ell\right)\left(P_0\right)$$
at $P_0$ on the quotient space $T^{\mathbb
R}_{M,P_0}\left/\left({\mathcal N}^{(\infty)}\right)_{P_0}\right.$
form the dual basis of the elements in $T^{\mathbb
R}_{M,P_0}\left/\left({\mathcal N}^{(\infty)}\right)_{P_0}\right.$
induced by
$\tilde\tau_1\left(P_0\right),\cdots,\tilde\tau_\ell\left(P_0\right)$.
Hence $\left(d\rho_j\right)\left(\zeta_p\right)$ at $P_0$ is $0$ for
$j\not=p-1$.

\medbreak For later use we need a slight variation of the above
discussion. Take a positive number $A$. From
$\tilde\tau_{j-1}=\xi_j+\overline{\xi_j}$ and
$\tau_j=\eta_j+\overline{\eta_j}$ we get
$$A\tilde\tau_j=\left[\tilde\tau_{j-1},A\tau_j\right]=\frac{1}{2}\left(\left[\xi_j+\sqrt{-1}A\eta_j,\overline{\xi_j+iA\eta_j}\right]-
\left[\xi_j-\sqrt{-1}A\eta_j,\overline{\xi_j-\sqrt{-1}A\eta_j}\right]\right),$$
where
$$
\displaylines{\xi_j+\sqrt{-1}A\eta_j=\frac{1}{2}\left(\tilde\tau_{j-1}+AJ\tau_j\right)
-\frac{\sqrt{-1}}{2}\left(\tilde\tau_{j-1}+AJ\tau_j\right),\cr
\xi_j-\sqrt{-1}A\eta_j=\frac{1}{2}\left(\tilde\tau_{j-1}-AJ\tau_j\right)
-\frac{\sqrt{-1}}{2}\left(\tilde\tau_{j-1}-AJ\tau_j\right).\cr}
$$
At the point $P_0$ we have
$$
\displaylines{A=A\left(Jd\rho_j\right)\left(
\tilde\tau_j\right)=\left(Jd\rho_j\right)\left(
A\tilde\tau_j\right)=\left(Jd\rho_j\right)\left(\left[\tilde\tau_{j-1},A\tau_j\right]\right)
\cr=
\frac{\sqrt{-1}}{2}\left(Jd\rho_j\right)\left(\left[\xi_j+\sqrt{-1}A\eta_j,\overline{\xi_j+\sqrt{-1}A\eta_j}\right]-
\left[\xi_j-\sqrt{-1}A\eta_j,\overline{\xi_j-\sqrt{-1}A\eta_j}\right]\right)\cr}
$$
and at least one of
$$
\left(Jd\rho_j\right)\left(
\left[\xi_j+\sqrt{-1}A\eta_j,\overline{\xi_j+\sqrt{-1}A\eta_j}\right]\right)\quad{\rm
and}\quad \left(Jd\rho_j\right)\left(
\left[\xi_j-\sqrt{-1}A\eta_j,\overline{\xi_j-\sqrt{-1}A\eta_j}\right]\right)
$$
has absolute value at least $A$ at $P_0$.  We set $\zeta_{j,A}$ to
be one of the two possibilities
$$
\displaylines{\xi_j+\sqrt{-1}A\eta_j=\frac{1}{2}\left(\tilde\tau_{j-1}+AJ\tau_j\right)
+\frac{\sqrt{-1}}{2}\left(A\tau_j-J\tilde\tau_{j-1}\right),\cr
\xi_j-\sqrt{-1}A\eta_j=\frac{1}{2}\left(\tilde\tau_{j-1}-AJ\tau_j\right)
-\frac{\sqrt{-1}}{2}\left(A\tau_j+J\tilde\tau_{j-1}\right)\cr}
$$
so that
$\left|\left(Jd\rho_j\right)\left(\left[\zeta_{j,A},\overline{\zeta_{j,A}}\right]\right)\right|\geq
A$ at $P_0$.  Note that when $A=1$ we have $\zeta_{j,A}=\zeta_j$ so
that for any value of $A>1$ the vector $\zeta_{j,A}-\zeta_j$ is
equal to
$$\pm\left(\frac{(A-1)}{2}J\tau_j+\frac{\sqrt{-1}(A-1)}{2}\tau_j\right).$$ From the way we define the $1$-jet
of $\rho_j$ at $P_0$ we know that among the following vectors
$$
\displaylines{ \tau_0,\tau_1,\cdots,\tau_\ell,
J\tau_0,J\tau_1,\cdots,J\tau_\ell,\cr
\tilde\tau_1,\cdots,\tilde\tau_\ell,
J\tilde\tau_1,\cdots,J\tilde\tau_\ell\cr}
$$
at $P_0$ the only one at which $d\rho_j$ is nonzero is
$J\tilde\tau_j$ where the value of $d\rho_j$ is $1$, because the
vectors
$$
\tau_0,\tau_1,\cdots,\tau_\ell, J\tau_0,J\tau_1,\cdots,J\tau_\ell,
\tilde\tau_1,\cdots,\tilde\tau_\ell
$$
all belong to the tangent space $T^{\mathbb R}_{M,P_0}$ of $M$ at
$P_0$ which is equal to $\left({\mathcal N}^{(\infty)}\right)_{P_0}$
and $\rho_j$ vanishes on $M$ and because the ${\mathbb R}$-linear
functionals defined by
$$\left(Jd\rho_1\right)\left(P_0\right),\cdots,\left(Jd\rho_\ell\right)\left(P_0\right)$$
at $P_0$ on the quotient space $T^{\mathbb
R}_{M,P_0}\left/\left({\mathcal N}^{(\infty)}\right)_{P_0}\right.$
form the dual basis of the elements in $T^{\mathbb
R}_{M,P_0}\left/\left({\mathcal N}^{(\infty)}\right)_{P_0}\right.$
induced by
$\tilde\tau_1\left(P_0\right),\cdots,\tilde\tau_\ell\left(P_0\right)$.
Hence $\left(d\rho_j\right)\left(\zeta_p\right)$ at $P_0$ is $0$ for
$j\not=p-1$. Moreover, at $P_0$
$$\left(d\rho_j\right)\left(\zeta_{p,A}-\zeta_p\right)=
\left(d\rho_j\right)\left(\pm\left(\frac{(A-1)}{2}J\tau_j+\frac{\sqrt{-1}(A-1)}{2}\tau_j\right)\right)=0.
$$
Thus
$\left(d\rho_j\right)\left(\zeta_{p,A}\right)=\left(d\rho_j\right)\left(\zeta_p\right)$
is independent of $A$ for all $p$ and $j$.

\medbreak Let $1<q\leq\ell$ be the minimum such that $\nu_q\geq 1$
and $\nu_q+\cdots+\nu_\ell=k$ and
$\sigma_{0,\cdots,0,\nu_q,\cdots,\nu_\ell}\left(P_0\right)\not=0$.
Since we have $\zeta_{q,A}\rho_j=0$ at $P_0$ for $j\not=q-1$, the
term of lowest vanishing order at $P_0$ which we can get is $k-1$
and either come with the factors
$\zeta_{q,A}\overline{\zeta_{q,A}}\rho_j$ at $P_0$ from
$$
\sum_{\nu_1+\cdots+\nu_\ell=k}\sigma_{\nu_1,\cdots,\nu_\ell}\left(\rho_1\right)^{\nu_1}\cdots\left(\rho_\ell\right)^{\nu_\ell}
$$
or come with the factor
$\zeta_{q,A}\rho_{q-1}\overline{\zeta_{q,A}}\rho_{q-1}$ at $P_0$
from
$$
\sum_{\nu_1+\cdots+\nu_\ell=k+1}\sigma_{\nu_1,\cdots,\nu_\ell}\left(\rho_1\right)^{\nu_1}\cdots\left(\rho_\ell\right)^{\nu_\ell}
$$
in the expansion of $r|_V$.  The sum of all these terms of vanishing
order $k-1$ at $P_0$ is identically zero only when
$$
\displaylines{\sum_{\nu_q+\cdots+\nu_\ell=k}\sigma_{0,\cdots,0,\nu_q,\cdots,\nu_\ell}\nu_q\left(\zeta_{q,A}
\overline{\zeta_{q,A}}
\rho_q\right)\left(\rho_q\right)^{\nu_q-1}\left(\rho_{q+1}\right)^{\nu_{q+1}}\cdots\left(\rho_\ell\right)^{\nu_\ell}\cr
=-\,\sum_{\nu_1+\cdots+\nu_\ell=k+1}\sigma_{\nu_1,\cdots,\nu_\ell}\nu_{q-1}\left(\nu_{q-1}-1\right)
\left(\zeta_{q,A}\rho_{q-1}\right)\left(\overline{\zeta_{q,A}}\rho_{q-1}\right)\cdot\cr
\hfill\cdot\left(\left(
\rho_1\right)^{\nu_1}\cdots\left(\rho_{q-2}\right)^{\nu_{q-2}}
\left(\rho_{q-1}\right)^{\nu_{q-1}-2}\left(\rho_{q+1}\right)^{\nu_{q+1}}
\cdots\left(\rho_\ell\right)^{\nu_\ell}\right).\cr}
$$
Since the contradiction comes from the change of the sign of the
Levi Form of $S=\partial\Omega$ when we approach $M$ from along some
appropriate path in $S=\partial\Omega$ which corresponds to a path
in $V$ up to order $k$, we have trouble only when for any choice of
$\nu_q\geq 1,\nu_{q+1}\geq 0,\cdots,\nu_\ell\geq 0$ with
$\nu_q+\cdots+\nu_\ell=k$ we have
$$
\displaylines{\qquad
\sigma_{0,\cdots,0,\nu_q,\cdots,\nu_\ell}\nu_q\left(\zeta_{q,A}\overline{\zeta_{q,A}}
\rho_q\right)\hfill\cr\hfill=-\,\sum_{1\leq i\leq
j<q}\sigma_{0,\cdots,0,\nu_{q-1}+2,\nu_q-1,\nu_{q+1}\cdots,\nu_\ell}\nu_{q-1}\left(\nu_{q-1}-1\right)
\left(\zeta_{q,A}\rho_{q-1}\right)
\left(\overline{\zeta_{q,A}}\rho_{q-1}\right)\cr}
$$
at the point $P_0$ of $M$ for any choice of $\nu_q\geq
1,\nu_{q+1}\geq 0,\cdots,\nu_\ell\geq 0$ with
$\nu_q+\cdots+\nu_\ell=k$, which is the same as
$$
\displaylines{\qquad
A^2\sigma_{0,\cdots,0,\nu_q,\cdots,\nu_\ell}\nu_q\left(\zeta_q\overline{\zeta_q}
\rho_q\right)\hfill\cr\hfill=-\,\sum_{1\leq i\leq
j<q}\sigma_{0,\cdots,0,\nu_{q-1}+2,\nu_q-1,\nu_{q+1}\cdots,\nu_\ell}\nu_{q-1}\left(\nu_{q-1}-1\right)
\left(\zeta_q\rho_{q-1}\right)
\left(\overline{\zeta_q}\rho_{q-1}\right),\cr}
$$
because $\left(\zeta_{q,A}\overline{\zeta_{q,A}}
\rho_q\right)\left(P_0\right)=A^2\left(\zeta_q\overline{\zeta_q}
\rho_q\right)\left(P_0\right)$ and
$\left(\zeta_{q,A}\rho_{q-1}\right)\left(P_0\right)=\left(\zeta_q\rho_{q-1}\right)\left(P_0\right)$
and
$\left(\overline{\zeta_{q,A}}\rho_{q-1}\right)\left(P_0\right)=\left(\overline{\zeta_q}\rho_{q-1}\right)\left(P_0\right)$.
Since $\left(\zeta_q\overline{\zeta_q}
\rho_q\right)\left(P_0\right)\not=0$, this trouble can simply be
handled with the choice of a sufficiently large $A$.

\medbreak Finally, in order to get a contradiction from the evenness
of the vanishing order $k$ of $r|_{V\cap U}$ at $M\cap U$, we
construct
\begin{itemize}
\item[(i)] a real-analytic
curve $\Gamma_S$ in $S\cap U$ containing $P_0$ which is transversal
to $M$,

\item[(ii)] a real-analytic
curve $\Gamma_V$ in $V\cap U$ containing $P_0$ which is transversal
to $M$,

\item[(iii)] a smooth bijection $\Psi$ from $\Gamma_V$ to
$\Gamma_S$,

\item[(iv)] a vector field $\zeta_S$ of type $(1,0)$ tangential to $S$ defined
only at points of the curve $\Gamma_S$ and smooth along $\Gamma_S$
whose value at $P_0$ is $\zeta_{q,A}$, and

\item[(v)] a vector field $\zeta_V$ of type $(1,0)$ tangential to $V$ defined
only at points of the curve $\Gamma_V$ and smooth along $\Gamma_V$
whose value at $P_0$ is $\zeta_{q,A}$
\end{itemize}
such that
\begin{itemize}
\item[(a)]
the distance between $P\in\Gamma_V$ and $\Psi(P)\in\Gamma_S$ is of
order of $\left({\rm dist}_{\Gamma_V}\left(P,P_0\right)\right)^k$,
where ${\rm dist}_{\Gamma_V}\left(P,P_0\right)$ is the distance
between $P$ and $P_0$ along $\Gamma_V$, and
\item[(a)]
the difference of the value of $\zeta_V$ at $P\in\Gamma_V$ and the
value of $\zeta_S$ at $\Psi(P)\in\Gamma_S$ is of order of
$\left({\rm dist}_{\Gamma_V}\left(P,P_0\right)\right)^k$.
\end{itemize}
Then the Levi form of $r$ in the direction $\zeta_S$ at a point $P$
in $\Gamma_S$ other than $P_0$ will change sign as $P$ moves along
$\Gamma_S$ to pass $P_0$ because the evenness of $k$ implies that
the Levi form of $r$ in the direction $\zeta_S$ at a point $P$
vanishes of odd order $k-1$ at $P_0$ along $\Gamma_S$. This
contradicts the weak pseudoconvexity of $S$.

\eject\noindent{\bf Appendix A: Some Techniques of Applying of
Skoda's Theorem on Ideal Generation}

\bigbreak In this Appendix we give some techniques of applying
Skoda's theorem on ideal generation [Sk72, Th.1, pp.555-556] which
involve derivatives and Jacobian determinants. The significance is
more in the techniques themselves than in the statements given here
to demonstrate their use.  Though these techniques are not directly
used in this note (except the use of (A.2) in (III.7) and the use of (A.3) in (III.8)), they may be useful in reducing the vanishing orders
of multiplier ideals in Kohn-type algorithms in the setting of more
general partial differential equations.

\bigbreak\noindent(A.1) {\it Proposition.}  Let $\Omega$ be a
bounded Stein open subset of ${\mathbb C}^n$.  Let $g_1,\cdots,g_n,
\rho$ be holomorphic functions on some open neighborhood
$\tilde\Omega$ of the topological closure $\bar\Omega$ of $\Omega$.
Let $Z$ be the common zero-set of $g_1,\cdots,g_n$ in
$\tilde\Omega$.  Assume that $\rho$ vanishes on $Z$.  Let $J$ be the
Jacobian determinant of $g_1,\cdots,g_n$. Then there exist
holomorphic $h_1,\cdots,h_n$ on $\Omega$ such that $\rho
J=\sum_{j=1}^n h_jg_j$.

\medbreak\noindent{\it Proof.}  For any $0<\gamma<1$ and any compact subset $K$ of ${\mathbb C}^n$ with coordinates $w=\left(w_1,\cdots,w_n\right)$ the integral
$$
\int_{w\in K}\frac{\prod_{j=1}^n\left(\sqrt{-1}dw_j\wedge d\overline{w_j}\right)}{\left(\sum_{j=1}^n\left|w_j\right|^2\right)^{\gamma n}}\leqno{({\rm A}.1.1)}
$$
is finite.
Since $\rho$ vanishes on $Z$, it follow that there exists some $0<\eta<1$ such that
$$
\frac{\left|\rho\right|^2}{\left(\sum_{j=1}^n\left|g_j\right|^2\right)^{\eta n}}\leqno{({\rm A}.1.2)}
$$
is bounded on some open neighborhood $U$ of $\bar\Omega$ in $\tilde\Omega$.  Let $\gamma=1-\frac{\eta}{2}$ and $\alpha=1+\frac{\eta}{2}$.  Since $J$ is the Jacobian
determinant of $g_1,\cdots,g_n$, by pulling back
(A.1.1) by $w_j=g_j$ for $1\leq j\leq n$ and using the uniform boundedness of (A.1.2) on $U$, we conclude that
$$
\int_\Omega\frac{\left|\rho\,J\right|^2}{\left(\sum_{j=1}^n\left|g_j\right|^2\right)^{\alpha
n}}<\infty.\leqno{({\rm A}.1.3)}
$$
By using (A.1.3) and applying Skoda's theorem [Sk72, Th.1,
pp.555-556] to the Stein domain $\Omega$ to
express $\rho\,J$ as a linear combination of $g_1,\cdots,g_n$ with
holomorphic functions, we obtain $h_1,\cdots,h_n$ satisfying the
requirements of the Proposition.  Q.E.D.

\bigbreak\noindent(A.2) {\it Proposition (Ideal Generated by Components of
Gradient).} Let $f$ be a holomorphic
function germ on ${\mathbb C}^n$ at the origin which vanishes at the
origin. Then $f^{n+1}$ belongs to the ideal ${\mathcal I}$ generated
by $\frac{\partial f}{\partial z_j}$ for $1\leq j\leq n$ at the
origin, where $z_1,\cdots,z_n$ are the coordinates of ${\mathbb
C}^n$.

\medbreak\noindent{\it Proof.}  Let $\pi:\tilde U\to U$ be the
simultaneous resolution of singularities for the ideal ${\mathcal
I}$ and the ideal ${\mathcal O}_{{\mathbb C}^n}\,f$ generated by $f$
with exceptional hypersurfaces $\left\{E_\ell\right\}_\ell$ in
normal crossing in $\tilde U$, where $U$ is an open neighborhood of
the origin in ${\mathbb C}^n$ on which the holomorphic function germ
$f$ is defined. We claim that
$$
\frac{\left|f\right|^2}{\sum_{j=1}^n\left|\frac{\partial f}{\partial
z_j}\right|^2}\leqno{({\rm A}.2.1)}
$$
is uniformly bounded in some relatively compact open neighborhood
$U^\prime$ of the origin in $U$. Otherwise, when we write the
divisor of $\pi^*f$ of $f$ as $\sum_\ell a_\ell E_\ell$ and write
$\pi^*{\mathcal I}$ as $\sum_\ell b_\ell E_\ell$ with $a_\ell$ and
$b_\ell$ being nonnegative integers, we have $b_\ell>a_\ell$ for
some $\ell$ with $0\in\pi\left(E_\ell\right)$ and we can find a
local holomorphic curve $\tilde\varphi:W\to\tilde U$ with $W$ being
an open neighborhood of the origin in ${\mathbb C}$ and
$\pi\tilde\varphi(0)=0$ such that $\varphi(W)$ is transversal to
$E_\ell$ and is disjoint from any $E_k$ with $k\not=\ell$.  Then
$d\left(f\circ\varphi\right)$ vanishes at $0$ to an order higher
than that $f\circ\varphi$, which is a contradiction, because
$f\circ\varphi$ vanishes at $0$.  This argument actually gives a
slightly higher vanishing order of $\left|f\right|^2$ than that of
$\sum_{j=1}^n\left|\frac{\partial f}{\partial z_j}\right|^2$ along
each $E_\ell$ when they are pulled back to $\tilde U$ so that
$$
\int_{U^\prime}\frac{\left|f^{n+1}\right|^2}{\left(\sum_{j=1}^n\left|\frac{\partial
f}{\partial z_j}\right|^2\right)^{\alpha(n+1)}}<\infty
$$
for some $\alpha>1$.  The conclusion of the Proposition now follows
from Skoda's theorem [Sk72, Th.1, pp.555-556].  Q.E.D.

\bigbreak\noindent(A.2.2) {\it Remark on the Relation Between Proposition (A.2) and l'H\^opital's Rule.}  The argument in the proof of Proposition (A.2) consists of the verification of the uniform bound of (A.2.1) on some open neighborhood $U^\prime$ of the origin in ${\mathbb C}^n$ and a straightforward application of Skoda's theorem [Sk72, Th.1, pp.555-556].  The argument used in the verification of the uniform bound of (A.2.1) on $U^\prime$ is actually the usual l'H\^opital's rule in calculus applied to the pullback of the quotient (A.2.1) to the open unit $1$-disk $\Delta$ in ${\mathbb C}$ by a holomorphic map $g:\Delta\to{\mathbb C}^n$ with $g(0)=0$ when one applies differentiation at the origin along any ray of $\Delta$ until one ends up with a nonzero derivative of the denominator.  The uniformity of the bound comes from the fact that one needs only consider a compact holomorphic family of such holomorphic maps $g:\Delta\to{\mathbb C}^n$, as is easily seen, for example, by using a resolution of singularities.  Another simple way of looking at the bound of (A.2.1) is the trivial observation that the vanishing order of an analytic function at a point of its zero-set is no more than the vanishing order of its gradient.

\bigbreak\noindent(A.2.3) {\it Remark on the Difference Between the Jacobian Determinants with Respect to All Variables and the Jacobian Determinants With Respect to All Variables with Respect to a Proper Subset of Variables.}  Let $F_1,\cdots,F_N$ be holomorphic function germs on ${\mathbb C}^2$ at the origin vanishing at the origin such that the ideal $I_1$ generated by $F_1,\cdots,F_N$ contains an effective power of the maximum ideal sheaf ${\mathfrak m}_{{\mathbb C}^2,0}$ of ${\mathbb C}^2$.  By Proposition (A.2) the ideal generated by the components of the gradients of $F_1,\cdots,F_N$, namely by $\frac{\partial F_j}{\partial z_k}$ for $1\leq j\leq N, 1\leq k\leq 2$, contains an effective power of ${\mathfrak m}_{{\mathbb C}^2,0}$.  We can regard each $\frac{\partial F_j}{\partial z_k}$ for $1\leq j\leq N, 1\leq k\leq 2$ as the Jacobian determinant of the single function $F_j$ with respect to the single variable $z_j$.  These first-order partial derivatives can be regarded as the Jacobian determinants with respect to a proper subset of all the variables.  Proposition (A.2) can be restated as follows.  The ideal generated by all such Jacobian determinants with respect to a proper subset of all the variables contains an effective power of ${\mathfrak m}_{{\mathbb C}^2,0}$.  The situation is very different from the ideal $I_2$ generated by all Jacobian determinants with respect to the full set of all the variables
$$
\frac{\partial\left(F_{j_1},F_{j_2}\right)}
{\partial\left(z_1,z_2\right)}\quad{\rm for\ \ }1\leq j_1,j_2\leq N.
$$
In general, the ideal $I_2$ does not contain an effective power of ${\mathfrak m}_{{\mathbb C}^2,0}$, as one can easily see in the special case where $N=2$ and the ideal $I_2$ is generated by a single holomorphic function germ.

\medbreak In general, for the complex Euclidean space ${\mathbb C}^n$ instead of ${\mathbb C}^2$, when we have holomorphic function germs $F_1,\cdots,F_N$ on ${\mathbb C}^n$ at the origin vanishing at the origin such that the ideal generated by $F_1,\cdots,F_N$ contains an effective power of the maximum ideal sheaf ${\mathfrak m}_{{\mathbb C}^n,0}$ of ${\mathbb C}^n$, we can consider for $1\leq\nu\leq n$ the ideal $I_\nu$ generated by the Jacobian determinants
$$
\frac{\partial\left(F_{j_1},\cdots,F_{j_\nu}\right)}
{\partial\left(z_{k_1},\cdots,z_{k_\nu}\right)}
$$
for $1\leq j_1,\cdots,j_\nu\leq N$ and $1\leq k_1,\cdots,k_\nu\leq n-1$.  As we see in Proposition (A.3) below, for $1\leq\nu\leq n-1$ the ideal $I_\nu$ contains an effective power of ${\mathfrak m}_{{\mathbb C}^n,0}$, though in general the ideal $I_n$ does not contain an effective power of ${\mathfrak m}_{{\mathbb C}^n,0}$.  It means that the situation for the ideal generated by
all Jacobian determinants with respect to a proper subset of all the variables is very different from the ideal $I_2$ generated by all Jacobian determinants with respect to the full set of all the variables.

\bigbreak\noindent(A.2.4) {\it Remark on a Generalization of the Special Case of Proposition (A.2) for Dimension Two.}  The special case of Proposition (A.2) for dimension two is used in this note in (III.7) to prove the effective termination of Kohn's algorithm for ${\mathbb C}^2$.  For the proof of the effective termination of Kohn's algorithm for ${\mathbb C}^n$ the corresponding statement which has to be used is not Proposition (A.2) for dimension $n$, but the following Proposition (A.3).

\bigbreak\noindent(A.3) {\it Proposition (Ideal Generated by Jacobian Determinants with Respect to a Proper Subset of Variables).}  Let $F_1,\cdots,F_N$ be holomorphic function germs on ${\mathbb C}^n$ at the origin vanishing at the origin such that the ideal generated by $F_1,\cdots,F_N$ contains an effective power of the maximum ideal of ${\mathbb C}^n$ at the origin.  Let $1\leq\nu<n$.  Let $J_\nu$ be the ideal generated by
$$
\frac{\partial\left(F_{j_1},\cdots,F_{j_\nu}\right)}{\partial\left(z_{k_1},\cdots,z_{k_\nu}\right)}.
$$
for $1\leq j_1<\cdots<j_\nu\leq N$ and $1\leq k_1<\cdots<k_\nu\leq n$.  Then the ideal $J_\nu$ contains an effective power of the maximum ideal of ${\mathbb C}^n$ at the origin.

\medbreak\noindent{\it Proof.}  Let us first introduce some notations.  For an ideal $I$ of ${\mathcal O}_{{\mathbb C}^n,0}$ we define
$$\left|s_I\right|=\left(\sum_{j=1}^{k_I}\left|s_{j,I}\right|^2\right)^{\frac{1}{2}},$$
where $s_{1,I},\cdots,s_{k_I,I}$ form a set of generators of $I$.  The expression $\left|s_I\right|$ is defined up to a choice of the set of generators.  We use this expression only in the context of determining whether one such expression is dominated by a constant times another such expression $\left|s_J\right|$ for another ideal $J$ of ${\mathcal O}_{{\mathbb C}^n,0}$.  For such a purpose the choices of generators in the definitions for
$\left|s_I\right|$ and $\left|s_J\right|$
are immaterial.   For our purpose, if $\lambda\in{\mathbb N}$ and $\hat I$ is $I^\lambda$, then we can use $\left|s_{\hat I}\right|=\left|s_I\right|^\lambda$.
For a holomorphic map $\psi:\Delta\to{\mathbb C}^n$ with $\psi(0)=0$ and an ideal $I$ of ${\mathcal O}_{{\mathbb C}^n,0}$ with generators $s_{1,I},\cdots,s_{k_I,I}$, by the vanishing order $a_{I,\psi}$ of $I$ on $\psi$ at $0$ we mean the minimum of ${\rm ord}_0\left(s_{j,I}\circ\psi\right)$ for $1\leq j\leq k_I$, where ${\rm ord}_0\left(\cdot\right)$ denotes the vanishing order on ${\mathbb C}$ at the origin.  For an $\ell$-jet $\xi$ of ${\mathbb C}^n$ at the origin which can be represented by $\psi$ we denote $a_{I,\psi}$ also by $a_{I,\xi}$.  (Here the convention is that a $1$-jet is a tangent vector.) If $a_{I,\psi}<\ell$, then $a_{I,\xi}=a_{I,\varphi}$ for any holomorphic map $\varphi:\Delta\to{\mathbb C}^n$ with $\varphi(0)=0$ which represents the $\ell$-jet $\xi$.

\medbreak Note that for our purpose we could also use alternatively the concept of the {\it normalized vanishing order} of $I$ on $\psi$ at $0$ (instead of the vanishing order $a_{I,\psi}$) by defining the normalized vanishing order of $I$ on $\psi$ at $0$ as the minimum of $$\frac{{\rm ord}_0\left(s_{j,I}\circ\psi\right)}{{\rm ord}_0\psi}$$ for $1\leq j\leq k_I$, where ${\rm ord}_0\psi$ is the minimum of the vanishing orders of the $n$ components of $\psi$ on ${\mathbb C}$ at the origin.

\medbreak Since all the main arguments in this proof occur already in the proof of the special case where $N=n=3$, for notational simplicity we will only present the proof of this special case.  The general case is completely analogous but with much more complicated notations.  We break down the proof into the following five steps.

\bigbreak\noindent{\it Step One.}  Let $G_1,G_2$ be holomorphic function germs on ${\mathbb C}^3$ at the origin vanishing at the origin such that the divisor $Z_1$ of $G_1$ is irreducible and of multiplicity $1$.  Assume that $dG_1\wedge dG_2$ is not identically zero.  Then there exists some positive constant $C$ such that
$$
\sum_{k_1,k_2=1}^3\left|G_2\left(dG_1\wedge dz_{k_1}\wedge dz_{k_2}\right)\right|^2\leq
C\sum_{j=1}^3\left|dG_1\wedge dG_2\wedge dz_j\right|^2
$$
on $Z_1=\left\{G_1=0\right\}$.

\medbreak\noindent Step One is verified by\begin{itemize}\item[(i)] taking any holomorphic curve $\varphi:\Delta\to G_1$ with $\varphi(0)=0$ and the image of $\varphi\left(\Delta\right)$ not contained in the zero-set of $G_2$, \item[(ii)]  using the fact that the vanishing order at the origin of the pullback $G_2\circ\varphi$ on $\Delta$ is no more than the minimum of the vanishing orders of its first-order partial derivatives at the origin, and \item[(iii)] observing that at a regular point of $Z_1$, where $z_{k_1}, z_{k_2}$ are used as local coordinates, the component of the gradient of the restriction of $G_2$ to $Z_1$ for the coordinate $z_{k_1}$ is equal to the quotient of $dG_1\wedge dG_2\wedge dz_{k_2}$ by
$dG_1\wedge dz_{k_1}\wedge dz_{k_2}$ as one can easily see by using the chain rule and the implicit differentiation for functions defined on $Z_1=\left\{G_1=0\right\}$.\end{itemize}

\bigbreak\noindent{\it Step Two.} Let $I$ and $J$ be ideals in ${\mathcal O}_{{\mathbb C}^3,0}$ contained in the maximum ideal ${\mathfrak m}_{{\mathbb C}^3,0}$ of ${\mathcal O}_{{\mathbb C}^3,0}$ such that $I$ contains $\left({\mathfrak m}_{{\mathbb C}^n,0}\right)^q$ for some positive integer $q$.  If $\left|s_I\right|$ is not dominated by a constant times $\left|s_J\right|$, then there exists some $(q+2)$-jet $\xi$ of ${\mathbb C}^3$ at the origin which is represented by some holomorphic map $\psi:\Delta\to{\mathbb C}^3$ with $\psi(0)=0$ such that $a_{I,\psi}\leq q$ and $a_{I,\psi}<a_{J,\psi}$.

\bigbreak\noindent{\it Step Three.}  Let $A$ be the ideal generated by elements $F_1, F_2, F_3$ of ${\mathfrak m}_{{\mathbb C}^3,0}$ such that $A$ contains $\left({\mathfrak m}_{{\mathbb C}^3,0}\right)^q$ for some positive integer $q$ .  Let $p\in{\mathbb N}$.  Then there exists a positive integer $q_1$ depending only on $q$ and there exists a positive integer $m$ depending on $q$ and $p$ with the following property.  For any $p$-jet $\xi$ of ${\mathbb C}^3$ at the origin, let $P_1\left(F_1,F_2,F_3\right)$ be a generic homogeneous polynomial of degree $m$ in $F_1,F_2,F_3$ whose divisor $V$ contains a holomorphic curve representing $\xi$ and let $\varphi:\Delta\to{\mathbb C}^3$ be a holomorphic curve germ with $\varphi(0)=0$ whose image is a generic curve germ in $V$ which represents $\xi$.  Then the minimum vanishing order of $\frac{\partial}{\partial z_\ell}P_1\left(F_1,F_2,F_3\right)$ on the holomorphic curve germ $\varphi:\Delta\to{\mathbb C}^3$ at the origin is no more than $(m-1)a_{A,\varphi}+q_1$ for $1\leq\ell\leq 3$.

\medbreak\noindent Note that $V$, as the divisor of a generic homogeneous polynomial of degree $m$ in $F_1,F_2,F_3$ which contains a holomorphic curve representing $\xi$, is a reduced and irreducible hypersurface germ in ${\mathbb C}^3$ at the origin.  Moreover, the image of $\varphi:\Delta\to{\mathbb C}^3$, as a generic
curve germ in $V$ which represents $\xi$, is contained in $\left\{0\right\}\cup\,{\rm Reg}(V)$, where ${\rm Reg}(V)$ is the regular part of $V$.

\bigbreak\noindent{\it Step Four.}  Let $F_1, F_2, F_3$ be from Step Three.
Let $J$ be the ideal generated by $$\frac{\partial\left(F_{j_1},F_{j_2}\right)}{\partial\left(z_{k_1},z_{k_2}\right)}$$
for $1\leq j_1<j_2\leq 3$ and $1\leq k_1<k_2\leq 3$.   Let $\lambda\in{\mathbb N}$ and $I=\left({\mathfrak m}_{{\mathbb C}^3,0}\right)^\lambda$.  Let $p=\lambda+2$.  Assume that $\left|s_I\right|$ is not dominated by any positive constant times $\left|s_J\right|$.  By Step Two, there exists a $p$-jet $\xi$ of ${\mathbb C}^3$ at origin such that, for any holomorphic map $\varphi:\Delta\to{\mathbb C}^3$ whose $p$-jet at the origin is equal to $\xi$,
the vanishing order $a_{J,\varphi}$ of $J$ on $\varphi$ at $0$ is greater than the vanishing order $a_{I,\varphi}$ of $I$ on $\varphi$ at the origin.   By Step Three we have positive integers $q_1, m$ (with $q_1$ depending only on $p$ and with $m$ depending only on $p$ and $q$) and we have a polynomial $P_1\left(F_1,F_2,F_3\right)$ homogeneous of degree $m$ in $F_1,F_2,F_3$ and a holomorphic curve germ $\varphi:\Delta\to{\mathbb C}^3$ at the origin such that
\begin{itemize}\item[(i)] the divisor of $P_1\left(F_1,F_2,F_3\right)$ is a reduced and irreducible hypersurface germ of ${\mathbb C}^3$ at the origin,\item[(ii)] the image of the holomorphic curve germ $\varphi:\Delta\to{\mathbb C}^3$ is contained in the divisor of
$P_1\left(F_1,F_2,F_3\right)$, \item[(iii)] the holomorphic curve germ $\varphi:\Delta\to{\mathbb C}^3$ represents the $p$-jet $\xi$ of ${\mathbb C}^3$ at the origin, \item[(vi)] the minimum vanishing order of $\frac{\partial}{\partial z_\ell}P_1\left(F_1,F_2,F_3\right)$ on the holomorphic curve germ $\varphi:\Delta\to{\mathbb C}^3$ is no more than $(m-1)a_{A,\varphi}+q_1$ for $1\leq\ell\leq 3$.
\end{itemize}
For any polynomial $P_2\left(F_1,F_2,F_3\right)$ of degree $m$ in $F_1, F_2$ and for $1\leq k_1<k_2\leq 3$, it follows from
$$
\frac{\partial P_j}{\partial z_k}=\sum_{\ell=1}^3\frac{\partial P_j}{\partial F_\ell}\frac{\partial F_\ell}{\partial z_k}
$$ by the chain rule that
$$\frac{\partial\left(P_1,P_2\right)}{\partial\left(z_{k_1},z_{k_2}\right)}
=\sum_{\ell_1,\ell_2=1}^3\frac{\partial P_1}{\partial F_{\ell_1}}\frac{\partial P_2}{\partial F_{\ell_2}}\frac{\partial\left(F_{\ell_1},F_{\ell_2}\right)}{\partial\left(z_{k_1},z_{k_2}\right)}.
$$
The vanishing order of
$$\frac{\partial\left(P_1,P_2\right)}{\partial\left(z_{k_1},z_{k_2}\right)}
$$
on the holomorphic curve germ $\varphi:\Delta\to{\mathbb C}^3$ is at least $2(m-1)a_{A,\varphi}+a_{J,\varphi}$.

\medbreak Applying Step One to the case of $G_1=P_1\left(F_1,F_2,F_3\right)$ and $G_2=P_2\left(F_1,F_2,F_3\right)$ with $P_2\left(F_1,F_2,F_3\right)$ being any generic polynomial homogeneous of degree $m$ in $F_1,F_2,F_3$, we get
$$
\sum_{k_1,k_2=1}^3\left|P_2\left(dP_1\wedge dz_{k_1}\wedge dz_{k_2}\right)\right|^2\leq
C\sum_{j=1}^3\left|dP_1\wedge dP_2\wedge dz_j\right|^2
$$
on $P_1=0$, where $C$ is a positive constant.  We restrict this inequality to the curve $\varphi$ and conclude that
$$ma_{A,\varphi}+(m-1)a_{A,\varphi}+q_1\geq 2(m-1)a_{A,\varphi}+a_{J,\varphi}.$$  By the choice of $P_1\left(F_1,F_2,F_3\right)$ and the holomorphic curve germ $\varphi:\Delta\to{\mathbb C}^3$, we have $a_{J,\varphi}>\lambda$.  Thus
$$ma_{A,\varphi}+(m-1)a_{A,\varphi}+q_1\geq 2(m-1)a_{A,\varphi}+\lambda$$
and we conclude that $\lambda\leq a_{A,\varphi}+q_1\leq q+q_1$, because
$a_{A,\varphi}\leq q$ from the fact that $A$ contains $\left({\mathfrak m}_{{\mathbb C}^3,0}\right)^q$.

\bigbreak\noindent{\it Step Five.} By setting $\lambda=q+q_1+1$, we conclude from Step Four that $\left|s_I\right|$ is dominated by a constant times $\left|s_J\right|$.  As in the last part of the proof of Proposition (A.2), by Skoda's theorem [Sk72, Th.1, pp.555-556] it follows from the local integrability of the quotient
$$
\frac{\left|s_I\right|^{2(n+2)}}{\left|s_J\right|^{2(n+2)}}
$$
on ${\mathbb C}^3$ at the origin that $\left({\mathfrak m}_{{\mathbb C}^3,0}\right)^{\left(q+q_1+1\right)\left(n+2\right)}$ is contained in the ideal $J$ generated by $$\frac{\partial\left(F_{j_1},F_{j_2}\right)}{\partial\left(z_{k_1},z_{k_2}\right)}$$
for $1\leq j_1<j_2\leq 3$ and $1\leq k_1<k_2\leq 3$.  This finishes the proof. 

\medbreak We would like to remark that the main point of this proof is to apply the argument for gradients given in Proposition (A.2) for ${\mathbb C}^n$ to the divisor of $P_1\left(F_1,F_2,F_3\right)$ in ${\mathbb C}^3$ instead of to ${\mathbb C}^n$.  Q.E.D.

\bigbreak\noindent(A.4) {\it Proposition.}  Let $h_1,\cdots,h_n$ be
holomorphic function germs on ${\mathbb C}^n$ at the origin so that
the origin is their only common zero. Let $dh_1\wedge\cdots\wedge
dh_n=J\left(dz_1\wedge\cdots\wedge dz_n\right)$.  Then $J$ does not
belong to the ideal generated by $h_1,\cdots,h_n$.

\medbreak\noindent{\it Proof.} Suppose the contrary.  Then there
exist holomorphic function germs $f_1,\cdots,f_n$ on ${\mathbb C}^n$
at the origin such that $J=\sum_{j=1}^n f_j h_j$.  We let
$\omega_j=f_j\left(dz_1\wedge\cdots\wedge dz_n\right)$ for $1\leq
j\leq n$ so that $$dh_1\wedge\cdots\wedge dh_n=\sum_{j=1}^n
h_j\omega_j.\leqno{({\rm A}.4.1)}$$  Since the origin is the only
common zero of $h_1,\cdots,h_n$, we can find connected open
neighborhoods $U$ and $W$ of the origin in ${\mathbb C}^n$ so that
the map $\pi:{\mathbb C}^n\to{\mathbb C}^n$ defined by
$$\left(z_1,\cdots,z_n\right)\mapsto\left(w_1,\cdots,w_n\right)
=\left(h_1\left(z_1,\cdots,z_n\right),\cdots,h_n\left(z_1,\cdots,z_n\right)\right)$$
maps $U$ properly and surjectively onto $W$ and makes $U$ a branched
cover over $W$ of $\lambda$ sheets.  By replacing $U$ and $W$ by
relatively compact open neighborhoods $U^\prime$ and $W^\prime$ of
the origin in $U$ and $W$ respectively, we can assume without loss
of generality that $\int_U\left|\omega_j\right|^2\leq C<\infty$ for
$1\leq j\leq n$. We take the direct image of the equation (A.4.1)
under $\pi$. The left-hand side of the equation (A.4.1) yields
$\lambda\left(dw_1\wedge\cdots\wedge dw_n\right)$, because the map
$\pi$ is defined by $w_j=h_j$ for $1\leq j\leq n$.  Let $\theta_j$
be the direct image of $\omega_j$ under $\pi$ for $1\leq j\leq n$.
Let $Z$ be the branching locus of $\pi$ in $W$.  For any simply
connected open subset $G$ of $W-Z$, $U\cap\pi^{-1}(G)$ is the
disjoint union of $\lambda$ open subsets $H_1,\cdots,H_\lambda$ of
$U$ and $\theta_j(Q)=\sum_{\ell=1}^\lambda\omega_j\left(\tilde
Q_j\right)$, where $U\cap\pi^{-1}\left(Q\right)=\left\{\tilde
Q_1,\cdots,\tilde Q_\lambda\right\}$ with $\tilde Q_j\in H_j$. Now
$$
\int_G\left|\theta_j\right|^2\leq\lambda\sum_{j=1}^\lambda\int_{H_j}\left|\omega_j\right|^2
\leq\lambda C.
$$
Since $W-Z$ can be covered by a finite number of simply connected
open subsets, it follows that
$$
\int_G\left|\theta_j\right|^2<\infty\quad{\rm for\ }1\leq j\leq n.
$$
Thus $\theta_j$ is a holomorphic $n$-form on $G$ and
$$
\lambda\,dz_1\wedge\cdots\wedge dz_n=\sum_{j=1}^n z_j\,\theta_j
$$
on $G$, which gives a contradiction, because the left-hand side does
not vanish at the origin whereas the right-hand side does. Q.E.D.

\bigbreak\noindent(A.5) {\it Remark.}  Proposition (A.4) uses only
the direct images of top-degree holomorphic forms and actually does
not use Skoda's theorem  [Sk72, Th.1, pp.555-556].  The significance
of Proposition (A.4) is that the coefficient $J$ in
$$dh_1\wedge\cdots\wedge dh_n=J\left(dz_1\wedge\cdots\wedge
dz_n\right)$$ cannot be contained in $\left({\mathfrak m}_{{\mathbb
C}^n,0}\right)^p$ if
$$\left({\mathfrak m}_{{\mathbb
C}^n,0}\right)^p\subset\sum_{j=1}^n{\mathcal O}_{{\mathbb
C}^n,0}h_j$$ so that the vanishing order of $J$ at $0$ is no more
than $p$.

\bigbreak\noindent(A.5) {\it Example to Show the Sharpness of
the Exponent in Skoda's Theorem.} The exponent in the denominator of the assumption in
Skoda's theorem [Sk72, Th.1, pp.555-556] plays a r\^ole in effective
bounds.  As stated in Skoda's theorem [Sk72, Th.1, pp.555-556] it is
sharp and cannot be lowered even in the case of Riemann surfaces. Let
$X$ be the Riemann sphere ${\mathbb P}_1$. Consider the hyperplane
section line bundle $H_{{\mathbb P}_1}$. Take two holomorphic
sections $g_1,\,g_2$ of $H_{{\mathbb P}_1}$ without common zeroes.
Take the holomorphic section $f$ of $2\,H_{{\mathbb
P}_1}+K_{{\mathbb P}_1}$ over ${\mathbb P}_1$ which corresponds to a
constant function on ${\mathbb P}_1$ via the isomorphism between
$K_{{\mathbb P}_1}$ and $-2\,H_{{\mathbb P}_1}$. If the exponent
used in the denominator of the assumption in Skoda's theorem [Sk72,
Th.1, pp.555-556] can be lowered so that $\alpha=1$, then $p=2$ and
$n=1$ and $q=\min\left(n,p-q\right)=1$ and $\alpha q+1=2$ and the
assumption
$$
\int_{{\mathbb
P}_1}\frac{\left|f\right|^2}{\left(\left|g_1\right|^2+\left|g_2\right|^2\right)^{\alpha
q+1}}<\infty
$$
is satisfied because $g_1,\,g_2$ have no common zeroes.  Note that
when $\alpha>1$, the integrand of the above inequality makes no
sense unless $$f\in\Gamma\left({\mathbb P}_1,\,m\,H_{{\mathbb
P}_1}+K_{{\mathbb P}_1}\right)$$ for some $m>2$.  If Skoda's theorem
[Sk72, Th.1, pp.555-556] holds with the lower exponent in the
denominator in its assumption, then we can write $f=h_1g_1+h_2g_2$
with
$$
h_1,\,h_2\in\Gamma\left({\mathbb P}_1,\,H_{{\mathbb
P}_1}+K_{{\mathbb P}_1}\right)
$$
which is impossible, because
$$
\Gamma\left({\mathbb P}_1,\,H_{{\mathbb P}_1}+K_{{\mathbb
P}_1}\right)=0
$$
from the isomorphism between $K_{{\mathbb P}_1}$ and
$-2\,H_{{\mathbb P}_1}$.

\bigbreak\centerline{ }

\bigbreak\noindent{\bf References}

\medbreak\noindent[Ca84] D. Catlin, Boundary invariants of
pseudoconvex domains. {\it Ann. of Math.} {\bf 120} (1984),
529--586.

\medbreak\noindent[DA79] J. P. D'Angelo, Finite type conditions for
real hypersurfaces. {\it J. Diff. Geom.} {\bf 14} (1979), 59--66.

\medbreak\noindent[DF77] K. Diederich and J. E. Fornaess,
Pseudoconvex domains: bounded strictly plurisubharmonic exhaustion
functions. {\it Invent. Math.} {\bf 39} (1977), 129--141.

\medbreak\noindent[DF78] K. Diederich and J. E. Fornaess,
Pseudoconvex domains with real-analytic boundary. {\it Ann. of
Math.} {\bf 107} (1978), 371-384.

\medbreak\noindent[GT83] D. Gilbarg and N. Trudinger, {\it Elliptic
partial differential equations of second order}. Second edition.
Grundlehren der Mathematischen Wissenschaften, Volume {\bf 224}.
Springer-Verlag, Berlin, 1983.

\medbreak\noindent[He62] S. Helgason, {\it Differential Geometry and
Symmetric Spaces}, Academic Press 1962.

\medbreak\noindent[Ko79] J.~J.~Kohn, Subellipticity of the $\bar
\partial $-Neumann problem on pseudo-convex domains: sufficient
conditions. {\it Acta Math.} {\bf 142} (1979), 79--122.

\medbreak\noindent[Ni07] Andreea C. Nicoara, Equivalence of types and Catlin boundary systems, arXiv:0711.0429 (November 2007).

\medbreak\noindent [Sk72] H.  Skoda, Application des techniques
$L^2$ \`a la th\'eorie des id\'eaux d'une alg\`ebre de fonctions
holomorphes avec poids, {\it Ann. Sci. Ec. Norm. Sup.} {\bf 5}
(1972), 548-580.

\medbreak\noindent[To72] J.-C. Tougeron, {\it Id\'eaux de fonctions
diff\'erentiables}. Ergebnisse der Mathematik und ihrer
Grenzgebiete, Volume {\bf 71}. Springer-Verlag, Berlin-New York,
1972.

\medbreak\noindent[ZS60] O. Zariski and P. Samuel, {\it Commutative
algebra}. Vol. II. The University Series in Higher Mathematics. D.
Van Nostrand Co., Inc., Princeton, N. J.-Toronto-London-New York
1960.

\bigbreak\noindent{\it Author's mailing address:} Department of
Mathematics, Harvard University, Cambridge, MA 02138, U.S.A.

\medbreak\noindent{\it Author's e-mail address:}
siu@math.harvard.edu

\end{document}